\documentclass[12pt]{article}
\usepackage{amsmath,amsthm}
\usepackage{amssymb,mathrsfs,dsfont}
\usepackage{xspace}
\usepackage{enumerate}

\usepackage[a4paper,
centering, 
left=94.27pt,right=94.27pt,
textheight=649.5pt,
headheight=3pt,
headheight=15pt
]{geometry}      

\usepackage[
unicode]{hyperref}
\hypersetup{ 
    bookmarksnumbered={true},
    bookmarksopen={true},
    unicode={true},
    pdftitle =Time inhomogeneous Generalized Mehler Semigroups,
    pdfauthor ={Ouyang,S.-X. and R\"ockner, M.},
    pdfsubject={Inhomogeneous generalized Mehler semigroup},
    pdfkeywords ={Inhomogeneous generalized Mehler semigroup,  additive process, Ornstein-Uhlenbeck process,  evolution system of measures, Harnack inequality, null controllability, strong Feller property},
    pdfdisplaydoctitle=true, 
    colorlinks = true,
    linkcolor = black,
    anchorcolor = black,
    citecolor = black,
    filecolor = black,
    pagecolor = black,
    urlcolor = black
}

\newtheorem{theorem}{Theorem}[section]
\newtheorem{prop}[theorem]{Proposition}
\newtheorem{lemma}[theorem]{Lemma}
\newtheorem{cora}[theorem]{Corollary}

\theoremstyle{definition}
\newtheorem{definition}[theorem]{Definition}
\newtheorem{assumption}[theorem]{Assumption}
\newtheorem{example}[theorem]{Example}
\newtheorem{remark}[theorem]{Remark}

\def\<{\langle} 
\def\>{\rangle}

\renewcommand{\H}{\ensuremath{\mathds{H}}\xspace}
\DeclareMathOperator{\e}{e}\def\E{\mathds E}
\DeclareMathOperator{\tr}{tr}
\DeclareMathOperator{\var}{var}
\DeclareMathOperator{\Span}{Span}
\DeclareMathOperator{\Prob}{Pr}

\def\m{{\textnormal {m}}}
\def\D{\mathscr D}
\def\P{\mathds P}
\def\Q{\mathds Q}

\def\BB{B} 
\def\B{\mathscr B}

\def\N{\mathbb N}
\def\F{\mathscr F}
\def\A{\mathscr A}
\def\R{\mathds R}

\DeclareFontFamily{U}{bbold}{}
\DeclareFontShape{U}{bbold}{m}{n}{<-5>bbold5<6>bbold6<7>bbold7<8>bbold8
<9>bbold9<9-11>bbold10<11-14>bbold12<14->bbold12}{}
\DeclareSymbolFont{bbold}{U}{bbold}{m}{n}
\DeclareSymbolFontAlphabet{\mathbbold}{bbold}
\DeclareMathSymbol{\Eins}{\mathord}{bbold}{`1}

\numberwithin{equation}{section}

\makeatletter
\DeclareRobustCommand*{\bfseries}{%
  \not@math@alphabet\bfseries\mathbf
  \fontseries\bfdefault\selectfont
  \boldmath
}
\makeatother

\title{\bf Time inhomogeneous Generalized Mehler Semigroups\footnote{Supported in part by the DFG through SFB--701 and IRTG--1132.
The support of Issac Newton Institute
for Mathematical Sciences in Cambridge is also gratefully acknowledged
where part of this work was done during the special semester on
``Stochastic Partial Differential Equations''.
}}
\author{
 {\bf Shun-Xiang Ouyang\footnote{souyang@math.uni-bielefeld.de (corresponding author).}
  and Michael R\"ockner\footnote{roeckner@math.uni-bielefeld.de.}}
\\[0.4em]
\footnotesize{Department of Mathematics, Bielefeld
University, D-33501 Bielefeld, Germany}}

\begin{document}

\maketitle

\begin{abstract} 
A time inhomogeneous generalized Mehler semigroup  on a real separable 
Hilbert space \H is defined through
\[
p_{s,t}f(x)=\int_\H f(U(t,s)x+y)\,\mu_{t,s}(dy), \quad t\geq s, \ x\in\H
\]
for every bounded measurable function $f$ on $\H$, 
where $(U(t,s))_{t\geq s}$ is an evolution family of bounded operators on $\H$  and $(\mu_{t,s})_{t\geq s}$ is a family of probability measures on $(\H, \B(\H))$ satisfying the time inhomogeneous skew convolution equations 
$$\mu_{t,s}=\mu_{t,r}*\left(\mu_{r,s}\circ U(t,r)^{-1}\right),\quad t\geq r\geq s.$$ 
This kind of semigroup is closely related with 
the ``transition semigroup'' of non-autonomous (possibly non-continuous) Ornstein-Uhlenbeck process driven by some proper additive process. We show the weak continuity,  infinite divisibility, associated ``additive processes", L\'evy-Khintchine type representation, construction and spectral representation of $(\mu_{t,s})_{t\geq s}$. We study the structure, existence and uniqueness of the corresponding evolution systems of measures (=space-time invariant measures) of $(p_{s,t})_{t\geq s}$. We also establish dimension free Harnack inequalities in the sense of Wang (1997, PTRF) for $(p_{s,t})_{t\geq s}$. As applications of the Harnack inequalities, we investigate the strong Feller property and contractivity etc. for $p_{s,t}$. 
Finally we prove a Harnack inequality and show the strong Feller property for the transition semigroup of a semi-linear non-autonomous Ornstein-Uhlenbeck process driven by a Wiener process.
\end{abstract} 


\medskip


\tableofcontents

\section{Introduction}\label{Sec:Intro}
 Let \(\H\) be a real separable Hilbert space with norm and inner product denoted by \(|\cdot|\) and \(\<\cdot, \cdot\>\) respectively. Let $\B(\H)$ be the space of Borel measurable subsets of \H, and $\BB_b(\H)$ the space of all bounded Borel measurable functions on \H.  A generalized Mehler semigroup \((p_t)_{t\geq 0}\) on  \(\H\) is defined by the formula
\begin{equation}\label{GMehlerSemigroup}
  p_t f(x)=\int_\H f(T_t x+y)\,\mu_t(dy), \quad t\geq 0,\ x\in\H,\ f\in \BB_b(\H).
\end{equation}
Here \((T_t)_{t\geq 0}\) is a strongly continuous semigroup on \(\H\) and \((\mu_t)_{t\geq 0}\) is a family of probability measures on  \((\H,\B(\H))\) satisfying the following skew convolution (semigroup) equation
 \begin{equation}\label{SkewConvolutionMeasure}
   \mu_{t+s}=\mu_s *(\mu_t\circ T_s^{-1}) ,\quad s,t\geq 0.	  
\end{equation}
Recall that for any two positive Borel measures $\mu$ and $\nu$ on $\H$, 
the convolution 
$\mu*\nu$  of $\mu$ and $\nu$ is a Borel measure on \((\H,\B(\H))\) 
such that 
\[
         \mu*\nu(B):=\int_\H\int_\H \Eins_B(x+y)\,\mu(dx)\nu(dy) = \int_\H \mu(B-x)\,\nu(dx),\quad B\in \B(\H).
\]
Condition \eqref{SkewConvolutionMeasure} is necessary and sufficient for the semigroup property of $(p_t)_{t\geq 0}$ 
(and the Markov property of the corresponding stochastic process respectively) to hold. That is, 
\[
\eqref{SkewConvolutionMeasure}\  
\textrm{holds}\ 
\textrm{if\ and\ only\ if}\quad 
p_tp_s=p_{t+s}\quad \textrm{on}\ \BB_b(\H),\quad t,s\geq 0.
\]

The semigroup \eqref{GMehlerSemigroup} is a generalization of the classical Mehler formula for the transition semigroup of Ornstein-Uhlenbeck process driven by Wiener process. See for instance \cite{Ouyang09} and references therein for some details on Ornstein-Uhlenbeck process, Mehler formula and other related topics.  
In \cite{BR95,BRS96} the second author and colleagues first studied this generalization for the Gaussian case. Then Fuhrman and  R\"ockner \cite{FR00} studied the non-Gaussian case.  It is clear (cf. \cite{BRS96,FR00}) that under some weak conditions 
there is a one to one correspondence between 
generalized Mehler semigroups and 
 transition semigroups of L\'evy driven Ornstein-Uhlenbeck processes. 
For the introductions to L\'evy driven Ornstein-Uhlenbeck processes, 
mainly from the point of view of stochastic process, 
we refer to \cite{App06, PZ07} etc..

Now generalized Mehler semigroup has been extensively studied. For instance, Schmuland and Sun \cite{SS01} investigated the infinite divisibility of $\mu_t$ ($t\geq 0$) and the continuity of $t\mapsto\log\hat{\mu}_{t}$; Lescot and R\"ockner considered in \cite{LR02} and \cite{LR04} the generator and perturbation of $(p_t)_{t\geq 0}$ respectively; Wang and R\"ockner established some useful functional inequalities and their applications for $p_t$; Neerven \cite{Van00},  Li and his colleagues \cite{DL04,DLSS04} carefully studied the representation of $(\hat{\mu}_t)_{t\geq 0}$ under non-differentiable condition and related topics. It was also noted that, see the recent monograph by Li \cite{Li11}, skew convolution semigroups \eqref{SkewConvolutionMeasure}
for measures $(\mu_t)_{t\geq 0}$ on the space of finite measures on a Lusin topological space
and Ornstein-Uhlenbeck processes are closely related with measure-valued branching processes.  Applebaum \cite{App06, App07a,App07b} also considered some quite interesting problems.

A probability measure $\mu$ on $(\H,\B(\H))$  is said to be operator self-decomposable if 
\begin{equation}\label{invariant_measure:equ:Xu}
\mu=(\mu\circ T_t^{-1}) *\mu_t,\quad t\geq 0
\end{equation} 
for a family of semigroups $(T_t)_{t\geq 0}$ and measures $(\mu_t)_{t\geq 0}$. 
In the setting of \eqref{GMehlerSemigroup}, any solution $\mu$ to the convolution equation \eqref{invariant_measure:equ:Xu} is just an invariant measure for the generalized Mehler semigroup  \eqref{GMehlerSemigroup}. 
The study of operator self-decomposable distributions has been done in many papers, see for example,  \cite{App06} and references therein. 

Generalized Mehler semigroups and the related skew convolution equations for measures have so many connections with various topics such that it is impossible to trace all literatures here. For example, we noted that the measure equation \eqref{SkewConvolutionMeasure} has also been studied under the name `` (twisted) cocycle equation" (cf. the note by Jurek \cite{Jur04}). 

Recently much work, for instance \cite{DaPL07, DaPR08, GL08,Kna09, Woo09}, has been devoted to the study of non-autonomous Ornstein-Uhlenbeck processes. In these papers,  linear stochastic partial differential equations with time-dependent drifts were studied. The noise in the equations is modeled by 
a stationary process, e.g. a Wiener process or L\'evy process. To get a full inhomogeneous Ornstein-Uhlenbeck process, it is natural to consider a more general noise modeled by (time) non-stationary processes such as additive processes. 
To be precise, let us first describe our framework in detail.  

Let \((A(t), \D(A(t)) )_{t\in\R}\) be a family of linear operators on \(\H\) with dense domains. Suppose that
the non-autonomous Cauchy problem 
\[
\left\{
\begin{aligned}
dx_t&=A(t)x_tdt,\quad t\geq s,\\
x_s&=x
\end{aligned}
\right.
\]
is well posed (see \cite{Paz83}).
That is, there exists an evolution family of bounded operators \((U(t,s))_{t\geq s}\) on \(\H\) such that \(x(t)=U(t,s)x\) for $x\in \D(A(s))$ is a classical unique solution of this Cauchy problem. 

Recall that a family of bounded linear operators 
\((U(t,s))_{t\geq s}\) on \(\H\) is said to be a (strongly continuous) evolution family if 
\begin{enumerate}
 \item For every $s\in \R$, \(U(s,s)=I\) and 
       for all $t\geq r\geq s$, \[U(t,r)U(r,s)=U(t,s).\]
 \item For every \(x\in \H\), the map \((t,s)\mapsto U(t,s)x\) is strongly continuous on 
	\(\{(t,s)\in\R^2\colon t\geq s\}\).
\end{enumerate}
An evolution family is also named evolution system, propagator etc.. 
For more details we refer e.g. to \cite{Paz83,EN00,DK74}.

Let $(Z_t)_{t\geq s}$ be an additive processes  in \H, i.e. an $\H$-valued stochastic continuous stochastic process with independent increments. Consider the following stochastic differential equation 
\begin{equation}\label{NonAutoLangivinEqu}
 \left\{
    \begin{aligned}
	dX_t&=A(t)X_tdt+dZ_t,\\
    	X_s&=x.
    \end{aligned} 
  \right.
\end{equation}
We call the following process 
\begin{equation}\label{NonAutoOUprocess}
 X(t,s,x)=U(t,s)x+ \int_s^t U(t,r)\,dZ_r,\quad 
 t\geq s,\ x\in\H
\end{equation}
a mild solution of \eqref{NonAutoLangivinEqu}
if the stochastic convolution integral in \eqref{NonAutoOUprocess} is well-defined for a proper additive process $(Z_t)_{t\in\R}$ (see \cite{Det83,Sat06}). 
For all $t\geq s$, let \({\mu}_{t,s}\) denote the distribution of the convolution integral \(\int_s^t U(t,r)\,dZ_r\). 
Then the transition semigroup of $X(t,s,x)$ is given by 
\begin{equation}\label{2GMehlerSemigroup0}
 P_{s,t}f(x)=\E f(X(t,s,x))=\int_\H f\left(U(t,s)x+y \right)
		          {\mu}_{t,s}(dy),\quad x\in \H, \ f\in B_b(\H)
\end{equation}
for all \(t\geq s\).

The aim of the present paper is to adopt the axiomatic approach as in \cite{BRS96} to study this non-autonomous process through its transition semigroup \eqref{2GMehlerSemigroup0}. 
That is, we shall consider the abstract form of  \eqref{2GMehlerSemigroup0} which is the time inhomogeneous version of the generalized Mehler semigroup \eqref{GMehlerSemigroup}: 
\begin{equation}\label{2GMehlerSemigroup0-p}
p_{s,t}f(x)=\int_\H f(U(t,s)x+y)\,\mu_{t,s}{dy}, \quad x\in\H, \ f\in \BB_b(\H).
\end{equation}
Here $(\mu_{t,s})_{t\geq s}$ is a family of probability measures on $(\H,\B(\H))$ satisfying 
 the following  skew convolution equation for measures with two parameters
\begin{equation}\label{Skew:Meas:Equ:0}
 \mu_{t,s}=\mu_{t,r}*\left(\mu_{r,s}\circ U(t,r)^{-1}\right),\quad
 s\leq r\leq t.
\end{equation}  
It is clear that this equation is a time inhomogeneous analog of \eqref{SkewConvolutionMeasure}.

We remark that,
analogous to the semigroup  
\eqref{2GMehlerSemigroup0-p} and  the measure equation \eqref{Skew:Meas:Equ:0} respectively,  
 ``Mehler hemigroup" and
``hemigroup of probability measures" which are 
 are studied in, 
for examples \cite{Haz05, BKH09} and \cite{Bor90,HP04} etc., under some special or different situations, particularly on locally compact groups.

The present paper is organized as follows. 

In Section \ref{Sec:2GMehlerSemigroup} we shall introduce the time inhomogeneous transition function \((p_{s,t})_{t\geq s}\) 
\eqref{2GMehlerSemigroup0-p} (rewritten as \eqref{2GMS})
of generalized Mehler type. 
In Proposition \ref{Prop:SemigroupChar} we state that the flow property for $(p_{s,t})_{t\geq s}$, 
 i.e. the Chapman-Kolmogorov equations holds: 
 for all \(s\leq r\leq t\), 
$p_{s,r}p_{r,t}=p_{s,t}$ if and only if \eqref{Skew:Meas:Equ:0} (rewritten as \eqref{Equ:SkewMeasure}) holds for the measures $(\mu_{t,s})_{t\geq s}$.

 We can put the time inhomogeneous generalized Mehler semigroup \eqref{2GMS} in a more general framework, i.e. the inhomogeneous skew convolution semigroups.  We characterize the Markov property of this more general skew convolution semigroups in Proposition \ref{Prop:SemigroupProperty:SCGg} which implies Proposition \ref{Prop:SemigroupChar} immediately. 

In Section \ref{Sec:On_the_Equation} we study the skew convolution equation \eqref{Skew:Meas:Equ:0} for measures in several subsections. We first warm up with some preliminaries and motivations in Subsections \ref{subsec:Pre} and \ref{subsec:Mot}. In the remaining Subsections, we always assume Assumption \ref{Ass:ContinuousMeasure'} which requires that the measure $(\mu_{t,s})_{t\geq s}$ is weakly continuous on the diagonal. 
In Subsection  \ref{subsec:cont}, among some other continuity results, in particular we show that for fixed $t\in\R$, $\mu_{t,s}$ is weakly continuous 
in $s$ for all $s\leq t$. We also investigate the continuity of $p_{s,t}$.
In Subsection \ref{Subsec:Inf_Div}, we prove that for every $t\geq s$, $\mu_{t,s}$ is infinitely divisible. 
In Subsection \ref{Subsec:additive_process} we show that there exists a natural ``additive process" associated with $(\mu_{t,s})_{t\geq s}$. 
Then the estimate \eqref{Equ:UniformStochContinuous:Mod} appears naturally. This provides us a second proof of the  infinite divisibility of $\mu_{t,s}$. 
Finally we study the structure and representations of the measures $\mu_{t,s}$ in Subsections \ref{sub:sec:LKrepr}, \ref{subsec:ConsRep1} and \ref{subsec:ConsRep2}. 

In Section \ref{Sec:Evolution_Measures}  
we study evolution systems of measures, i.e. space-time invariant measures, for the semigroup $(p_{s,t})_{t\geq s}$. 
Concerning this topic, some work has been done in  \cite{BRS96, FR00, DaPL07, Kna09, Woo09} etc.. 
We first show some basic properties of the evolution systems of measures. Then we give sufficient and necessary conditions for the existence and uniqueness of evolution systems of measures. In particular we emphasize that  Theorem \ref{EMeasure:Levy2GMS}, Theorem \ref{EMeasure:Levy2GMS:Converse} and Corollary \ref{Cora:Evolu-Measure-uniquness:tight} not only 
generalize the results in \cite{Woo09} for finite dimensional L\'evy driven non-autonomous Ornstein-Uhlenbeck processes but also contain some new results even in the finite dimension case. As applications, we show an example using periodic conditions (in time) (cf. \cite{DaPL07, Kna09}).

In Section \ref{Sec:Harnack inequalities} we prove Harnack inequalities for \((p_{s,t})_{t\geq s}\)
 using much simpler arguments than the methods used in the previous papers 
\cite{RW03, Kna09, ORW09, Ouyang09} where Harnack inequalities for generalized Mehler semigroups or Ornstein-Uhlenbeck semigroup driven by L\'evy processes were shown.
The method in \cite{Kna09} and \cite{RW03} relies on taking the derivative of a proper functional; the method in  \cite{ORW09,Ouyang09}  is based on coupling of stochastic processes and Girsanov transformation. 
Our approach in this paper is based on a convolution decomposition of $p_{s,t} $ for each $t\geq s$.  As applications of the Harnack inequality, we prove that null controllability implies the strong Feller property and that for the  Gaussian case,  null controllability, Harnack inequality and strong Feller property are in fact equivalent to each other as in the time homogeneous case. 

In Section \ref{Sec:Semi-linear} we apply Girsanov's theorem to study the existence of martingale solutions of semi-linear non-autonomous Ornstein-Uhlenbeck process driven by a Wiener process
for possibly non-Lipschitz non-linearities. 
For the Lipschitz case we refer to \cite{Ver09}. 
Our approach is an adaption of the standard procedure when the linear part, i.e the operator $A$ does not depend on time 
(see \cite[Chapter 10]{DZ92}).
Our main contribution here is to establish a Harnack inequality and hence show the strong Feller property for the transition semigroup.

In Section \ref{Sec:NullControll}
we append a brief introduction to the control theory of non-autonomous linear control systems and null controllability. This  is closely related to the strong Feller property of the corresponding Ornstein-Uhlenbeck processes. The minimal energy representation also proves useful for the estimates of the constants in the Harnack inequalities.

\section{Time inhomogeneous generalized Mehler\\ semigroups and skew convolution semigroups}\label{Sec:2GMehlerSemigroup}
Let \((U(t,s))_{t\geq s}\) be an evolution family of operators on \(\H\) and \((\mu_{t,s})_{t\geq s}\) a family of probability measures on \((\H, \B(\H))\). 
For every \(f\in\BB_b(\H)\) and \(t\geq s\), define
\begin{equation}\label{2GMS}
\begin{aligned} 
 p_{s,t}f(x)
=&\int_\H f(U(t,s)x+y)\,\mu_{t,s}(dy),\quad x\in \H.
\end{aligned}
\end{equation}
In terms of convolution, we may write 
\begin{equation}\label{p_st:convolution_form}
p_{s,t}f(x)=(\delta_{U(t,s)x}*\mu_{t,s}) f, \quad x\in\H.
\end{equation}
Here for any $x\in\H$, $\delta_x$ is the Dirac measure concentrating at the single point $x$. 

We are interested in the case when the flowing Chapman-Kolmogorov equations holds for $(p_{s,t})_{t\geq s}$ in \eqref{2GMS}: 

\begin{prop}\label{Prop:SemigroupChar}
For all \(s\leq r\leq t\), 
  \begin{equation}\label{SemigroupProperty}
	p_{s,r}p_{r,t}=p_{s,t}   \quad \textnormal{(``Chapman-Kolmogorv\ equations")}
  \end{equation}
holds on \(\BB_b(\H)\) if and only if for all
\(s\leq r\leq t\), 
\begin{equation}\label{Equ:SkewMeasure}
 \mu_{t,s}=\mu_{t,r}*\left(\mu_{r,s}\circ U(t,r)^{-1}\right).
\end{equation}
\end{prop}

The proof is given at the end of this section after the proof of Proposition \ref{Prop:SemigroupProperty:SCGg} which gives a more general result for skew convolution semigroups.

Later on, we shall always assume that \eqref{Equ:SkewMeasure}  or equivalently the following equation holds
\begin{equation}\label{Equ:SkewMeasure:Fourier}
 \hat{\mu}_{t,s}(\xi)=\hat{\mu}_{t,r}(\xi)
		\hat{\mu}_{r,s}(U(t,r)^*\xi),\quad \xi\in \H.
\end{equation}
Here for any linear operator $U$ on $\H$, $U^*$ denotes the adjoint operator of $U$, and for every probability measure \(\mu\) on \((\H,\B(\H))\), $\hat{\mu}$ denotes its Fourier transform (the characteristic functional of $\mu$), i.e., 
\[\hat{\mu}(\xi)=\int_\H \e^{i\<x,\xi\>}\mu(dx),\quad \xi\in  \H.\]
For the study of probability measures in finite and infinite dimensional spaces and their characteristic functional we recommend the monographs \cite{Sat99, VTC87, Lin86, Par67} etc..
In particular, 
the equivalence of \eqref{Equ:SkewMeasure} and \eqref{Equ:SkewMeasure:Fourier} follows from the fact that 
probability measures on Hilbert spaces  are determined by their characteristic functionals (see e.g. \cite[Section IV.2.2, Theorem 2.2, Page 200]{VTC87}). 

By Proposition \ref{Prop:SemigroupChar} we call the family of probability kernels \((p_{s,t})_{t\geq s}\) defined in \eqref{2GMS} with 
$(\mu_{t,s})_{t\geq s}$ satisfying \eqref{Equ:SkewMeasure}
 is called a \emph{time inhomogeneous generalized Mehler semigroup}. 

It is also interesting to look at the connection between the transition functions determined by the transition semigroup and the family of probability measures satisfying \eqref{Equ:SkewMeasure}.

If \eqref{SemigroupProperty} holds for $(p_{s,t})_{t\geq s}$, then we have the following associated  transition function $(p_{s,t}(\cdot,\cdot))_{t\geq s}$  given by
\begin{equation}\label{p_st_translation_fnt_2_semigroup}
	p_{s,t}(x,B)=p_{s,t}\Eins_B(x) = \int_\H \Eins_B(U(t,s)x+y)\,\mu_{t,s}(dy)
\quad x\in\H,\  B\in\B(\H)	
	 .
\end{equation}
The transition function $(p_{s,t}(\cdot,\cdot))_{t\geq s}$ is not only time inhomogeneous but also space inhomogeneous. In fact we have for all $t\geq s$, 
\begin{equation}\label{equ:mu_ts_2_p_st}
	p_{s,t}(x,B)=p_{s,t}(0,B-U(t,s)x)=\mu_{t,s}(B-U(t,s)x),\quad x\in\H,\ B\in\B(\H).
\end{equation}

On the other hand, given an evolution family of operators \((U(t,s))_{t\geq s}\) on \(\H\) and a family of transition function $(p_{s,t}(\cdot,\cdot))_{t\geq s}$, $x\in\H$, $B\in\B(\H)$ or equivalently a semigroup $(p_{s,t})_{t\geq s}$ satisfying $p_{s,t}=p_{s,r}p_{r,t}$ for all $s\leq r\leq t$, we can define a family of probability measures $(\mu_{t,s})_{t\geq s}$ on \(\H,\B(\H)\)  by
\[
	\mu_{t,s}(B):=p_{s,t}(0,B),\quad B\in\B(\H).
\]
If $(p_{s,t}(\cdot,\cdot))_{t\geq s}$ satisfies the  space inhomogeneous condition 
$$p_{s,t}(x,B)=p_{s,t}(0,B-U(t,s)x),\quad x\in\H,\ B\in\B(\H)$$ 
as in \eqref{equ:mu_ts_2_p_st}, then it is easy to verify that $(p_{s,t})_{t\geq s}$
 is the time inhomogeneous generalized Mehler semigroups given in \eqref{2GMS} and $(\mu_{t,s})_{t\geq s}$ satisfies \eqref{Equ:SkewMeasure}.

Taking $s=r=t$ in \eqref{Equ:SkewMeasure} we obtain  
\begin{equation}\label{delta_tt=0}
 \mu_{t,t}=\mu_{t,t}*\mu_{t,t},\quad t\in\R.
\end{equation}
That is, for every $t\in\R$, $\mu_{t,t}$ is an idempotent probability measure on $(\H,\B(\H))$. 
By \cite[Section I.4.3, Proposition 4.7, Page 67, see also Section IV.2.2, Corollary 1, Page 203]{VTC87} (or 
\cite[Section III.3, Theorem 3.1, Page 62]{Par67} and noting that there is no nontrivial compact subgroup in \H), 
the trivial measure 
$\delta_0$ is the only idempotent measure on $(\H, \B(\H))$. 
Therefore, \eqref{delta_tt=0} yields 
\begin{equation}\label{Equ:mu_tt=delta0}
	\mu_{t,t}=\delta_0
\end{equation}
for every $t\in\R$.

As noted by Li et al. (see \cite{Li06} for a survey), a generalized Mehler semigroup is a special case of the so called skew convolution semigroup. In the rest of this section we shall briefly discuss time inhomogeneous skew convolution semigroups which constitute a more general framework than time inhomogeneous generalized Mehler semigroups.
But in the following sections of this paper we shall not work in this general framework. 

Let $(S,+)$ be an abelian semigroup with identity (neutral) element $0$. 
The operation $+\colon S^2\rightarrow S$ is associative, commutative and for every $x\in S$, $x+0=x$. 
Let $(u_{s,t})_{t\geq s}$ be a Borel Markov transition function on $S$
satisfying 
\begin{equation}\label{Delta0}
u_{s,t}(0,\cdot)=\delta_0
\end{equation} and 
\begin{equation}\label{u_tran_function}
  u_{s,t}(x+y,\cdot)=u_{s,t}(x,\cdot)*u_{s,t}(y,\cdot)
\end{equation}
for every $t\geq s$ and $x,y\in S$. 

As we discussed previously, if $S$ is an Hilbert space, then \eqref{Delta0} is a simple consequence of \eqref{u_tran_function}. But in general \eqref{u_tran_function} doesn't imply \eqref{Delta0}. We refer to \cite{HM11} for the discussion of idempotent measures on a locally compact Hausdorff second countable topological semigroup. 

Since $(u_{s,t})_{t\geq s}$ is a family of Markov transition functions, we have the following Chapman-Kolmogorov equations
\begin{equation}\label{Equ:CKequ:onS}
  u_{s,t}f(x)=u_{s,r}(u_{r,t}f)(x),\quad x\in S,\  f\in \BB_b(S)
\end{equation}
 for all $t\geq r\geq  s$.
Here $B_b(S)$ denote the space of all bounded measurable function on $S$. 
Writing the equation above \eqref{Equ:CKequ:onS}  in integral form, we have
\begin{equation}\label{CKEqu}
   \int_{S^2} f(z)u_{r,t}(y,dz)u_{s,r}(x,dy)=\int_S f(z)u_{s,t}(x,dz).
\end{equation} 

Let $\B(S)$ denote the space of Borel $\sigma$-algebra over $S$. 
For every probability measure $\mu$ on $(S,\B(S))$ we associate with $u_{s,t}$ $(t\geq s)$ a new probability measure $\mu u_{s,t}$ by 
\[
\mu u_{s,t}(A)=\int_S u_{s,t}(x,A)\,\mu(dx),\quad A\in\B(S)
\]
for every $t\geq s$. 

It is easy to show the following result. 

\begin{prop}
For any two probability measures $\mu$ and $\nu$ on $(S,\B(S))$, we have 
\[(\mu*\nu) u_{s,t}=(\mu u_{s,t})*(\nu u_{s,t})\]
for all $t\geq s$. 
\end{prop}

\begin{proof}
  For every $f\in\BB_b(S)$, by \eqref{Delta0} we have
  \[
    \begin{aligned}
          (\mu*\nu) u_{s,t} f 
       &= \int_{S^2} f(z) u_{s,t}(x,dz) \, (\mu*\nu)(dx)\\
       & = \int_{S^3} f(z) u_{s,t}(x+y,dz) \, \mu(dx)\nu(dy)\\             
       &= \int_{S^3} f(z) 
              \bigl[u_{s,t}(x,\cdot)*u_{s,t}(y,\cdot)\bigr] (dz) 
       \, \mu(dx)\nu(dy) \\
       &= \int_{S^4} f(z_1+z_2) u_{s,t}(x,dz_1) \mu(dx)\, u_{s,t}(y,dz_2) \nu(dy)\\
       &= \int_{S^2} f(z_1+z_2) (\mu u_{s,t})(dz_1)\, (\nu u_{s,t})(dz_2)\\
       &= \int_S f(z) ((\mu u_{s,t})* (\nu u_{s,t}))(dz)
       = ((\mu u_{s,t})* (\nu u_{s,t}) ) f.
    \end{aligned}
  \]
This completes the proof.   
\end{proof}

Now let $(\mu_{t,s})_{t\geq s}$ be a family of probability measures on $(S,\B(S))$. For all $t\geq s$, define a family of functions $q_{s,t}(\cdot,\cdot)\colon S\times \B(S)\to \R$ by 
\[
	q_{s,t}(x,\cdot)=u_{s,t}(x,\cdot)*\mu_{t,s}(\cdot),\quad x\in S.
\]
 Associated with $q_{s,t}(\cdot,\cdot)$ we can define an operator $q_{s,t}$ on $B(S)$ by 
\[
	q_{s,t}f(x)=\int_S f(y)\,q_{s,t}(x,dy),\quad x\in S,\ f\in B_b(S).
\]

We have the following characterization of the Markov property of $(q_{s,t})_{t\geq s}$.

\begin{prop}\label{Prop:SemigroupProperty:SCGg}
The family of  operators 
 $(q_{s,t})_{t\geq s}$ 
 has the following property
\begin{equation}\label{SemigroupProperty:SCGg}
 q_{s,t}=q_{s,r}q_{r,t},\quad t\geq r\geq s
\end{equation}
if and only if 
\begin{equation}\label{Non_auto_SCSg}
  \mu_{t,s}=\mu_{t,r}*(\mu_{r,s} u_{r,t}),\quad t\geq r\geq s,
\end{equation}
or equivalently
\[
  \hat{\mu}_{t,s}(\xi)=\hat{\mu}_{t,r}(\xi) 
  \widehat{(\mu_{r,s} u_{r,t}) }(\xi),\quad \xi\in S,
        \ t\geq r\geq s.
\]
\end{prop}

\begin{proof}
  For every $f\in\BB_b(S)$, $x\in S$, we have
\begin{equation}\label{Thm:Proof:SC1}
\begin{aligned}
  &q_{s,r}q_{r,t}f(x)\\
=&\int_S q_{r,t}f(y) q_{s,r}(x,dy)\\
=&\int_{S^2} q_{r,t}f(y_1+y_2) u_{s,r}(x,dy_1)\mu_{r,s}(dy_2) \\
=&\int_{S^4} f(z) q_{r,t}(y_1+y_2,dz) u_{s,r}(x,dy_1)\mu_{r,s}(dy_2) \\
=&\int_{S^4} f(z_1+z_2) 
	u_{r,t}(y_1+y_2,dz_1) 
        \mu_{t,r}(dz_2) u_{s,r}(x,dy_1)\mu_{r,s}(dy_2) \\
=&\int_{S^5} f(z_{11}+z_{12}+z_2) 
        u_{r,t}(y_1,dz_{11}) u_{r,t}(y_2,dz_{12}) 
        \mu_{t,r}(dz_2) u_{s,r}(x,dy_1)\mu_{r,s}(dy_2) \\
=& \int_{S^4} f(z_{11}+z_{12}+z_2) 
        u_{s,t}(x,dz_{11}) u_{r,t}(y_2,dz_{12}) 
        \mu_{t,r}(dz_2) \mu_{r,s}(dy_2) \\
=& \int_{S^3} f(z_{11}+z_{12}+z_2) 
        u_{s,t}(x,dz_{11}) (\mu_{r,s} u_{r,t})(dz_{12}) 
        \mu_{t,r}(dz_2) \\
=& \int_{S} f(z) 
        (u_{s,t}(x,\cdot)* (\mu_{r,s} u_{r,t}) 
        *\mu_{t,r})(dz).
\end{aligned}
\end{equation}
Here we have used \eqref{CKEqu} to get the sixth identity in the calculation above.
If \eqref{Non_auto_SCSg} holds, then by \eqref{Thm:Proof:SC1} we obtain
\[q_{s,r}q_{r,t}f(x)=\int_{S} f(z) 
        [u_{s,t}(x,\cdot)* 
        \mu_{t,s}](dz)=q_{s,t}f(x).\]
That is, \eqref{SemigroupProperty:SCGg} holds. 

Conversely, if \eqref{SemigroupProperty:SCGg} holds, then by taking $x=0$ in  \eqref{Thm:Proof:SC1} and using 
\eqref{Delta0}, 
we get 
\[\int_S f(z) [(\mu_{r,s}u_{r,t})*\mu_{t,r}](dz)
=\int_S f(z) \mu_{t,s} (dz)
\]
for every $f\in\BB_b(S)$. This implies \eqref{Non_auto_SCSg}. So the proof is complete. 
\end{proof}

\begin{proof}[Proof of Proposition \ref{Prop:SemigroupChar}]
We shall apply Proposition \ref{Prop:SemigroupProperty:SCGg}. 
Let $S=\H$ be a real separable Hilbert space and 
$$u_{s,t}(x,\cdot)=\delta_{U(t,s)x}$$
 for every $t\geq s$ and $x\in \H$. 
 According to \eqref{p_st:convolution_form}, 
\((q_{s,t})_{t\geq s}\) coincides with the inhomogeneous generalized Mehler semigroup 
$(p_{s,t})_{t\geq s}$ defined in \eqref{2GMS} and the equivalence of 
\eqref{SemigroupProperty:SCGg} and \eqref{Non_auto_SCSg} in Proposition \ref{Prop:SemigroupProperty:SCGg}
is exactly the equivalence of \eqref{SemigroupProperty} and \eqref{Equ:SkewMeasure} in Proposition \ref{Prop:SemigroupChar}. The latter is thus proved. 
\end{proof}

\begin{example}\label{skew_sgExa2}
 Let $S=M(E)$ be the space of all finite Borel measures on a Lusin topological space $E$. Let $(u_{s,t})_{t\geq s}$ be the transition semigroup of some measure-valued branching process and $(\mu_{t,s})_{t\geq s}$ be a family of probability measures on $M(E)$ satisfying \eqref{Non_auto_SCSg}. Then $(q_{s,t})_{t\geq s}$ is called an immigration process in \cite{Li02}.
\end{example}

\section{On the equation $ \mu_{t,s}=\mu_{t,r}*\left(\mu_{r,s}\circ U(t,r)^{-1}\right)$}\label{Sec:On_the_Equation}
As we have seen in the previous section, the time inhomogeneous skew convolution equation \eqref{Equ:SkewMeasure} for probability measures is  essential for the study of the semigroup $(p_{s,t})_{t\geq s}$. Hence in this section, we concentrate on the equation  \eqref{Equ:SkewMeasure}. We shall study the weak continuity, infinite divisibility, constructions, associated additive processes and representations of $(\mu_{t,s})_{t\geq s}$.

\subsection{Preliminaries}\label{subsec:Pre}
In order to fix some notations and recall some results, we put here some preliminaries which can be found for example in the monographs \cite{Bil99, VTC87, Lin86, Par67,Sko91,JS87,Sat99} and paper \cite{Rit88} etc..

\textbf{Convergence of probability measures} 
We say that a sequence of probability measures $(\mu_n)_{n\geq 1}$ converges weakly to a probability measure $\mu$, written as
\[
      \mu_n \Rightarrow \mu\quad \textrm{as}\ n\to \infty
\]
if, for every $f\in C_b(\H)$
\[
   \lim_{n\to\infty} \int_\H f(x)\,\mu_n(dx)=\int_\H f(x)\,\mu(dx). 
\]
Here $C_b(\H)$ denotes the space of bounded continuous functions on \H. Sometimes we also write shortly 
\[
	\lim_{n\to \infty} \mu_n =\mu.
\]
Similarly we may define weak convergence of probability measures $\mu_{t} \Rightarrow \mu_{t_0}$ as $t\to t_0$, and 
$\mu_{t,s}\Rightarrow \mu_{t_0,s_0}$ as $t\to t_0$ and $s\to s_0$ for 
$t\geq s$ and $t_0\geq s_0$.

We say that a sequence of $\H$-valued random variables $(X_n)_{n\geq 1}$ converges 
stochastically, or converges in probability, to a $\H$-valued random variable $X$, 
written as 
\[
      X_n \stackrel{\Prob}{\longrightarrow} X,
\]
if, for each $\varepsilon >0$, 
\[
    \lim_{n\to\infty}\P(|X_n-X|>\varepsilon)=0.
\]
It is equivalent to say that for each $\varepsilon >0$ and $\eta>0$, there exists a constant $N>0$ such that for each $n>N$, we have
\[
     \P(|X_n-X|>\varepsilon)<\eta.
\]
Assume that $I$ is an uncountable index set. Then the convergence in probability of $\{X_t\}_{t\in I}$ is defined similarly. 


Let $\mu_n$ and $\mu$ denote the distributions of $X_n$ and $X$ respectively. Then it is well known that as $n\to\infty$, $X_n \stackrel{\Prob}{\longrightarrow}  X$ implies 
$\mu_n\Rightarrow\mu$ (in other words, $X_n$ converges to $X$ in distribution).  On the other hand, if in particular, $X=x\in\H$ is deterministic, then $\mu_n\Rightarrow \delta_x$ implies $X_n \stackrel{\Prob}{\longrightarrow}  X$. Therefore, we have 
$X_n \stackrel{\Prob}{\longrightarrow}  x$ if and only if
$\mu_n\Rightarrow \delta_x$. 

This simple observation above is helpful for the understanding of the relation between  \eqref{additive_process:stoch_continuous} and \eqref{additive_process:stoch_continuous''}, as well as for the understanding of Assumption \ref{Ass:ContinuousMeasure'}, Lemma \ref{lmm:Ass:ContinuousMeasure'} and Proposition \ref{Ass:ContinuousMeasure}.

\textbf{Infinitely divisible probability measures}
A probability measure $\mu$ on $(\H,\B(\H))$ is said to be infinitely divisible if for any $n\in\N$, there exists a probability measure $\mu_n$ on $(\H,\B(\H))$ such that 
\[
\mu=\mu_n^{*n}:=\underbrace{\mu_n*\mu_n*\cdots*\mu_n}_{n\ \textrm{times}}. 
\]
In this case we call $\mu$ an infinitely divisible probability measure (or distribution). 

It is well known that 
a finite Borel measure $\mu$ is infinitely divisible if and only if 
its characteristic functional $\mu$ has a characteristic exponent 
\[\psi(\xi):=-\log\hat{\mu}(\xi),\quad \xi\in\H\]
satisfying the following conditions
\begin{enumerate}
 \item $\psi(0)=0$.
 \item $\psi$ is negative definite, i.e.\ 
 for all  $n$-tuple of elements $(\xi_1,\xi_2,\cdots,\xi_n)$ from $\H$, $n=1,2,\cdots$,
 the $n\times n$ matrix 
\[
   (\psi(\xi_i)+\overline{\psi(\xi_j)}-\psi(\xi_i-\xi_j))_{i,j}
\]
is positive hermitian.
 \item $\psi$ is Sazonov continuous, i.e.\ it is continuous with respect to the locally convex topology on $\H$ defined by the seminorms $y\mapsto \|  Ay\|$ for all $y\in \H$, where $A$ runs over all Hilbert-Schmidt operators on $\H$. 
\end{enumerate}

Recall that an operator $A$ on \H is called a Hilbert-Schmidt operator if
$$\sum_{n=1}^\infty |Ae_n|^2<\infty,$$
 where 
$\{e_n\}_{n=1}^\infty$ is an orthonormal basis of \H.
And $A$ is called a trace class operator on \H if 
$$\sum_{n=1}^\infty \<(A^*A)^{1/2}e_n,e_n\><\infty. $$


The celebrated L\'evy-Khintchine formula (\cite[Theorem VI.4.10, Page 182]{Par67}) asserts that 
the characteristic exponent $\psi$ can be written as
\begin{equation}\label{Lambda_ts}
   \psi(\xi)=-i\<a,\xi\>+\frac12\<R\xi,\xi\>
		-\int_\H \left(\e^{i\<\xi,x\>}
		-1-\frac{i\<\xi,x\>}{1+|x|^2} \right)\,\m(dx),
\end{equation}
for some element \(a\in\H\), 
some non-negative definite, symmetric  trace class operator  \(R\)  on \H,  
and some  L\'evy measure \(\m\) on \H, i.e.\ a $\sigma$-finite measure concentrating on $\H\setminus \{0\}$ satisfying 
\[
   \int_{\H} (1\wedge |x|^2) \,\m(dx)<+\infty. 
\]

We call the triplet $(a,R,\m)$ the characteristic of the measure $\mu$. 
For simplicity, we shall write 
\[
        \mu=[a,R,\m]. 
\]

\subsection{Motivations}\label{subsec:Mot}
We recall some backgrounds on additive process  
and convolution equations (special cases of \eqref{Equ:SkewMeasure}) etc..

\textbf{Additive processes, L\'evy processes and convolution equations}.
Let $(X_t)_{t\in\R}$ be a stochastic process taking values in \H. Assume that $X_0=0$. The process $(X_t)_{t\in\R}$ is called an additive process if
it has independent increments, i.e.\ if for any $t>s$, $X_t-X_s$ is independent of $\sigma(\{X_r\colon r\leq s\})$. 
Let $\mu_{t,s}$ denote the distribution of $X_t-X_s$. 
For  $t\geq r\geq s$, we have  $$X_t-X_s=(X_t-X_r)+(X_r-X_s).$$
Note that $X_t-X_r$, $X_r-X_s$ are independent, 
so we have 
\begin{equation}\label{convolution_equation:additiv_process}
   \mu_{t,s}=\mu_{t,r}*\mu_{r,s},\quad t\geq r\geq s.
\end{equation}

Usually we require that the additive process is stochastic continuous (in other words, continuous in probability): for every $t\in\R$ and $\varepsilon>0$, 
\begin{equation}\label{additive_process:stoch_continuous}
     \lim_{h\to 0}\P(|X_{t+h}-X_t|\geq \varepsilon)=0. 
\end{equation}
This condition means exactly that  
 \begin{equation}\label{additive_process:stoch_continuous''}
 \begin{aligned}
    \mu_{s,t}&\Rightarrow \delta_0, \quad \textrm{as}\ s\uparrow t,\\
    \mu_{s,t}&\Rightarrow \delta_0, \quad \textrm{as}\ t\downarrow s.
 \end{aligned}
\end{equation}

If in addition the increments of $(X_t)_{t\in\R}$ is stationary, i.e.\ if for any $t>s$ the distribution of $X_t-X_s$  only depends on $t-s$, then we call it a L\'evy process. In this case we shall only consider $X_t$ for ${t\geq 0}$. For each $t\geq 0$
let $\mu_t$ denote the distribution of $X_t$. Then obviously 
we have the following convolution equation 
\begin{equation}\label{convolution_equation:Levy_process}
  \mu_{t+s}=\mu_t*\mu_s,\quad t,s\geq 0.
\end{equation}
This is the most simple case of \eqref{Equ:SkewMeasure}.
The stochastic continuity condition \eqref{additive_process:stoch_continuous} is reduced to 
\begin{equation}\label{additive_process:stoch_continuous'}
     \lim_{t\downarrow 0}\P(|X_t|\geq \varepsilon)=0. 
\end{equation}
It is equivalent to 
	$$\mu_t\Rightarrow \delta_0\quad  \textrm{as}\ t\downarrow 0.$$ 

We turn to consider a family of probability measures $(\mu_t)_{t\geq 0}$ satisfying  \eqref{convolution_equation:Levy_process}. 
It is clear that for every $t\geq 0$, the measure $\mu_t$ is infinitely divisible. Moreover, 
if $\mu_t$ is weakly continuous, i.e.\ $\mu_t\Rightarrow \delta_0$ as $t\to 0$, then we have
(cf. \cite[Section 8]{BF75} or \cite[Theorem 5.2.3]{Hey10})
\begin{equation}\label{Representation:mu_t}
  \hat{\mu}_t(\xi)=\exp(-t \psi(\xi)),\quad t\geq 0, \ \xi\in \H
\end{equation}
with $\psi=-\log\hat{\mu}_1$. Note that \eqref{Representation:mu_t} is no longer true without the assumption of continuous conditions (cf. \cite[Section 14.4]{Bre92}). 
Conversely, given a function $\psi$ on \H which is the characteristic exponent of some infinitely divisible measure, \eqref{Representation:mu_t} determines a convolution semigroup $(\mu_t)_{t\geq 0}$ satisfying \eqref{convolution_equation:Levy_process}. This is called Schoenberg correspondence (cf. \cite[Theorem 5.2.3]{Hey10}). 

\textbf{Skew convolution equations}. 
We now proceed to consider the skew convolution equation \eqref{SkewConvolutionMeasure} 
\begin{equation}\label{SkewConvolutionMeasure'} 
  \mu_{t+s}=\mu_s *(\mu_t\circ T_s^{-1}) ,\quad s,t\geq 0
  \tag{\ref{SkewConvolutionMeasure}$'$}
\end{equation}
which  generalizes \eqref{convolution_equation:Levy_process}. 
It is proved in \cite{SS01} that $(\mu_t)_{t\geq 0}$ satisfying (\ref{SkewConvolutionMeasure}$'$)
is also infinitely divisible. Moreover, by 
\cite[Theorem 2.3]{DLSS04}, $(\hat{\mu}_t)_{t\geq 0}$  has the following representation
\[
    \hat{\mu}_{t}(\xi) =\exp\left( - \int_0^t \psi_s(\xi) \,ds \right),\quad \xi\in \H,\  t\geq 0,
\]
where $(\psi_s)_{s>0}$ is the characteristic exponents of a family of infinitely divisible measures such that 
$$\psi_{t+s}(\xi)=\psi_s(T_t^*\xi),\quad\xi\in\H,\ t,s> 0.$$ 
If in particular there exists an infinitely divisible measure with characteristic exponent $\psi_0$ such that $\psi_{t}(\xi)=\psi_0(T_t^*\xi)$ for any $t\geq 0$, then we have 
\[
    \hat{\mu}_{t}(\xi) =\exp\left( - \int_0^t \psi_0(T_s^*\xi) \,ds \right),\quad\xi\in \H,\ t\geq 0. 
\]
The formula above was first derived in \cite{BRS96} with  
$$\psi_0(\xi)=-\frac{d}{dt}\hat{\mu}_t(\xi)|_{t=0},\quad \xi\in\H$$ 
under the assumption that for every $\xi\in\H$ the map
$t\mapsto \hat{\mu}_{t}(\xi)$  is absolutely continuous on $[0,+\infty)$ and differentiable at $t=0$.
  
\textbf{Time inhomogeneous convolution equations}.  
Now we look at the two-parameter convolution equation \eqref{convolution_equation:additiv_process}.

\emph{Infinite divisibility}.
First we consider the finite dimensional case when $\H=\R^d$. 
It is known (see \cite{Ito06} or \cite[Theorem 9.1 and Theorem 9.7]{Sat99}) that if $(\mu_{t,s})_{t\geq s}$
 satisfies (cf. \eqref{additive_process:stoch_continuous''})
 \begin{equation}\label{equ:weak:continuity:limit}
 \lim_{t\downarrow s}\mu_{t,s}=\lim_{s\uparrow t}\mu_{t,s}= \delta_0,
 \end{equation}
 then  for any $t\geq s$, $\mu_{t,s}$ is infinitely divisible. 
 The idea of the proof given in \cite{Sat99} is described as follows. 
 
 Using the weak continuity of $(\mu_{t,s})_{t\geq s}$ on the diagonal, we first show that $(\mu_{t,s})_{t\geq s}$ is uniformly weak continuous on any finite interval $[s_0,t_0]$. That is, for every $\varepsilon>0$ and $\eta>0$, there is $\delta>0$ such that for all $s$ and $t$ in $[s_0,t_0]$ satisfying $0\leq t-s\leq \delta$, we have  (cf. Lemma \ref{lmm:Ass:ContinuousMeasure'})
 $$\mu_{t,s}(|x|>\varepsilon)<\eta.$$
By the celebrated Kolmogorov-Khintchine limit theorems on sums of independent random variables (see \cite[Theorem 9.3]{Sat99}), we obtain that $\mu_{t_0,s_0}$ is infinitely divisible. 

The uniform weak continuity of $(\mu_{t,s})_{t\geq s}$ on $[s_0,t_0]$ can be proved by constructing a stochastic continuous additive process $(X_t)_{t\in \R}$ such that for any $t\geq s$ the increment $X_t-X_s$ has the distribution $\mu_{t,s}$ (see \cite[Theorem 9.7 (ii)]{Sat99} and \cite[Lemma 9.6]{Sat99}).

The arguments above can be easily generalized to the infinite dimensional case.  However， in  Subsection \ref{Subsec:Inf_Div} below we shall go further to study the infinite divisibility of $(\mu_{t,s})_{t\geq s}$ satisfying the skew convolution equation \eqref{Equ:SkewMeasure}. We still use similar strategy but we have to deal with the problem more carefully. The main difficulty is getting the uniform weak continuity for some proper measures. 
To this aim,  we provide in Subsection \ref{Subsec:Inf_Div} a direct proof (see Lemma \ref{Lemma:UniformStochContinuous}), while in Subsection \ref{Subsec:additive_process} we present another proof for Lemma \ref{Lemma:UniformStochContinuous} by constructing a proper stochastic continuous additive process. 
  
\emph{Natural additive processes and their factoring.}
Let $(Z_t)_{t\in\R}$ be an stochastic continuous additive process taking values in \H. For every $t>s$, let $\mu_{t,s}$ denote the distribution of $Z_t-Z_s$ and  
suppose that $\mu_{t,s}=[a_{t,s}, R_{t,s}, \m_{t,s}]$. 

\begin{definition}
We say $Z_t$ is \emph{natural} if for every $t\in\R$ and $\xi\in\H$, the function
$s\mapsto \<a_{t,s},\xi\>$ is locally of bounded variation. 
\end{definition}

\begin{definition}
Let $(\lambda_s)_{s\in \R}$ be a family of functions on \H such that
\begin{enumerate}
\item For every $s\in\R$, $\lambda_s$ is the characteristic exponent of some infinitely divisible probability measure on \H.
\item For every $\xi\in\H$, the map $\R\ni s\mapsto \lambda_s(\xi)$ is measurable in $s$. 
\end{enumerate}
\end{definition}

Let $\sigma$ be a continuous locally of finite measure on $(\R,\B(\R))$, i.e. a measure on $(\R,\B(\R))$
such that  $\sigma(\{t\})=0$ for all $t\in\R$ and 
$\sigma(B)<\infty$ for all finite interval $B$ in $\R$.  

If we have 
\begin{equation}
      \hat{\mu}_{t,s}(\xi)=\exp\left( -\int_s^t \lambda_u (\xi)\,\sigma(du) \right),\quad \xi\in\H,\ t\geq s,
\end{equation}
then we call $((\lambda_s)_{s\in\R}, \sigma)$ a \emph{factoring} of $(\mu_{t,s})_{t\geq s}$ (or a factoring of the additive process if the distribution of its increments has this factoring).

\begin{example}
Let $(Z_t)_{t\geq 0}$ be a stochastic continuous L\'evy process with $Z_0\equiv 0$. 
Then for all $t\geq s\geq 0$, the distribution $\mu_{t,s}$ of $Z_t-Z_s$ (i.e. the distribution of $Z_{t-s}$) has the form
$$\hat{\mu}_{t,s}(\xi)=\exp(-(t-s)\lambda_1(\xi)),$$
 where $\lambda_1$ is the characteristic exponent of $Z_1$.
So obviously  $(\mu_{t,s})_{t\geq s}$ has a factoring $(\lambda_1, ds)$, 
where $ds$ is the Lebesgue measure. 
\end{example}

It is proved in \cite{Sat06} that  a finite dimensional additive process  admits a factoring if and only if it is natural. The generalization of this statement to the infinite dimensional case will be studied in Sub-subsection \ref{subsubsec:time_homo_factoring}. Moreover, we shall study similar factoring for the time inhomogeneous case in   Subsection 
\ref{subsec:ConsRep2}


For the use of the factoring of additive process, we refer to Example \ref{Exa:represation:factoring}. 

\subsection{Weak continuity}\label{subsec:cont}
In the rest of this chapter we assume that $(\mu_{t,s})_{t\geq s}$  satisfies the measure equation \eqref{Equ:SkewMeasure} and we shall use the following assumption.

\begin{assumption}\label{Ass:ContinuousMeasure'}
Assume that for all $t\geq s$, 
 \begin{equation}\label{equ:weak:continuity:limit-assumption}
 \lim_{t\downarrow s}\mu_{t,s}=\lim_{s\uparrow t}\mu_{t,s}= \delta_0.
 \end{equation}
\end{assumption}


In order to express  \eqref{equ:weak:continuity:limit-assumption} more explicitly, we need the following simple fact. 

\begin{lemma}\label{lmm:Ass:ContinuousMeasure'}
Let $(\mu_n)_{n\geq 1}$ be a sequence of probability measures on $(\H,\B(\H))$. 
Then $\mu_n\Rightarrow \delta_0$ as $n\rightarrow \infty$ if and only if\,  
for all $\varepsilon >0$, 
\begin{equation}\label{equ:mu_n20}
	\lim_{n\to\infty}\mu_n(\{x\in\H\colon |x|\geq \varepsilon\})=0.
\end{equation}
\end{lemma}
\begin{proof}
     Suppose that $\mu_n\Rightarrow \delta_0$ as $n\rightarrow \infty$. Then by 
     the Portmanteau theorem (see for instance \cite[Theorem 2.1, Page 16]{Bil99}), 
     \[
          \limsup_{n\to\infty} \mu_n(F)\leq \delta_0(F)
     \]
     for all closed set $F$ in \H. Obviously $\{x\in\H\colon |x|\geq \varepsilon\}$ is closed and 
     $$\delta_0(\{x\in\H\colon \{|x|\geq \varepsilon\})=0.$$ Hence \eqref{equ:mu_n20} follows immediately. 
     
     Now we assume that \eqref{equ:mu_n20} holds for all $\varepsilon>0$. Let $f$ be an arbitrary continuous 
     bounded functions on \H. Suppose that $|f|\leq M$ for some $M>0$.
     We are going to show $\mu_n(f)\rightarrow \delta_0(f)$ as $n\rightarrow\infty$.
     
     By the continuity of $f$ we get that for any $\eta>0$, there exists a $\varepsilon_0>0$ such that for all
     $|x|<\varepsilon_0$, 
     \begin{equation}\label{f(x)-f(0)}
          |f(x)-f(0)|<\frac{\eta}{2}. 
     \end{equation}
     By \eqref{equ:mu_n20} there exists  a constant $N>0$ such that for all $n>N$,
     \begin{equation}\label{mu_nAeps}
         \mu_n(\{x\in\H\colon  |x|\geq \varepsilon_0\})<\frac{\eta}{4M}. 
     \end{equation}

    From \eqref{f(x)-f(0)} and \eqref{mu_nAeps} we obtain
    \[
        \begin{aligned}
                   &   \left| \int_\H f\,d\mu_n -  \int_\H f\,d\delta_0 \right|    \\         \leq &  \int_\H |f(x)-f(0)|\,d\mu_n \\
                =&   \int_{\{x\in\H\colon |x|\geq \varepsilon_0\}} |f(x)-f(0)|\,d\mu_n 
               	  +
			\int_{\{x\in\H\colon |x|< \varepsilon_0\}} |f(x)-f(0)|\,d\mu_n      
                \\
            \leq& 2M\mu_n(\{x\in\H\colon |x|\geq \varepsilon_0\}) + \sup_{|x|< \varepsilon_0}|f(x)-f(0)|\cdot \mu_n
            			(\{x\in\H\colon |x| < \varepsilon_0\}) \\
	   <&2M\cdot \frac{\eta}{4M} + \frac{\eta}{2}=\eta.		
        \end{aligned}
    \]
    This completes the proof. 
\end{proof}

By Lemma \ref{lmm:Ass:ContinuousMeasure'} we have the following equivalent description for Assumption
\ref{Ass:ContinuousMeasure'}.
  
\begin{prop}\label{Ass:ContinuousMeasure}
Condition  \eqref{equ:weak:continuity:limit-assumption} is equivalent to 
\begin{equation}
     \lim_{t\downarrow s}  \mu_{t,s}
            (\{x\in\H\colon |x|>\varepsilon\})
=\lim_{s\uparrow t}\mu_{t,s}
            (\{x\in\H\colon |x|>\varepsilon\})            
                        = 0
\end{equation}
for all $\varepsilon>0$.
More precisely, they are equivalent to the following two conditions: 
for every $\varepsilon, \eta>0$, and for every ${u}\in \R$, there exists a constant $\delta_{u}$ such that 
 \begin{enumerate}
  \item For every $t \in (u, u+\delta_{u})$ we have  
        \begin{equation}\label{eq1:Ass:ContinuousMeasure}
              \mu_{t,{u}}(\{x\in\H \colon |x|>\varepsilon \})<\eta.
        \end{equation}
  \item For every $s \in (u-\delta_{u},u)$ we have  
      \begin{equation}\label{eq2:Ass:ContinuousMeasure}
         \mu_{{u},s}(\{x\in\H \colon |x|>\varepsilon \})<\eta . 
      \end{equation}
 \end{enumerate}
\end{prop}


It is easy to see that  \eqref{equ:weak:continuity:limit-assumption} implies (hence is equivalent to) the following condition: 
for every $u\in\R$, 
\begin{equation}\label{equ:Ass:ContinuousMeasure'-ts-2u}
 \mu_{t,s}\Rightarrow \delta_0\quad 
 \textrm{as}\ t \downarrow u\ \textrm{and}\ s\uparrow u.
\end{equation}
Indeed, for all $t\geq u\geq s$, by \eqref{Equ:SkewMeasure} 
we have 
\[
	\lim_{t\downarrow u,s\uparrow u}\mu_{t,s}
         =\lim_{t\downarrow u,s\uparrow u}\mu_{t,u}*\mu_{u,s}\circ U(t,u)^{-1}
         =\delta_0*\delta_0=\delta_0. 	 
\]
Here we have omitted some details concerning the second equality for which we refer to the proof of Theorem \ref{thm:mu_ts:continuity}.

Note that $\mu_{u,u}=\delta_0$ for all $u\in\R$.
So Assumption \ref{Ass:ContinuousMeasure'} simply says that $\mu_{t,s}$ is weakly continuous on
the diagonal $\{(u,u)\colon u\in\R\}$. We are able to show that 
the continuity of $\mu_{t,s}$ on the diagonal implies in particular that for every fixed $t\in\R$, the map $s\mapsto \mu_{t,s}$ over $(-\infty,t]$ is weak continuous besides some other results. This will be used in subsection \ref{subsec:ConsRep2} for the spectral representation of $\mu_{t,s}$. Before we show it, we need to prove a  lemma. 

First of all, for convenience we include here two results from \cite{Par67}. A set of probability measures on 
$(\H,\B(\H))$
is said to be shift (relatively) compact if for every sequence $(\mu_n)_{n\geq 1}$  there is a sequence $(\nu_n)_{n\geq 1}$ such that 
\begin{enumerate}
\item $(\nu_n)_{n\geq 1}$ is a translate of $(\mu_n)_{n\geq 1}$. That is, there exists a sequence $\{(x_n)_{n\geq 1}\}$ in $\H$ such that $\nu_n=\mu_n*\delta_{x_n}$ holds for all $n\geq 1$.
\item $(\nu_n)_{n\geq 1}$ has a convergent subsequence. 
\end{enumerate}

\begin{theorem}[\textnormal{\cite{Par67}}]\label{Theorem_Par67_III_2_1}
Let $(\sigma_n)_{n\geq 1}$, $(\mu_n)_{n\geq 1}$ and $(\nu_n)_{n\geq 1}$ be three sequences of measures on \H such that $\sigma_n=\mu_n*\nu_n$ for all $n\in \N$. 
\begin{enumerate}
\item \textnormal{(\cite[Theorem III.2.1, Page 58]{Par67})} If the sequences $(\sigma_n)_{n\geq 1}$ and $(\mu_n)_{n\geq 1}$ are relatively compact, then so is the sequence $(\nu_n)_{n\geq 1}$. 
\item \textnormal{(\cite[Theorem III.2.2, Page 59]{Par67})} If the sequence $(\sigma_n)_{n\geq 1}$ is relatively compact then the sequences $(\mu_n)_{n\geq 1}$ and $(\nu_n)_{n\geq 1}$ are shift compact, respectively. 
\end{enumerate}
\end{theorem}

Now we can show the following lemma.

\begin{lemma}\label{lmm:convolution:continuity}
  Let $\mu_n, \nu_n, \sigma_n$ with $n\geq 1$,  $\mu,\nu$ and $\sigma$  be measures on a complete separable metric space such that
  $\sigma_n=\mu_n*\nu_n$. 
  \begin{enumerate}
   \item If $\mu_n\Rightarrow \mu$ and $\nu_n\Rightarrow \nu$ as $n\to \infty$, then 
              $\sigma_n\Rightarrow \mu*\nu$ as $n\to\infty$. 
   \item 
            Suppose that $\sigma_n\Rightarrow \sigma$ and $\mu_n\Rightarrow \mu$ as $n\to \infty$. Then there exists a probability measure $\nu$ such that 
            \begin{equation}\label{equ:sigma=mu*nu}
            	\sigma=\mu*\nu.
            \end{equation}
If $\nu$ is a unique measure such that \eqref{equ:sigma=mu*nu} holds, then 
            	    $\nu_n\Rightarrow \nu$ as $n\to\infty$. 
  \end{enumerate}
\end{lemma}
\begin{proof}
 The first conclusion says that the convolution operation keeps the weak continuity. The proof can be found, for example, in \cite[Proposition 2.3]{HM11} or \cite[Theorem III.1.1, Page 57]{Par67}. 
 
 Now we show the second conclusion.  Take an arbitrary subsequence $(\nu_{n_i})_{i\geq 1}$ from 
 $(\nu_n)_{n\geq 1}$ and consider 
 $$\sigma_{n_i}=\mu_{n_i}*\nu_{n_i},\quad i\geq 1.$$ 
 Since $\sigma_n\Rightarrow \sigma$ and $\mu_n\Rightarrow \mu$ as $n\to \infty$, both $(\sigma_{n_i})_{i\geq 1}$ and $(\mu_{n_i})_{i\geq 1}$ are relatively compact. 
 By \cite[Theorem III.2.1, Page 58]{Par67} (see Theorem \ref{Theorem_Par67_III_2_1}), the sequence $(\nu_{n_i})_{i\geq 1}$ is also relatively compact. Let $(\nu_{n'_i})_{i\geq 1}$ be a weak convergent subsequence of  $(\nu_{n_i})_{i\geq 1}$ with limit $\nu'$. Then by the first conclusion of this theorem, we have
 \[
     \sigma_{n'_i}=\mu_{n'_i}*\nu_{n'_i}
       \Rightarrow \mu*\nu',\quad n_i'\to\infty. 
 \]
 Since we also have $\sigma_{n_i'}\Rightarrow \sigma$ as $n_i'\to\infty$, we get 
 $\sigma=\mu*\nu'$. 
 This shows that there exists a probability measure $\nu=\nu'$ such that \eqref{equ:sigma=mu*nu} holds. If there is only one measure  $\nu$ such that \eqref{equ:sigma=mu*nu} holds, then  any subsequence of $(\nu_n)_{n\geq 1}$ contains a further subsequence converging weakly to $\nu$. 
 This is sufficient  to conclude that $(\nu_n)_{n\geq 1}$ converges weakly to $\nu$ (cf.  \cite[Theorem 2.6, Page 20]{Bil99}). Hence the proof is complete. 
\end{proof}

\begin{remark}\label{Remark:cancellation:law} 
In the second part of the previous thereom, the assumption that $\nu$ is  the unique solution to convolution equation $\sigma=\mu*\nu$ amounts to say that the following 
cancellation law for convolution operation holds: Let $\nu,\nu'$ be two measures on $\H$, if 
\begin{equation}\label{Equ:mu_nu=mu_nu'}
	\mu*\nu=\mu*\nu',
\end{equation}
then the measure $\mu$ can be canceled on both sides of the equation and get $\nu=\nu'$. It is obvious that this cancellation law holds provided  $\hat{\mu}\neq 0$. Indeed, from \eqref{Equ:mu_nu=mu_nu'} we get 
$$\hat{\mu}\hat{\nu}=\hat{\mu}\hat{\nu'}.$$
Since $\hat{\mu}\neq 0$ we have $\hat{\nu}=\hat{\nu'}$. So $\nu=\nu'$. 
It is well known that $\hat{\mu}\neq 0$ if $\mu$ is an infinitely divisible distribution. 
\end{remark}

\begin{remark}
	A long time after the proof of the second assertion of Lemma \ref{lmm:convolution:continuity}, we found that there is a similar result in \cite[Corollary 2.2.4, Page 40]{Hey10} where condition $\hat{\mu}\neq 0$ is used. Our proof is different. The example after \cite[Theorem 5.1.1]{Luk70} shows that there exist probability measures $\mu,\nu$ and $\nu'$ on $\R$ with $\nu\neq \nu'$ such that \eqref{Equ:mu_nu=mu_nu'} holds. The fact which called Khintchine phenomenon shows that the cancellation law for convolution operation is not valid in general.  So the condition that \eqref{Equ:mu_nu=mu_nu'} implies $\nu=\nu'$ is necessary for  the second assertion of Lemma {lmm:convolution:continuity}. Otherwise, if $\nu\neq \nu'$, then $(\mu*\nu_n)_{n\geq 1}$ with $\nu_{2k-1}=\nu$ and $\nu_{2k}=\nu'$ for all $k\geq 1$ converges weakly, but  $(\nu_n)_{n\geq 1}$ does not converge weakly.  
\end{remark}

As a summary of the discussion above, we have the following result. 

\begin{cora}\label{cora:convolution:continuity''}
  Let $\mu_n, \nu_n, \sigma_n$  with $n\geq 1$  and $\sigma,\mu$ be  measures on $\H$ with the following properties
  \begin{enumerate}
  \item 
   For all  $n\geq 1$, 
  $\sigma_n=\mu_n*\nu_n$;
  \item As $n\to \infty$, 
	$\sigma_n\Rightarrow \sigma$ and $\mu_n\Rightarrow \mu$. 
\end{enumerate}	
If $\mu$ is an infinitely divisible distribution, 
then  $(\nu_n)_{n\geq 1}$ converges weakly to 
some  measure $\nu$ on $\H$
such that $\sigma=\mu*\nu$ as $n\to \infty$.
\end{cora}


In particular, we have the following result. 

\begin{cora}\label{cora:convolution:continuity}
  Let $\mu_n, \nu_n, \sigma_n$ with $n\geq 1$, and $\sigma$ be measures on \H. Suppose that for all $n\geq 1$, 
  $\sigma_n=\mu_n*\nu_n$. 
If 
$\sigma_n\Rightarrow \sigma$ and $\mu_n\Rightarrow \delta_0$ as $n\to \infty$, 
then $\nu_n\Rightarrow \sigma$ as $n\to \infty$. 
\end{cora}

We shall need the following fact.

\begin{lemma}\label{lem:U_boundedness}
The evolution family 
$(U(t,s))_{t\geq s}$ is uniformly bounded on every compact interval. That is,
for every $s_0<t_0$, there exists some constant $c\geq 1$ such that for all $s_0\leq s\leq t\leq t_0$, 
\begin{equation}\label{equ:U_boundedness}
                   |U(t,s)x|\leq c|x|,\quad x\in\H, \ s_0\leq s\leq t\leq t_0.
\end{equation}
\end{lemma}
\begin{proof}
For every $x\in\H$, $|U(t,s)x|$ is a continuous function in $(t,s)$ on $\Lambda_{t_0,s_0}:=\{(t,s)\colon s_0\leq s\leq t\leq t_0\}$. Hence $|U(t,s)x|$ is uniformly bounded on $\Lambda_{t_0,s_0}$ for every $x\in\H$. By the
Banach-Steinhaus theorem we have 
$$\sup_{(t,s)\in \Lambda_{t_0,s_0}} \|U(t,s)\| <\infty.$$
That is, there exists some $c>0$ such that \eqref{equ:U_boundedness} holds. Noting  that $U(t,t)=I$ for all $t\in\R$, it is clear that $c\geq 1$.
\end{proof}

We shall use the following simple fact several times, so we formulate it as  a Lemma. 

\begin{lemma}\label{Lemma:measure:contraction}
	Let $T$ be a bounded operator on \H such that for all $x\in\H$, $|Tx|\leq c |x|$ for some constant $c>0$. Let $\mu$ be a measure on $(\H,\B(\H))$ and $\varepsilon>0$ be a constant. Then we have
	\begin{equation}\begin{aligned}
		  \mu\circ T^{-1}(\{x\in\H\colon |x|>\varepsilon\})
	      =&\mu(\{x\in\H\colon |Tx|>\varepsilon\})\\
	      \leq & \mu(\{x\in\H\colon |x|>\varepsilon/c\}).
	 \end{aligned}     
	\end{equation}	
\end{lemma}

Now we can prove the following main result. 

\begin{theorem}\label{thm:mu_ts:continuity}
Suppose that Assumption \ref{Ass:ContinuousMeasure'} holds. The following assertions hold.  
\begin{enumerate}
\item For every $t\in\R$, the map $s\mapsto \mu_{t,s}$ with $s\leq t$ is weakly continuous. 
\item  For every $t,s\in\R$ with $t\geq s$ we have 
	  \begin{equation}\label{equ1:pf:mu_ts:continuity}
	     \mu_{t+\varepsilon,s} \Rightarrow \mu_{t,s}\quad
		\textrm{as}\ \varepsilon \downarrow 0.
	  \end{equation}
\item Then we have for every $t,s\in\R$ with $t>s$
	\begin{equation}\label{Equ:mu_t-epsilon-continuity}
		\mu_{t-\varepsilon,s}\circ U(t,t-\varepsilon)^{-1}\Rightarrow
		 \mu_{t,s},\quad
		 \varepsilon\downarrow 0.
	\end{equation}
	  
\end{enumerate}
\end{theorem}

\begin{proof}
(i)  Let $s<t$. 
We need to show both 
  \begin{equation}\label{equ4:pf:mu_ts:continuity_new}
	  \mu_{t,s-\varepsilon} \Rightarrow \mu_{t,s},\quad \textrm{as}\ \varepsilon \downarrow 0 
 \end{equation}	  
and
  \begin{equation}\label{equ5:pf:mu_ts:continuity_new}
	  \mu_{t,s+\varepsilon} \Rightarrow \mu_{t,s},\quad \textrm{as}\ \varepsilon \downarrow 0 .  
 \end{equation}	  

Following from \eqref{Equ:SkewMeasure} we have for every $\varepsilon\in(0,t-s)$
  \begin{equation}\label{equ1:pf:mu_ts:continuity_new}
	     \mu_{t,s-\varepsilon} =  \mu_{t,s}*\left( \mu_{s,s-\varepsilon}\circ U(t,s)^{-1}  \right)
  \end{equation}
and
  \begin{equation}\label{equ2:pf:mu_ts:continuity_new}
	     \mu_{t,s} =  \mu_{t,s+\varepsilon}*\left( \mu_{s+\varepsilon,s}\circ U(t,s+\varepsilon)^{-1}  \right).
  \end{equation}

By Lemma \ref{lem:U_boundedness} there exists some constant $c\geq 1$ such that  for all  $\varepsilon\in[0,t-s]$, we have 
$$\| U(t,s+\varepsilon)\|\leq c.$$ 
Hence by   Lemma
\ref{Lemma:measure:contraction} we have for all $\eta>0$
  \begin{equation}\label{equ3:pf:mu_ts:continuity_new}
	\mu_{s,s-\varepsilon}\circ U(t,s)^{-1}(\{x\in\H\colon |x|>\eta\})\leq 
		\mu_{s,s-\varepsilon} (\{x\in\H\colon |x|>\eta/c\}).
\end{equation}
Because  $\mu_{s,s-\varepsilon}\Rightarrow \delta_0$ as $\varepsilon\to 0$, by Lemma \ref{lmm:Ass:ContinuousMeasure'} we get 
\[
	\lim_{\varepsilon\downarrow 0}\mu_{s,s-\varepsilon} (\{x\in\H\colon |x|>\eta/c\})=0.
\]
Hence it follows from \eqref{equ3:pf:mu_ts:continuity_new} we obtain that 
\[
	\lim_{\varepsilon\downarrow 0} \mu_{s,s-\varepsilon}\circ U(t,s)^{-1}(\{x\in\H\colon |x|>\eta\}) =0.
\]
By Lemma \ref{lmm:Ass:ContinuousMeasure'}, we obtain 
$$\mu_{s,s-\varepsilon}\circ U(t,s)^{-1} \Rightarrow \delta_0,\quad \varepsilon \downarrow 0.$$
Therefore, applying the first result of Lemma 
\ref{lmm:convolution:continuity} to \eqref{equ1:pf:mu_ts:continuity_new} we get \eqref{equ4:pf:mu_ts:continuity_new}.

By the same arguments, it is easy to show that
 \[
 	\mu_{s+\varepsilon,s}\circ U(t,s+\varepsilon)^{-1} \Rightarrow \delta_0 \quad \textrm{as}\ \varepsilon \downarrow 0 . 
 \]
Then by Corollary \ref{cora:convolution:continuity}, \eqref{equ5:pf:mu_ts:continuity_new} follows from  
 \eqref{equ2:pf:mu_ts:continuity_new}.

So \eqref{equ4:pf:mu_ts:continuity_new} and \eqref{equ5:pf:mu_ts:continuity_new} is proved and hence the proof is complete. 

(ii) 				
According to \eqref{Equ:SkewMeasure} we have for all $t\geq s$, $\varepsilon \geq  0$, 
  \[
     \mu_{t+\varepsilon,s} =
        \mu_{t+\varepsilon,t}*(\mu_{t,s}\circ U(t+\varepsilon,t)^{-1}).
  \]
By assumption we have  $ \mu_{t+\varepsilon,t}\Rightarrow \delta_0$. 
Hence by applying the first assertion of Lemma \ref{lmm:convolution:continuity}, 
we get \eqref{equ1:pf:mu_ts:continuity} providing 
\begin{equation}\label{Equ:mu_ts_U_t+epsilon}
	\mu_{t,s}\circ U(t+\varepsilon,t)^{-1}\Rightarrow \mu_{t,s},\quad
		\textrm{as}\  \varepsilon \downarrow 0.
\end{equation}
Now we show \eqref{Equ:mu_ts_U_t+epsilon}. 

Let $f$ be a  continuous and bounded function on $\H$. For every $\varepsilon>0$ set 
$$f_{\varepsilon}(x)=f(U(t+\varepsilon,t)x),\quad x\in\H.$$
It is clear that $f_\varepsilon$ converges to $f$ pointwisely and $f_\varepsilon$ is bounded since $f$ is bounded. Hence by bounded convergence theorem we have 
\[
	\lim_{\varepsilon\downarrow 0} \int_{\H} f\,d\mu_{t,s}\circ U(t+\varepsilon,t)^{-1} 
	=\lim_{\varepsilon\downarrow 0} \int_{\H} f(U(t+\varepsilon,t)x)\,d\mu_{t,s}(x)
	= \int_{\H} f(x)\,d\mu_{t,s}(x).
\]
This proves \eqref{Equ:mu_ts_U_t+epsilon}.

	%

(iii)  By \eqref{Equ:SkewMeasure} we have for all $t\geq t-\varepsilon >s$, 
  \[
     \mu_{t,s} =
        \mu_{t,t-\varepsilon}
           *(\mu_{t-\varepsilon,s}\circ U(t,t-\varepsilon)^{-1}).
  \]
Since $\mu_{t,t-\varepsilon}\Rightarrow \delta_0$ as $\varepsilon\downarrow 0$, by the second conclusion of Lemma \ref{lmm:convolution:continuity},
we have \eqref{Equ:mu_t-epsilon-continuity}.


 \end{proof}

By Theorem \ref{thm:mu_ts:continuity} we have the following result concerning the space-homogeneous case. 
\begin{prop}
Let $(\mu_{t,s})_{t\geq s}$ be a family of  probability measures on $(\H,\B(\H))$ satisfying 
\begin{equation}\label{Equ:space-homo-convolution-equ-mu_ts}
	\mu_{t,s}=\mu_{t,r}*\mu_{r,s},\quad t\geq r\geq s
\end{equation}
and 
\[
	\lim_{t\downarrow s}\mu_{t,s}=
	\lim_{s\uparrow t}\mu_{t,s}=
	\delta_0.
\]
Then $\mu_{t,s}$ is weak continuous in $t\geq s$.
\end{prop}

\begin{proof}
Let $\varepsilon_n\downarrow 0$ and $\delta_n\downarrow 0$. By \eqref{Equ:space-homo-convolution-equ-mu_ts} we have
\begin{equation}\label{Equ:mu_ts_e-d1}
	\mu_{t,s\pm\delta_n}=\mu_{t,t-\varepsilon_n}*\mu_{t-\varepsilon_n,s\pm\delta_n}
\end{equation}
and 
\begin{equation}\label{Equ:mu_ts_e-d2}
	\mu_{t+\varepsilon_n,s\pm\delta_n}
	=\mu_{t+\varepsilon_n,t}*\mu_{t,s\pm\delta_n}.
\end{equation}
As $n\rightarrow \infty$,  we have 
\(
	\mu_{t,t-\varepsilon_n}\Rightarrow \delta_0
\)
by assumption
and 
\(
	\mu_{t,s\pm\delta_n}\Rightarrow \delta_0
\)
by Theorem \ref{thm:mu_ts:continuity}.
Hence by applying the second assertion of  Lemma \ref{lmm:convolution:continuity} to Equation \eqref{Equ:mu_ts_e-d1} we obtain 
\[
	\mu_{t-\varepsilon_n,s\pm\delta_n}\Rightarrow \mu_{t,s},\quad
	\textrm{as}\ n\rightarrow \infty. 
\]
Similarly, it follows 
\[
	\mu_{t+\varepsilon_n,s\pm\delta_n} 
		\Rightarrow \mu_{t,s},\quad
			\textrm{as}\ n\rightarrow \infty. 
\]
from \eqref{Equ:mu_ts_e-d2} by using the continuity of the convolution  operator. So we have 
\[
	\mu_{t\pm\varepsilon_n,s\pm\delta_n} 
		\Rightarrow \mu_{t,s},\quad
			\textrm{as}\ n\rightarrow \infty
\]
and the proof is complete 
\end{proof}

As a direct application of Theorem \ref{thm:mu_ts:continuity} we obtain  the continuity of the characteristic functional (see Proposition \ref{char_functional_continuity}). 
In the following we study the continuity of $(p_{s,t})_{t\geq s}$.

 For every $t\geq s$, 
it is clear that
$p_{s,t}$ is Feller, i.e.
$p_{s,t}(C_b(\H))\subset C_b(\H)$. 
Now we look at the continuity of the map $(s,x)\mapsto p_{s,t}f(x)$ for every $f$ in $C_b(\H)$.
The proposition below is a direct generalization of \cite[Lemma 2.1]{BRS96}. The proof is quite similar to the proof in \cite{BRS96}.

\begin{prop}\label{Prop:GMS:Continuity}
Let $s_n, t_n\in\R$, $x_n\in\H$, $s_n\leq t_n$ with $n\geq 1$ such that 
$(s_n,t_n)\to (s,t)\in\R^2$ and $x_n\to x\in\H$  as $n\to \infty$.
If $\mu_{t_n,s_n}\Rightarrow \mu_{t,s}$ as $n\to \infty$,
 then for any $f\in C_b(\H)$,  $p_{s_n,t_n}f(x_n)\to p_{s,t}f(x)$ as $n\to \infty$.
\end{prop}

\begin{proof} 
(1) Since 
 $\mu_{t_n,s_n}\Rightarrow \mu_{t,s}$ as $n\to \infty$,  
by Prohorov's theorem,  for every $\varepsilon>0$, there exists a compact set $K\subset \H$ such that 
\begin{equation}\label{Equ1:Prop:GMS:Continuity}
 \mu_{r,\sigma}(K)\geq 1-\varepsilon,\quad \textrm{for\ all}\ (r,\sigma)\in\{(t,s), (t_n,s_n)\colon n\in\N\}.
\end{equation}

For abbreviation, we set $z_n=U(t_n,s_n)x_n$ and $z=U(t,s)x$. By the strong continuity of the evolution family $(U(t,s))_{t\geq s}$, the set
 $Z:=\{z, z_n\colon n\in\N\}$
is compact. Hence $Z+K$ is also compact. So there exists $N\in\N$ such that for any $n>N$ and for any $y\in K$, 
\begin{equation}\label{Equ2:Prop:GMS:Continuity}
	|f(z_n+y)-f(z+y)|<\varepsilon, 
\end{equation}
since $f$ is uniformly continuous on compacts.

Because  $\mu_{t_n,s_n}\Rightarrow \mu_{t,s}$ as $n\to \infty$,  
(taking $N$ larger if necessary) we have for all $n>N$
\begin{equation}\label{Equ3:Prop:GMS:Continuity}
	\left|\int_\H  f(z+y)\, \mu_{t_n,s_n}(dy)
			- 	\int_\H  f(z+y)\, \mu_{t,s}(dy)\right|< \varepsilon.
\end{equation}
From \eqref{Equ1:Prop:GMS:Continuity}, \eqref{Equ2:Prop:GMS:Continuity} and  \eqref{Equ3:Prop:GMS:Continuity} 
we get 
\[
\begin{aligned}
& \left|  \int_\H f(z_n+y)\,\mu_{t_n,s_n}(dy)
	-  \int_\H f(z+y)\,\mu_{t,s}(dy)\right|\\
\leq & \left|\int_\H  f(z+y)\, \mu_{t_n,s_n}(dy)
			- 	\int_\H  f(z+y)\, \mu_{t,s}(dy)\right| \\
	&\qquad\quad\quad\qquad +\int_K \left| f(z_n+y)-f(z+y)  \right|\,\mu_{t_n,s_n}(dy)+2\|f\|_{\infty} \mu(\H\setminus K) \\		
<  & 2\varepsilon (1+\|f\|_\infty).
\end{aligned}
\]
Hence the result is proved since $\varepsilon$ was arbitrary. 

\end{proof}

\subsection{Infinite divisibility}\label{Subsec:Inf_Div}
Here we use similar method introduced in Subsection \ref{subsec:Mot} to show the infinite divisibility of $\mu_{t,s}$ for all $t\geq s$. 
Moreover, we shall give  another proof in Corollary \ref{cora:inf_div}.

The following estimate is crucial for the the proof of  Theorem \ref{Thm:mu_ts_IID} which is the main result of this section. 

\begin{lemma}\label{Lemma:UniformStochContinuous}
Suppose that Assumption \ref{Ass:ContinuousMeasure'} holds. 
On every compact interval $[s_0,t_0]$, for all $\varepsilon, \eta>0$, there exists a constant $\delta>0$ such that for all $s,t\in[s_0,t_0]$ with $0\leq t-s< \delta$, 
\begin{equation}\label{Equ:UniformStochContinuous:Mod}
\mu_{t,s}\circ U(t_0,t)^{-1}(\{x\in\H\colon |x|>\varepsilon\})<\eta.
\end{equation}
\end{lemma}

\begin{proof}
It is  trivial to see that \eqref{Equ:UniformStochContinuous:Mod} holds for the case $t=s$. So we shall assume $t>s$.
By Lemma \ref{lem:U_boundedness}, there exists a constant $c\geq 1$ such that 
\begin{equation}\label{Equ:U:bound:Lemma:UniformStochContinuous}
	|U(t,s)x|\leq c|x|,\quad \ x\in\H,\ s_0\leq s<t\leq t_0.
\end{equation}
Let us set
$$\varepsilon'=\varepsilon/c$$
and denote for every $r>0$
$$A(r):=\{x\in\H\colon |x|>r\}.$$

By Assumption \ref{Ass:ContinuousMeasure'} and Equations \eqref{eq1:Ass:ContinuousMeasure},  \eqref{eq2:Ass:ContinuousMeasure} in Proposition \ref{Ass:ContinuousMeasure}, we have   
for every $\varepsilon,\eta>0$,  $t\in[s_0,t_0]$, there exists a constant $\delta_t\geq 0$ such that
\begin{equation}\label{equ1:mu_ts:Avar_eps_2c}
\mu_{t,s}\left(A\left(\frac{\varepsilon'}{2c}\right)\right)<\eta/2,\quad 
s\in(t-\delta_t,t)
\end{equation}
and
\begin{equation}\label{equ2:mu_ts:Avar_eps_2c}
\mu_{r,t}\left(A\left(\frac{\varepsilon'}{2c}\right)\right)<\eta/2,\quad 
r\in(t, t+\delta_t).
\end{equation}
Since $c\geq 1$ we have $\frac{\varepsilon'}{2c}\leq \frac{\varepsilon'}{2}$. Hence it follows from 
estimates \eqref{equ1:mu_ts:Avar_eps_2c} and \eqref{equ2:mu_ts:Avar_eps_2c} we obtain
\begin{equation}\label{equ1:mu_ts:Avar_eps_2}
\mu_{t,s}\left(A\left({\varepsilon'}/{2}\right)\right)<\eta/2,
\quad s\in(t-\delta_t,t)
\end{equation}
and
\begin{equation}\label{equ2:mu_ts:Avar_eps_2}
\mu_{r,t}\left(A\left({\varepsilon'}/{2}\right)\right)<\eta/2,\quad r\in(t, t+\delta_t).
\end{equation}

Moreover, according to Lemma \ref{Lemma:measure:contraction} and \eqref{Equ:U:bound:Lemma:UniformStochContinuous}, we obtain from  estimates \eqref{equ1:mu_ts:Avar_eps_2c} and \eqref{equ2:mu_ts:Avar_eps_2c}  that  
\begin{equation}\label{equ1:mu_ts:Avar_eps_U2}
\mu_{t,s} \circ U(t',t)^{-1} \left(A\left({\varepsilon'}/{2}\right)\right)<\eta/2,
\quad 
t-\delta_t\leq s\leq t,
\ 
t\leq t'\leq t_0
\end{equation}
and
\begin{equation}\label{equ2:mu_ts:Avar_eps_U2}
\mu_{r,t} \circ U(r',r)^{-1}  \left(A\left({\varepsilon'}/{2}\right)\right)<\eta/2,\quad
t\leq r\leq  t+\delta_t, 
\ r\leq r'\leq t_0.
\end{equation}



For every $t\in[s_0,t_0]$, let  $$I_t:=(t-\delta_t,t+\delta_t).$$ 
It is obviously that $\{I_t\colon t\in[s_0,t_0]\}$ covers the interval $[s_0,t_0]$. Hence there is a finite sub--covering $\{I_{t_j}\colon j=1,2\cdots,n\}$ of $[s_0,t_0]$. 
Then for every $t\in[s_0,t_0]$, we have $t\in I_{t_j}$ for some $j\in \{1,2,\cdots,n\}$. 
Let $\delta$ be the minimum of $\{\delta_{t_j}/2\colon j=1,2,\cdots, n\}$.
For every $s\in[s_0,t_0]$ such that $0<t-s<\delta$,  we have 
\[|s-t_j|\leq |s-t|+|t-t_j|<\delta+\delta_{t_j}/2\leq \delta_{t_j}.\]
Therefore we get that both $t$ and $s$ are in the same sub-interval $I_{t_j}$. 
We need to consider the following three cases respectively: 1. $s\leq t_j<t$; 2. $s<t\leq t_j$; 3. $t_j<s<t$. 

\textbf{Case 1 ($s\leq t_j<t$)}. Note that for all $x,y\in\H$, if $|x+y|>\varepsilon'$, then either $|x|>\varepsilon'/2$ or $|y|>\varepsilon'/2$. That is, the following inequality holds
\begin{equation}\label{Ineq:Eins_Axy_plus}
	\Eins_{A(\varepsilon')}(x+y)\leq 	\Eins_{A(\varepsilon'/2)}(x) + \Eins_{A(\varepsilon'/2)}(y).
\end{equation}

By  \eqref{Equ:SkewMeasure},  
\eqref{equ2:mu_ts:Avar_eps_2}, \eqref{equ1:mu_ts:Avar_eps_U2}
and \eqref{Ineq:Eins_Axy_plus} we have
\[
\begin{aligned}
	    \mu_{t,s}(A(\varepsilon'))
   & =    \mu_{t,t_j}*(\mu_{t_j,s}\circ U(t,t_j)^{-1})(A(\varepsilon'))\\
   & =    \int_\H\int_\H \Eins_{A(\varepsilon')} (x+y)\, \mu_{t,t_j}(dx) ( \mu_{t_j,s}\circ U(t,t_j)^{-1})(dy)\\
   & \leq \int_\H\int_\H (\Eins_{A(\varepsilon'/2)} (x)+\Eins_{A(\varepsilon'/2)} (y))\,
   				 \mu_{t,t_j}(dx) (\mu_{t_j,s}\circ U(t,t_j)^{-1})(dy)\\ 
   &= \mu_{t,t_j}(A(\varepsilon'/2))+ (\mu_{t_j,s}\circ U(t,t_j)^{-1})(A(\varepsilon'/2))\\
   &<\frac{\eta}{2}+\frac{\eta}{2}=\eta.				 
\end{aligned}
\]
Therefore by Lemma \ref{Lemma:measure:contraction} and \eqref{Equ:U:bound:Lemma:UniformStochContinuous} 
we have 
\[
	\mu_{t,s}\circ U(t_0,t)^{-1}(\{x\in\H\colon |x|>\varepsilon\}) < \eta. 
\]	


\textbf{Case 2 ($s<t\leq t_j$)}. We first show 
\[(\mu_{t,s}\circ U(t_j,t)^{-1})(A(\varepsilon'))<\eta\]
 by contradiction.  
If otherwise the following inequality 
\begin{equation}\label{Equ:contradict}
   (\mu_{t,s}\circ U(t_j,t)^{-1})(A(\varepsilon'))\geq \eta
\end{equation} 
holds. 
Then by \eqref{Equ:SkewMeasure},
\eqref{equ1:mu_ts:Avar_eps_2}
 and  \eqref{Equ:contradict} 
we obtain 
\[
\begin{aligned}
	\frac{\eta}{2}
   &>\mu_{t_j,s}(A(\varepsilon'/2))
       = \mu_{t_j,t}*(\mu_{t,s}\circ U(t_j,t)^{-1})(A(\varepsilon'/2))\\
   & =   \int_\H\int_\H \Eins_{A(\varepsilon'/2)} (x+y)\, \mu_{t_j,t}(dx) ( \mu_{t,s}\circ U(t_j,t)^{-1})(dy)\\
   &\geq \int_\H\int_\H \Eins_{A(\varepsilon'/2)^c} (x)\cdot \Eins_{A(\varepsilon')} (y)\,
			\, \mu_{t_j,t}(dx) ( \mu_{t,s}\circ U(t_j,t)^{-1})(dy)\\
   &= \mu_{t_j,t}(A(\varepsilon'/2)^c)\cdot (\mu_{t,s}\circ U(t_j,t)^{-1})(A(\varepsilon'))\\
   &\geq\eta\left(1-\frac{\eta}{2}\right)=\eta-\frac{\eta^2}{2}.   
\end{aligned}   
\]
Here we used the fact that for all $x,y\in\H$,  
$|x+y|\geq |y|-|x|>\frac{\varepsilon'}{2}$
if $|y|>\varepsilon'$ and $|x|\leq \varepsilon'/2$. 
Now we have
 $$\frac{\eta}{2}>\eta-\frac{\eta^2}{2}.$$ 
Consequently, we have $\eta>1$. By \eqref{Equ:contradict} this means that 
$$ (\mu_{t,s}\circ U(t_j,t)^{-1})(A(\varepsilon'))>1.$$ 
This is impossible since $\mu_{t,s}$ is a probability measure.  

Then by Lemma \ref{Lemma:measure:contraction} and 
\eqref{Equ:U:bound:Lemma:UniformStochContinuous} 
we have 
\[
\begin{aligned}
     & \mu_{t,s}\circ U(t_0,t)^{-1}(\{x\in\H\colon |x|>\varepsilon\}) \\
 =&    
    \left( \mu_{t,s}\circ U(t_j,t)^{-1}\right) \circ U(t_0,t_j)^{-1} (\{x\in\H\colon |x|>\varepsilon\}) \\
 \leq &
	 \mu_{t,s}\circ U(t_j,t)^{-1} (\{x\in\H\colon |x|>\varepsilon/c\}) \\
 = &   \mu_{t,s}\circ U(t_j,t)^{-1} (A(\varepsilon')) <\eta.
\end{aligned}
\]

\textbf{Case 3 ($t_j<s<t$)}. 
Similar to Case 1 we only need to show $\mu_{t,s}(A(\varepsilon'))<\eta$ whose proof turn to be 
similar to the proof in Case 2. 
Indeed, if 
\begin{equation}\label{Equ:mu_tsA_geq_eta}
	\mu_{t,s}(A(\varepsilon'))\geq \eta,
\end{equation}
then
by  \eqref{Equ:SkewMeasure}, \eqref{equ2:mu_ts:Avar_eps_2}
 and  \eqref{equ2:mu_ts:Avar_eps_U2}
\[
\begin{aligned}
	\frac{\eta}{2}
   &>\mu_{t,t_j}(A(\varepsilon'/2))
       = \mu_{t,s}*(\mu_{s,t_j}\circ U(t,s)^{-1})(A(\varepsilon'/2))\\
   & =   \int_\H\int_\H \Eins_{A(\varepsilon'/2)} (x+y)\, \mu_{t,s}(dx) ( \mu_{s,t_j}\circ U(t,s)^{-1})(dy)\\
   &\geq  \mu_{t,s}(A(\varepsilon'))\cdot (\mu_{s,t_j}\circ U(t,s)^{-1})(A(\varepsilon'/2)^c)\\
   &\geq \eta\left(1-\frac{\eta}{2}\right).
\end{aligned} 
\]
This implies $\eta>1$ which contradicts \eqref{Equ:mu_tsA_geq_eta} because $\mu_{t,s}$ is a probability measure.

Combining the three cases discussed above, we obtain \eqref{Equ:UniformStochContinuous:Mod} and hence the proof is complete. 
\end{proof}

\begin{remark}
It is worth to point out that the proof used the fact that $(\mu_{t,s})_{t\geq s}$ is a sequence of probability measures. Another proof given by 
Corollary \ref{cora:inf_div} is also based on this fact. It might be interesting to look at the convolution equation \ref{Equ:SkewMeasure} for general measures.  
\end{remark}

\begin{remark}
Suppose that Assumption \ref{Ass:ContinuousMeasure'} holds and
for every fixed $s_0<t_0$, there exists some constant $c>0$ such that for every $s_0\leq s\leq t\leq t_0$ (cf. \eqref{equ:U_boundedness}), 
\begin{equation*}
                1/c |x|   \leq |U(t,s)x|\leq c|x|,\quad x\in\H, \ s_0\leq s\leq t\leq t_0.
\end{equation*}
The most simple case is the homogeneous convolution equation when $U\equiv I$. 
Then from the proof of Lemma \ref{Lemma:UniformStochContinuous} or directly from \eqref{Equ:UniformStochContinuous:Mod} we obtain that 
$(\mu_{t,s})_{t\geq s}$ is uniformly stochastically continuous on  compact intervals. 
That is, for every $s_0< t_0$ and every $\varepsilon,\eta>0$, there exists a $\delta>0$ such that for all $s,t\in[s_0,t_0]$ 
satisfying
$0\leq t-s<\delta$, we have $$\mu_{t,s}(\{x\in\H\colon |x|>\varepsilon\})<\eta.$$
\end{remark}

Now we are ready to prove the following main result in this subsection. 

\begin{theorem}\label{Thm:mu_ts_IID}
The measures $(\mu_{t,s})_{t\geq s}$ satisfying \eqref{Equ:SkewMeasure} and Assumption \ref{Ass:ContinuousMeasure'} 
 are infinitely divisible. 
\end{theorem}

\begin{proof}
For simplicity we only show that $\mu_{1,0}$ is infinitely divisible since the proof for $\mu_{t,s}$ with arbitrary $t\geq s$ is similar. 
  By \eqref{Equ:SkewMeasure}, we can write for every $m\in\N$, 
  \begin{equation}\label{Equ:convolution-decomp-mu10}
 	 \mu_{1,0}=\Pi_{j=0}^{*(2^m-1)} \mu_{\frac{j+1}{2^m},
	  \frac{j}{2^m}} \circ U\left(1,\frac{j+1}{2^m}\right)^{-1} .
  \end{equation}
  Here $\Pi^*$ denotes the convolution product. Indeed, 
  By \eqref{Equ:SkewMeasure} we have 
  \[
   	 \mu_{1,0}=\mu_{1,1/2}*\left(\mu_{1/2,0}\circ  U\left(1,\frac12\right)^{-1}\right) . 
\]
  So \eqref{Equ:convolution-decomp-mu10} holds for $m=1$. Now we assume that \eqref{Equ:convolution-decomp-mu10} holds for some $m\geq 1$. Then by \eqref{Equ:SkewMeasure} we have 
  for all $j=0,1,\cdots,2^m-1$, 
    \[
  \begin{aligned}
  	  \mu_{\frac{j+1}{2^m},
	         \frac{j}{2^m}} 
   =
   	\mu_{\frac{(2j+1)+1}{2^{m+1}},\frac{2j+1}{2^{m+1}}}
	*\left(\mu_{\frac{2j+1}{2^{m+1}},\frac{2j}{2^{m+1}}}
	     \circ U\left(\frac{(2j+1)+1}{2^{m+1}},\frac{2j+1}{2^{m+1}}\right)^{-1}
	  \right)
  \end{aligned}
  \]
Note that for any probability measures $\mu,\nu$ on $\H$ and measurable map $T$ on \H, it is easy to check that 
\begin{equation}\label{Equ:mu*nu_T_inverse}
	(\mu*\nu)\circ T^{-1} =(\mu\circ T^{-1} ) *(\nu\circ T^{-1} ).
\end{equation}
So we have 
  \[
  \begin{aligned}
  	  &\mu_{\frac{j+1}{2^m},
	         \frac{j}{2^m}} \circ U\left(1,\frac{j+1}{2^m}\right)^{-1} \\
   =&
   	\mu_{\frac{(2j+1)+1}{2^{m+1}},\frac{2j+1}{2^{m+1}}}
	\circ U\left(1,\frac{(2j+1)+1}{2^{m+1}}\right)^{-1} 
	*\left(\mu_{\frac{2j+1}{2^{m+1}},\frac{2j}{2^{m+1}}}
	     \circ U\left(1,\frac{2j+1}{2^{m+1}}\right)^{-1}
	  \right) \\
=&\Pi_{k=2j}^{*(2j+1)} \mu_{\frac{k+1}{2^{m+1}},
	  \frac{k}{2^{m+1}}} \circ U\left(1,\frac{k+1}{2^{m+1}}\right)^{-1} .	  
  \end{aligned}
  \]
Therefore by assumption we have 
  \[
  \begin{aligned}
   	 \mu_{1,0}
=&\Pi_{j=0}^{*(2^m-1)}\Pi_{k=2j}^{*(2j+1)} \mu_{\frac{k+1}{2^{m+1}},
	  \frac{k}{2^{m+1}}} \circ U\left(1,\frac{k+1}{2^{m+1}}\right)^{-1} \\
=&\Pi_{k=0}^{*(2^{m+1}-1)}\mu_{\frac{k+1}{2^{m+1}},
	  \frac{k}{2^{m+1}}} \circ U\left(1,\frac{k+1}{2^{m+1}}\right)^{-1} .	 
\end{aligned}
\]
Now we get that \eqref{Equ:convolution-decomp-mu10} holds for $m\geq 1$. 

By Lemma \ref{Lemma:UniformStochContinuous} we get that $\mu_{1,0}$ is the limit of an infinitesimal triangular array. 
Then by \cite[Corollary VI.6.2]{Par67} we obtain that $\mu_{1,0}$ is   infinitely divisible.
\end{proof}

\subsection{Associated additive processes}\label{Subsec:additive_process}
The following theorem shows that there is a natural additive process associated with the famlily of measures satisfying the skewed convolution equation. 
It could be partially regarded as a generalization of \cite[Theorem 9.7 (ii)]{Sat99} where (non-skewed) convolution equations are studied. The proof follows the  hints given in \cite{Sat99}. That is, it is similar to the proof of  \cite[Theorem 7.10 (ii)]{Sat99} where time homogeneous convolution equations are studied.

\begin{theorem}\label{thm:mu_ts2process}
Fix $t_0\in\R$ and let  $(\mu_{t,s})_{t_0\geq t\geq s}$
be a system of probability measures on $(\H,\B(\H))$ 
such that for all $s\leq r\leq t\leq t_0$ 
 \begin{equation}\label{equ:thm:mu_ts2process:mu_cond1}
   \mu_{t,s} = \mu_{t,r}*(\mu_{r,s}\circ U(t,r)^{-1}), 
 \end{equation}
and 
\begin{equation}\label{equ:thm:mu_ts2process:mu_cond3}
 \begin{aligned}
    \mu_{t,s}&\Rightarrow \delta_0, \quad \textrm{as}\ s\uparrow t,\\
    \mu_{t,s}&\Rightarrow \delta_0, \quad \textrm{as}\ t\downarrow s.
 \end{aligned}
\end{equation}
Set for all $s\leq t\leq t_0$
 $$\tilde{\mu}_{t,s}=\mu_{t,s}\circ U(t_0,t)^{-1}.$$
Then 
\begin{enumerate}
\item  For all $s\leq r\leq t\leq t_0$,
	\begin{equation}\label{equ:thm:mu_ts2process:result1}
	       \tilde{\mu}_{t,s}= \tilde{\mu}_{t,r}* \tilde{\mu}_{r,s}, 
	\end{equation}
 \begin{equation}\label{equ:thm:mu_ts2process:result3}
  \tilde{\mu}_{s,s}=\delta_0  
\end{equation}
and  
\begin{equation}\label{equ:thm:mu_ts2process:result2}
 \begin{aligned}
    \tilde{\mu}_{t,s}&\Rightarrow \delta_0, \quad \textrm{as}\ s\uparrow t,\\
    \tilde{\mu}_{t,s}&\Rightarrow \delta_0, \quad \textrm{as}\ t\downarrow s.
 \end{aligned}
\end{equation} 
\item
There is an stochastic continuous additive process 
$(X_t)_{t_0\geq t}$ associated with $(\tilde{\mu}_{t,s})_{t_0\geq t\geq s}$ in the following sense:
\begin{enumerate}
\item For all $t\leq t_0$, $X_t$ has the distribution $\mu_{t_0,t}$. In particular, $X_{t_0}=0$ almost surely.
\item For all $t_0\geq t_1>t_2>\cdots >t_n$, the increments $X_{t_j}-X_{t_{j-1}}$ with $j=1,2,\cdots,n$ are independent to each other. Moreover, for all $s\leq t\leq t_0$, the increment $X_s-X_t$ has the distribution $\tilde{\mu}_{t,s}$.
\end{enumerate}
\end{enumerate}
\end{theorem}

\begin{proof}
(1)  Let $t_0$ be fixed. 
First we show \eqref{equ:thm:mu_ts2process:result1}. By \eqref{equ:thm:mu_ts2process:mu_cond1} we have 
\[
   \begin{aligned}
	\tilde{\mu}_{t,s}&=\mu_{t,s}\circ U(t_0,t)^{-1} \\
	&=\left(\mu_{t,r} * \left(\mu_{r,s}\circ U(t,r)^{-1}\right)\right)\circ U(t_0,t)^{-1} \\
	&=\left(\mu_{t,r}\circ U(t_0,t)^{-1}\right) *
	\left(\mu_{r,s}\circ U(t_0,r)^{-1}\right) \\
	&=\tilde{\mu}_{t,r}*\tilde{\mu}_{r,s}.
   \end{aligned}
\]
This proves \eqref{equ:thm:mu_ts2process:result1}. Hence 
 we have for all $s\leq t_0$
\[
	\tilde{\mu}_{ss}=	\tilde{\mu}_{ss}*	\tilde{\mu}_{ss}.
\]
As we have noted that the unique idempotent measure  on Hilbert space is the Dirac measure $\delta_0$. Hence 
\eqref{equ:thm:mu_ts2process:result3} follows immediately.

Fix some $s_0<t_0$.  By Lemma \ref{lem:U_boundedness} there exists some constant $c>0$ such that for all $x\in\H$ and $\ s_0\leq s\leq t\leq t_0$, $|U(t,s)x|\leq c|x|$. Hence for any $\varepsilon>0$, as mentioned in Lemma \ref{lem:U_boundedness} we have 
\[
\begin{aligned}
      \tilde{\mu}_{t,s}(\{x\in\H\colon |x|>\varepsilon\})
     &= \mu_{t,s}(\{x\in\H\colon |U(t,s)x|>\varepsilon\}) \\
     & \leq \mu_{t,s}(\{x\in\H\colon |x|>\varepsilon/c\}). 
\end{aligned}
\]
Therefore by Lemma \ref{lmm:Ass:ContinuousMeasure'} we obtain 
\eqref{equ:thm:mu_ts2process:result2} from \eqref{equ:thm:mu_ts2process:mu_cond3}. 
    		

(2)
 Let $\Omega=\H^{(-\infty,t_0]}$ be the collection of all functions 
 $\omega=(\omega(t))_{t\in(-\infty,t_0]}$ from $(-\infty,t_0]$ into $\H$. Let $X_t(\omega)=\omega(t)$, $t\leq t_0$, be the canonical process. Let $\F$ be the $\sigma$-algebra generated by cylinder sets on $\H^{(-\infty,t_0]}$. For any $n\in\N$, 
 $t_0\geq t_1>t_2>\cdots >t_n$,  and for all $B_j\in\B(\H)$, $j=1,2,\cdots,n$, define 
\begin{equation}\label{mu_t0tn}
\begin{aligned}
    &\m_{t_1,t_2\cdots,t_n}(B_1\times B_2\times \cdots\times B_n)\\
   =&\int_{\H}\cdots  \int_{\H}  \Eins_{B_1}(y_1)
               \tilde{\mu}_{t_0,t_1}(dy_1)
 	               \Eins_{B_2}(y_1+y_2) \tilde{\mu}_{t_1,t_2}(dy_2) \\               
               & \quad \quad \quad \quad \quad \quad \times \cdots \times 
	               \Eins_{B_n}(y_1+y_2+\cdots+y_n) \tilde{\mu}_{t_{n-1},t_n}(dy_n).
\end{aligned}
\end{equation}
Then $\m_{t_1,t_2,\cdots,t_n}$ is extended to be a probability measure on $(\H,\B(\H^{\otimes n}))$. Moreover, 
by using \eqref{equ:thm:mu_ts2process:result1}, it is easy to check that the family of probability measures 
$(\m_{t_1,t_2,\cdots,t_n})_{t_0\geq t_1>t_2>\cdots >t_n}$ satisfies the consistency condition.  
We shall only show that for all 
$t_0\geq t_1>t_2>t_3 >t_4$, 
 $B_j\in\B(\H)$, $j=1,2,3,4$ with $B_2=\H$, we have 
$$\m_{t_1,t_2,t_3,t_4}(B_1\times \H\times B_3\times B_4 )=\m_{t_1,t_3,t_4}(B_1 \times B_3 \times B_4)$$ 
below to illustrate the point. In fact, by  \eqref{equ:thm:mu_ts2process:result1} and the definition in  \eqref{mu_t0tn} we have
\[
\begin{aligned}
	       &  \m_{t_1,t_2,t_3,t_4}(B_1\times \H\times B_3 \times B_4) \\
	     =&\int_\H \int_\H
	     			\Eins_{B_1}(y_1)\tilde{\mu}_{t_0,t_1}(dy_1)   
				 \Eins_{\H}(y_1+y_2) \tilde{\mu}_{t_1,t_2}(dy_2) 
				   \\
	     	& \quad\quad	\int_\H\int_\H 					    
					      \Eins_{B_3}(y_1+y_2+y_3) \tilde{\mu}_{t_2,t_3}(dy_3) 
						     \Eins_{B_4 }(y_1+y_2+y_3+y_4)   	\tilde{\mu}_{t_3,t_4}(dy_4)
  					      \\
	     =&\int_\H \int_\H
	     			\Eins_{B_1}(y_1)\tilde{\mu}_{t_0,t_1}(dy_1)   
				\tilde{\mu}_{t_3,t_4}(dy_4)
				   \\
	     	& \quad\quad					      
	    		 \int_\H\int_\H 
					      \Eins_{B_3}(y_1+y_2+y_3)  \Eins_{B_4 }(y_1+y_2+y_3+y_4) \,	 \tilde{\mu}_{t_1,t_2}	(dy_2)  \tilde{\mu}_{t_2,t_3}(dy_3) \\ 		
	     =&\int_\H \int_\H
	     			\Eins_{B_1}(y_1)\tilde{\mu}_{t_0,t_1}(dy_1)   
				\tilde{\mu}_{t_3,t_4}(dy_4)
				   \\
	     	&					      
	    	\quad\quad	 \int_\H 
					      \Eins_{B_3}(y_1+z)  \Eins_{B_4 }(y_1+z+y_4) 
					      \,	 (\tilde{\mu}_{t_1,t_2} * \tilde{\mu}_{t_2,t_3})(dz) \\
             =&\int_\H \int_\H \int_\H 
	     			\Eins_{B_1}(y_1)\tilde{\mu}_{t_0,t_1}(dy_1)   
%
%
					      \Eins_{B_3}(y_1+z) 
					       \tilde{\mu}_{t_1,t_3}(dz) 
					       \Eins_{B_4 }(y_1+z+y_4) 
							\,	\tilde{\mu}_{t_3,t_4}(dy_4)\\
	 %
	     =&\m_{t_1,t_3,t_4}(B_1\times B_3\times B_4 ).					      
\end{aligned}
\]

Therefore, by Kolmogorov's extension theorem we get a unique probability measure $\P$ on $(\Omega,\F)$ such that 
\begin{equation}\label{P2mu_t0tn}
 \P(X_{t_1}\in B_1,X_{t_2}\in B_2, \cdots, X_{t_n}\in B_n)
 =\m_{t_1,t_2,\cdots,t_n}(B_1\times B_2\times \cdots\times B_n).
\end{equation}

Let us show that $(X_t)_{t_0\geq t}$ is a stochastic continuous additive process on $(\Omega,\F,\P)$ in the sense stated in this theorem. 

Note that for any $f\in\BB_b(\H^{\otimes n})$, following from \eqref{mu_t0tn} and \eqref{P2mu_t0tn} we obtain
\begin{equation}\label{Ef}
\begin{aligned}
    &\E[f(X_{t_1},X_{t_2},\cdots,X_{t_n})]\\
   =&\int_{\H}\cdots  \int_{\H} 
               f(y_1,y_1+y_2,\cdots,y_1+y_2\cdots+y_n)\,\tilde{\mu}_{t_0,t_1}(dy_1) \\               
               & \quad \quad \quad \quad \quad \quad \times 
               \tilde{\mu}_{t_1,t_2}(dy_2)\times\cdots\times\tilde{\mu}_{t_{n-1},t_n}(dy_n).
\end{aligned}    
\end{equation}

In particular from \eqref{Ef} we get that $(X_t)_{t\leq t_0}$ is distributed as $\tilde{\mu}_{t_0,t}=\mu_{t_0,t}$. 
Hence $\P(X_{t_0}=0)=1$ since $X_{t_0}\sim \tilde{\mu}_{t_0,t_0}=\mu_{t_0,t_0}=\delta_0$. 

Let $x_0=0$. For fixed $z_1,\cdots,z_n\in\H$, consider
\[
   f(x_1,\cdots,x_n)=\exp\left(i\sum_{j=1}^n\<z_j,x_j-x_{j-1}\>\right),\quad
   x_1,\cdots,x_n\in\H.
\]
It follows from \eqref{Ef} that 
\[
\begin{aligned}
 & \E\left[\exp\left(i\sum_{j=1}^n\left\<z_j,X_{t_j}-X_{t_{j-1}}\right\> \right)\right] \\
=& \int_{\H}\cdots \int_{\H}
\exp\left(i\sum_{j=1}^n\<z_j,y_j\>\right)\tilde{\mu}_{t_0,t_1}(dy_1)\cdots\tilde{\mu}_{t_{n-1},t_n}(dy_n) \\
=&\prod_{j=1}^n \int_{\H} \exp(i\<z_j,y_j\>)\,\tilde{\mu}_{t_{j-1},t_j}(dy_j).
\end{aligned}
\]
This implies for every $j=1,2,\cdots,n$ 
\begin{equation}\label{Xj-Xj-1}
\begin{aligned}
  \E\left[\exp\left(i\<z_j,X_{t_j}-X_{t_{j-1}}\> \right)\right] 
= \int_{\H}
\exp\left(i\<z_j,y_j\>\right)\tilde{\mu}_{t_j,t_{j-1}}(dy_j)
\end{aligned}
\end{equation}
and 
\begin{equation}\label{sumXj-Xj-1}
\begin{aligned}
  \E\left[\exp\left(i\sum_{j=1}^n\<z_j,X_{t_j}-X_{t_{j-1}}\> \right)\right] 
= \prod_{j=1}^n \E\left[\exp\left(i\<z_j,X_{t_j}-X_{t_{j-1}}\> \right)\right]. 
\end{aligned}
\end{equation}
Equation \eqref{Xj-Xj-1} shows that $X_{t_j}-X_{t_{j-1}}$ has the distribution $\tilde{\mu}_{t_j,t_{j-1}}$, while Equation \eqref{sumXj-Xj-1} shows that $(X_t)_{t_0\geq t}$ has independent increments. 

We have shown that for any $t_0\geq t\geq s$, the increment $X_s-X_t$ is  distributed as $\tilde{\mu}_{t,s}$. Therefore from \eqref{equ:thm:mu_ts2process:result2} we get that $X_s-X_t$ converges to 0 in probability as $t$ tends to $s$ or $s$ tends to $t$. This proves that $(X_t)_{t_0\geq t}$ is stochastic continuous. 
\end{proof}

\begin{example}\label{Example:X_t:process:mu_ts} 
For some fixed $t_0$, 
    consider a stochastic process $(X_{t_0,t})_{t\leq t_0}$ over some probability space $(\Omega, \P)$, defined by
    \[
     	X_{t_0,t}:=\int_{t}^{t_0} U(t_0,\sigma)\,dZ_\sigma,
    \]
where $U$ and $Z$ (suppose that $Z$ is stochastic continuous) are the same as in \eqref{NonAutoOUprocess}. Let $\mu_{t,s}$ be the distribution of 
$$X_{t,s}:=\int_{s}^{t} U(t,\sigma)\,dZ_\sigma,\quad t\geq s.
$$
 Then $(\mu_{t,s})_{t_0\geq t\geq s}$ fulfills 
conditions \eqref{equ:thm:mu_ts2process:mu_cond1}
and \eqref{equ:thm:mu_ts2process:mu_cond3}. For all $t_0\geq t> s$, the increment $X_{t_0,s}-X_{t_0,t}$ is distributed as
$$\tilde{\mu}_{t,s}=\mu_{t,s}\circ U(t_0,t)^{-1}.$$ 
The proof is direct.  First we note that 
\[
X_{t_0,s}-X_{t_0,t}=\int_{s}^{t_0} U(t_0,\sigma)\,dZ_\sigma-\int_{t}^{t_0} U(t_0,\sigma)\,dZ_\sigma=U(t_0,t)\int_s^t U(t,\sigma)\,dZ_\sigma. 
\]
So for all $B\in\B(\H)$ we have 
\[
\begin{aligned}
\P(X_{t_0,s}-X_{t_0,t}\in B)&=\P\left(\int_s^t U(t,\sigma)\,dZ_\sigma\in U(t_0,t)^{-1}B\right) \\
                           &=\mu_{t,s}(U(t_0,t)^{-1}B)=\tilde{\mu}_{t,s}(B).
\end{aligned}
\]

Moreover, it is easy to check that the increments of $(X_{t_0,t})_{t_0\geq t}$ are independent to each other. Hence following from 
$$X_{t_0,s}-X_{t_0,t}=(X_{t_0,s}-X_{t_0,r})+(X_{t_0,r}-X_{t_0,t}), \quad
s\leq r\leq t\leq t_0,$$
we obtain  \eqref{equ:thm:mu_ts2process:result1}. 
Finally, \eqref{equ:thm:mu_ts2process:result3} and \eqref{equ:thm:mu_ts2process:result2} are obvious.
\end{example}

Based on \eqref{equ:thm:mu_ts2process:result1} and \eqref{equ:thm:mu_ts2process:result2}, as the proof of  Lemma 
\ref{Lemma:UniformStochContinuous}, we can prove the following result:
on every compact interval $[s_0,t_0]$ with $s_0\leq t_0$, for all
$\varepsilon,\eta>0$, there is a constant $\delta>0$ such that for all $s_0\leq s< t\leq t_0$ satisfying  $|s-t|\leq \delta$, we have 
\begin{equation*}
   \tilde{\mu}_{t,s}(\{x\in\H\colon |x|>\varepsilon\})<\eta. 
\end{equation*}
According to the definition of $\tilde{\mu}_{t,s}$, it is exactly 
\eqref{Equ:UniformStochContinuous:Mod}. 
This is a slightly simpler proof of  Lemma \ref{Lemma:UniformStochContinuous}.
However,  using the stochastic process constructed in Theorem \ref{thm:mu_ts2process}, we have a more transparent 
probabilistic proof.

First let us recall the following lemma which has been proved in \cite[Lemma 9.6, Page 51]{Sat99} for finite dimensional case.  We include the same  proof here for completeness. Note that we  have used similar argument in the proof of 
Lemma \ref{Lemma:UniformStochContinuous}.

\begin{lemma}\label{lmm:unif_cont:stoch_conti_process}
 A stochastically continuous process $(X_t)_{t\in\R}$ taking values in \H is uniformly stochastically continuous on any finite interval $[s_0,t_0]$. That is, for every $\varepsilon>0$ and $\eta>0$, there is $\delta>0$ such that for any $s$ and $t$ in $[s_0,t_0]$ satisfying $|t-s|<\delta$, we have 
  $$\P(|X_t-X_s|>\varepsilon)<\eta.$$
\end{lemma}

\begin{proof}
Let $\varepsilon>0$ and $\eta>0$. Since $(X_t)_{t\in\R}$ is stochastic continuous, we obtain for every $t\in\R$ there is a constant $\delta_t>0$ such that for all $t,s\in [s_0,t_0]$ satisfying  $|t-s|<\delta_t$, we have 
\[\P(|X_t-X_s|>\varepsilon/2)<\eta/2.\] 
It is clear that  $\{I_t:=(t-\delta_t/2,t+\delta_t/2)\colon t\in[s_0,t_0]  \}$ covers the compact interval $[s_0,t_0]$. Then there is a sub-covering $\{I_{t_j}\colon j=1,2,\cdots,n\}$ for $[s_0,t_0]$. Let 
\[
	\delta:=\min\{\delta_{t_j}/2\colon j=1,2,\cdots,n\}. 
\]
For all $t\in[s_0,t_0]$, we have $t\in I_{t_j}$ for some $j=1,2,\cdots,n$. Thus for all  
$s\in[s_0,t_0]$ such that $0\leq |t-s| <\delta$, we get $|s-t_j|<\delta_j$ since 
\[
	|s-t_j|\leq |s-t| + |t-t_j| <\delta +\delta_{t_j}/2\leq \delta_{t_j}.
\]
Hence we obtain
\[
\begin{aligned}
	\P(|X_t-X_s|>\varepsilon)& \leq \P(|X_t-X_{t_j}|>\varepsilon/2) 
								+ \P(|X_s-X_{t_j}|>\varepsilon/2) \\
					     &< \eta/2+\eta/2=\eta.
\end{aligned}
\]
\end{proof}

\begin{cora}\label{cora:inf_div}
 In the framework of Theorem \ref{thm:mu_ts2process} the
estimate 
\eqref{Equ:UniformStochContinuous:Mod} holds and hence $\mu_{t,s}$ is infinitely divisible.
\end{cora}

\begin{proof}
The proof is the same as \cite[Theorem 9.1, Page 51]{Sat99}. 
By Theorem \ref{thm:mu_ts2process} 
there is a stochastic continuous additive process 
$(X_t)_{t_0\geq t\geq s_0}$ such that $\tilde{\mu}_{t,s}$ is the distribution of the increment $X_s-X_t$ for all $t_0\geq t\geq s\geq s_0$. 
By Lemma \ref{lmm:unif_cont:stoch_conti_process} we obtain that $(X_t)_{t_0\geq t\geq s_0}$ is uniformly stochastic continuous. 
This means that for every $\varepsilon>0$ and $\eta>0$, there is a $\delta>0$ such that for all $s,t\in[s_0,t_0]$ satisfying $|t-s|<\delta$, we have 
\[
\P(|X_s-X_t|> \varepsilon)<\eta.
\]
In other words $$\tilde{\mu}_{t,s}(|x|>\varepsilon)<\eta.$$
 This proves \eqref{Equ:UniformStochContinuous:Mod} since $\tilde{\mu}_{t,s}=\mu_{t,s}\circ U(t_0,t)^{-1}$ by definition. The remainder of the proof is the same with the proof of Theorem \ref{Thm:mu_ts_IID}. 
%
\end{proof}

Using the stochastic process constructed in Theorem \ref{thm:mu_ts2process}, we have also another proof for the weak continuity of $s\mapsto \mu_{t,s}$ with $s\leq t$ shown in the first assertion of Theorem 
\ref{thm:mu_ts:continuity}. 


\begin{proof}[Another proof of Theorem 
\ref{thm:mu_ts:continuity} (i)]
For every fixed $t\in\R$, we have shown in Theorem 
\ref{thm:mu_ts2process} 
that there is a stochastic continuous process 
$(X_s)_{t\geq s}$ such that $\mu_{t,s}=\tilde{\mu}_{t,s}$ is the distribution of the increments $X_s-X_t=X_s$ for all $s\leq t$. 
Since $(X_s)_{t\geq s}$ is stochastic continuous, 
$X_{s'}$ converges to $X_s$ in probability as $s'$ tends to $s$  
for all $s,s'\leq t$. 
This implies that the distribution of 
$X_{s'}$ converges to the distribution of $X_s$ weakly. That is, 
$\mu_{t,s'}\Rightarrow \mu_{t,s}$  as $s'\to s$. The proof is complete. 
\end{proof}

\begin{remark}
It is standard to show that  there exists a Markov process $(X_{t,s}^x)_{t\geq s}$ for all $x\in\H$ with $X_{s,s}^x=x$ associated with \((p_{s,t})_{t\geq s}\) by Kolmogorov's consistency theorem. 
One  may mimic the idea in  \cite{BRS96} and \cite{FR00} to construct corresponding Markov processes with c\`adl\`ag paths and even show that the process solves the stochastic equation in certain sense. This will be done separately. 
Here we mention that according to \eqref{equ:mu_ts_2_p_st}, for every $t\geq s$, the random variable $\tilde{X}_{t,s}:=X_{t,s}^x-U(t,s)x$ is  distributed as $\mu_{t,s}$ (independent of $x$). Indeed, let $\Q$ be the law of the process $(X_{t,s}^x)_{t\geq s}$. Then for all $B\in\B(\H)$, we have 
\[
	\Q(\tilde{X}_{t,s}\in B)=\Q(X_{t,s}^x\in B+U(t,s)x) = p_{s,t}(x,B+U(t,s)x) = \mu_{t,s}(B). 
\]
So the process $(\tilde{X}_{t_0,t})_{t_0\geq t}$ is one version of the process $(X_t)_{t_0\geq t}$ constructed in Theorem \ref{thm:mu_ts2process}. If both processes $(\tilde{X}_{t_0,t})_{t_0\geq t}$ and $(X_t)_{t_0\geq t}$ are right continuous, then they are indistinguishable. In this note we try to understand the process $(X_t)_{t_0\geq t}$ via the method of spectral representation (Ref. Remark \ref{Remark:convolution-integral-representation}).
\end{remark}

\subsection{L\'evy-Khintchine representation}\label{sub:sec:LKrepr}
Now we assume that for every $t\geq s$, the measure $\mu_{t,s}$  is infinitely divisible. 
Then by the L\'evy-Khintchine theorem \cite[Theorem VI.4.10, Page 182]{Par67},   there exists a negative definite, Sazonov continuous function $\psi_{t,s}$ on $\H$
such that 
\[
\hat{\mu}_{t,s}(\xi)=\exp(-\psi_{t,s}(\xi)),\quad 
\xi\in\H,
\]
and $\psi_{t,s}$ has the following form
\begin{equation}\label{Lambda_ts}
   \begin{aligned}
   \psi_{t,s}(\xi)=-i\<a_{t,s},\xi\>&+\frac12\<\xi,R_{t,s}\xi\>\\
		&{}-\int_\H \left(\e^{i\<\xi,x\>} 
		-1-\frac{i\<\xi,x\>}{1+|x|^2} \right)\,\m_{t,s}(dx),\quad 
		\xi\in\H, 
\end{aligned}		
\end{equation}
where \(a_{t,s}\in\H\), 
 \(R_{t,s}\) is a non-negative definite, symmetric trace class operator on \H,  
and \(\m_{t,s}\) is a L\'evy measure on \H.  That is, 
$$\mu_{t,s}=[a_{t,s}, R_{t,s}, \m_{t,s}],\quad t\geq s.$$ 

In terms of the  characteristic exponent $\psi_{t,s}$ of $\mu_{t,s}$, condition \eqref{Equ:SkewMeasure} is equivalent to  
\begin{equation}\label{Equ:SkewMeasure:Symbol}
	\psi_{t,s}(\xi)=\psi_{t,r}(\xi)+\psi_{r,s}(U(t,r)^*\xi),\quad \xi\in\H
\end{equation}
for every $t\geq r\geq s$.

According to \eqref{Lambda_ts} the right hand side of \eqref{Equ:SkewMeasure:Symbol} is given by
\[
\begin{aligned}
& \psi_{t,r}(\xi)+\psi_{r,s}(U(t,r)^*\xi) \\
 = &  -i\<a_{t,r},\xi\>+\frac12\<\xi,R_{t,r}\xi\>
		-\int_\H \left(\e^{i\<\xi,x\>}
		-1-\frac{i\<\xi,x\>}{1+|x|^2} \right)\,\m_{t,r}(dx) \\
   &	- i\<U(t,r)a_{r,s},\xi\>+\frac12\<\xi,U(t,r)R_{r,s}U(t,r)^*\xi\> \\
   &	\phantom{-i\<a_{t,s},\xi\>+\frac12\<\xi,R_{t,s}\xi\>}
            -\int_\H \left(\e^{i\<\xi,U(t,r)x\>} 
    	      -1-\frac{i\<\xi,U(t,r)x\>}{1+|x|^2} \right)\,\m_{r,s}(dx) \\
 =  &	 -i\<a_{t,r}+U(t,r)a_{r,s},\xi\>+\frac12\<\xi,(R_{t,r}+U(t,r)R_{r,s}U(t,r)^*)\xi\> 
                   \\
    &   -\int_\H \left(\e^{i\<\xi,x\>}
		-1-\frac{i\<\xi,x\>}{1+|x|^2} \right)\,
		(\m_{t,r}+ \m_{r,s}\circ U(t,r)^{-1} )
		         (dx) \\ 
    &  +\int_\H \frac{i\<\xi,U(t,r)x\>}{1+|x|^2}\,\m_{r,s}(dx)
       -\int_\H \frac{i\<\xi,U(t,r)x\>}{1+|U(t,r)x|^2}\,\m_{r,s}(dx).
\end{aligned}
\]

Therefore, by the uniqueness of the canonical representation for infinitely divisible distributions we have the following identities (cf. also the proof of  \cite[Corollary 1.4.11]{Ouyang09})
\begin{equation}\label{Equ:Structure_char_triplet}
\begin{aligned}
a_{t,s}&=a_{t,r}+U(t,r)a_{r,s}\\
&\quad\quad\quad +\int_{\H} U(t,r)x\left( 
		\frac{1}{1+|U(t,r)x|^2}-\frac{1}{1+|x|^2}
	 \right)\m_{r,s}(dx), \\
R_{t,s}&=R_{t,r}+U(t,r)R_{r,s}U(t,r)^*,\\
\m_{t,s}&=\m_{t,r}+\m_{r,s}\circ U(t,r)^{-1}
\end{aligned}
\end{equation}
for every $t\geq r\geq s$. 

In particular, from \eqref{Equ:Structure_char_triplet} (or directly from \eqref{Equ:mu_tt=delta0}) we have 
\begin{equation}\label{Equ:a_ttR_ttM_tt=0}
	a_{t,t}=0, \quad R_{t,t}=0,\quad \m_{t,t}=0
\end{equation}
for all $t\in\R$.

According to \cite[Theorem VI.5.5, Page 189]{Par67} (see also \cite[Theorem 8.7, Page 41]{Sat99} for the finite dimensional version), 
the following results concerning the continuity of the generating triplets follows from the weak continuity of the $\mu_{t,s}$ in $s$ over $(-\infty,t]$ proved in Theorem \ref{thm:mu_ts:continuity}.

\begin{prop}\label{char_functional_continuity}
Suppose that Assumption \ref{Ass:ContinuousMeasure'} holds.
Then for every $t\in\R$, we have  the following assertions. 
\begin{enumerate}
\item The map $s\mapsto a_{t,s}$ on $(-\infty,t]$ is continuous in \H;
\item The map $s\mapsto \m_{t,s}$ on $(-\infty,t]$ is weakly continuous outside every closed neighborhood of the origin. That is, for every $B\in\B(\R^d)$ with $B\subset \{x\colon |x|>\varepsilon \}$ for some $\varepsilon>0$, $s\mapsto \m_{t,s}(B)$ is continuous. 
\item Let  $(s_n)_{n\geq 1}$ be a sequence of real numbers such that $s_n\leq t$ and $s_n\to s$ as $n\to \infty$. For all $\xi\in\H$ define a trace-class operator $T_n$ by
\[
        \<T_n\xi,\xi \>=\<\xi,R_{t,s_n}\xi\> 
             + \int_{|x|\leq 1} \<x,\xi\>\,\m_{t,s_n}(dx).
\]
Then $\{T_n\}_{n\geq 1}$ is compact and 
         \begin{equation}\label{equ:cont_levy_measure}
              \lim_{\varepsilon\to 0} \varlimsup_{n\to\infty} 
          \left|
              \int_{|x|\leq \varepsilon } 	      
	    \<x,\xi\>^2\,\m_{t,s_n}(dx) +\<\xi,R_{t,s_n}\xi\>
	 - \<\xi,R_{t,s} \xi\>                         \right| =0.
         \end{equation}         
\end{enumerate}
\end{prop}

Using the structure information of  \(\mu_{t,s}=[a_{t,s},R_{t,s},\m_{t,s}] \)
 in \eqref{Equ:Structure_char_triplet}, Proposition \ref{char_functional_continuity} implies the following result. 

\begin{cora}\label{cora:continuity_R_t,s}
Suppose that Assumption \ref{Ass:ContinuousMeasure'} holds. 
For every $\xi\in\H$ and $t\in\R$, 
the map $s\mapsto \<\xi,R_{t,s}\xi\>$ with $s\leq t$ is continuous.  
\end{cora}

\begin{proof}
Let $(s_n)_{n\geq 1}$ with $s_n\leq t$ for all $n\geq 1$ be a sequence of real numbers converging to $s$. By \eqref{Equ:Structure_char_triplet} we have 
\[
         \m_{t,s_n}\leq \m_{t,s'},
\]
where $s'$ is any real number such that $s'\leq s_n$ for all $n\geq 1$. 
Hence 
         \[
             0\leq  \lim_{\varepsilon\to 0} \varlimsup_{n\to\infty}
                	 \int_{|x|\leq \varepsilon } \<x,\xi\>^2\,\m_{t,s_n}(xd)
	 \leq \lim_{\varepsilon\to 0} 
                	 \int_{|x|\leq \varepsilon } \<x,\xi\>^2\,\m_{t,s'}(dx) =0.
         \]
It implies  
\[
      \lim_{\varepsilon\to 0} \varlimsup_{n\to\infty}
                	 \int_{|x|\leq \varepsilon } \<x,\xi\>^2\,\m_{t,s_n}(dx) =0.    
\]         
From \eqref{equ:cont_levy_measure} we get 
\[
 \begin{aligned}
  & \varlimsup_{n\to\infty}
        \left| \<\xi,R_{t,s_n}\xi\> -  \<\xi,R_{t,s}\xi\> \right|  \\
=&  \varlimsup_{n\to\infty}
    \left| \<\xi,R_{t,s_n}\xi\> -  \<\xi,R_{t,s}\xi\> \right| 
     -  \lim_{\varepsilon\to 0} \varlimsup_{n\to\infty}
                	 \int_{|x|\leq \varepsilon } \<x,\xi\>^2\,d\m_{t,s_n}(x) \\
\leq &
       \lim_{\varepsilon\to 0} \varlimsup_{n\to\infty}
       \left|
                	 \int_{|x|\leq \varepsilon } \<x,\xi\>^2\,d\m_{t,s_n}(x) + \<\xi,R_{t,s_n}\xi\> -  \<\xi,R_{t,s}\xi\>  
 	\right|  \\
=&0.       	 
\end{aligned}	 
\]
This proves 
\[
   \lim_{n\to \infty}   \<\xi,R_{t,s_n}\xi\> =  \<\xi,R_{t,s}\xi\>.
\]
So the proof is finished.
\end{proof}

%

\subsection{Constructions}\label{subsec:ConsRep1}
The following proposition shows a typical form of the family of probability measures $(\mu_{t,s})_{t\geq s}$ satisfying the skew convolution equation \eqref{Equ:SkewMeasure}.

\begin{prop}\label{Prop:Lambda_ts}
Assume that $\lambda_{t,s}$ for every $t>s$ is the characteristic exponent of some infinitely divisible probability measure on $\H$.  Let $\pi$ be  an atomless $\sigma$-finite measure on $\R$. Suppose that for every $\xi\in\H$ the function  $s\mapsto \lambda_{t,s}(\xi)$ with $t>s$ is locally integrable with respect to $\pi$. Then 
\begin{enumerate}
 \item For every $t>s$,  the following characteristic functional 
 \begin{equation}\label{Char_rep}
\hat{\mu}_{t,s}(\xi)=\exp\left(-\int_s^t \lambda_{t,\sigma}(\xi)\,\pi(d\sigma) \right),\quad \xi\in\H, 
\end{equation}
defines a family of infinitely divisible probability measures $\mu_{t,s}$ on \H.
\item In addition, if for every $t>r>s$, 
\begin{equation}\label{Equ:Lambda_ts}
    \lambda_{t,s}(\xi)=\lambda_{r,s}(U(t,r)^*\xi),\quad \xi\in\H,
\end{equation}
then  \eqref{Equ:SkewMeasure} holds.
\end{enumerate}
\end{prop}

\begin{proof} 
(1)  Since for every $t>s$, $\lambda_{t,s}$ is the characteristic exponent of an
infinitely divisible probability measure on $\H$, there exists an element 
$a^{\lambda}_{t,s}$ in \H, a non-negative definite, self-adjoint trace class operator $R^{\lambda}_{t,s}$ on 
\H, and a L\'evy measure $\m^{\lambda}_{t,s}$ on \H such that
\begin{equation}\label{Equ:lambda_ts_aRm} 
\begin{aligned}
   \lambda_{t,s}(\xi)=-i\<a^{\lambda}_{t,s},\xi\>&+\frac12\<\xi,R^{\lambda}_{t,s}\xi\> \\
		&{}-\int_\H \left(\e^{i\<\xi,x\>}
		-1-\frac{i\<\xi,x\>}{1+|x|^2} \right)\,\m^{\lambda}_{t,s}(dx),\quad \xi\in\H.
\end{aligned}
\end{equation}

It is clear that $\int_s^t \lambda_{t,\sigma}(\xi)\,\pi(d\sigma)$ represents the characteristic exponent 
of the infinitely divisible measure with characteristic triplet
\begin{equation}\label{Equ:int_mu_ts_lambda_aRm}
\left[
\int_s^t a^{\lambda}_{t,\sigma} \pi(d\sigma),
\int_s^t R^{\lambda}_{t,\sigma} \pi(d\sigma),
\int_s^t  \m^{\lambda}_{t,\sigma}  \pi(d\sigma)
\right].
\end{equation}
This proves that \eqref{Char_rep} defines a family of infinitely divisible probability measures 
$(\mu_{t,s})_{t\geq s}$ on $(\H,\B(\H))$. 

(2) 
For every $\xi\in\H$, $t> r> s$, by \eqref{Equ:Lambda_ts} we have 
\[
\begin{aligned}
   &\hat{\mu}_{t,r}(\xi)\hat{\mu}_{r,s}\bigl(U(t,r)^*\xi\bigr) \\
  =& \exp\left(-\int_r^t \lambda_{t,\sigma}(\xi)\,\pi( d\sigma ) -\int_s^r \lambda_{r,\sigma}(U(t,r)^*\xi)\, \pi( d\sigma )\, \right) \\
  =& \exp\left(-\int_r^t \lambda_{t,\sigma}(\xi)\,\pi( d\sigma ) -\int_s^r \lambda_{t,\sigma}(\xi)\, \pi( d\sigma )\, \right) \\
  =&\exp\left(-\int_s^t \lambda_{t,\sigma}(\xi)\,\pi( d\sigma )\right)\\
  =& \hat{\mu}_{t,s}(\xi).
\end{aligned}
\]
This proves \eqref{Equ:SkewMeasure}. 
\end{proof}

\begin{cora}
For every $s\in\R$, assume that $\lambda_{s}$ 
be the characteristic exponent of some infinitely divisible probability measure on $\H$.
Then for all $t,s \in\R$ with $t\geq s$, 
\begin{equation}\label{Equ:Lambda_s:closed_el}
  \lambda_{t,s}(\xi):=\lambda_s\bigl(U(t,s)^*\xi\bigr)
\end{equation}
is the characteristic exponent of some infinitely divisible probability measure on $\H$. Moreover,   $(\lambda_{t,s})_{t\geq s}$
satisfies \eqref{Equ:Lambda_ts} and hence 
\begin{equation}\label{Char_rep_closed}
\hat{\mu}_{t,s}(\xi)=\exp\left(-\int_s^t \lambda_\sigma(U(t,\sigma)^*\xi)\,\pi(d\sigma) \right),\quad \xi\in\H, \ t>s
\end{equation}
defines a family of  infinite divisible measures $(\mu_{t,s})_{t\geq s}$ on $(\H,\B(\H))$ satisfying  \eqref{Equ:SkewMeasure}.
\end{cora}

\begin{proof} 
By the L\'evy-Khintchine formula,  for every $s\in \R$, the characteristic exponent
 \(\lambda_s\) is of the following form
\begin{equation}\label{Lambda_r}
\begin{aligned}
   \lambda_s(\xi)=-i\<a^{\lambda}_s,\xi\>&{}+\frac12\<\xi,R^{\lambda}_s\xi\> \\
	    & {}	-\int_\H \left(\e^{i\<\xi,x\>}
		-1-\frac{i\<\xi,x\>}{1+|x|^2} \right)\,\m^{\lambda}_s(dx),\quad \xi\in\H,
\end{aligned}		
\end{equation}
where \(a^{\lambda}_r\in\H\), \(R^{\lambda}_r\) is a non-negative definite, self-adjoint trace class operator on \H, and \(\m^{\lambda}_r\)
is a L\'evy measure on \H. 

Therefore for all $t\geq s$, $\lambda_{t,s}(\xi)$
has the form \eqref{Equ:lambda_ts_aRm}
with characteristic triplet $a^{\lambda}_{t,s}, R^{\lambda}_{t,s}$ and 
$\m^{\lambda}_{t,s}$ given respectively by
\[
 \begin{aligned}
	& a^{\lambda}_{t,s}= U(t,s) a^{\lambda}_{s}  
 		  + \int_\H U(t,s)x\left[ 
				\frac{1}{1+|U(t,s)x|^2}-\frac{1}{1+|x|^2}
			  \right]\,\m^{\lambda}_r(dx)	\\
        & R^{\lambda}_{t,s}= U(t,s)R^{\lambda}_{s}U(t,s)^*, \\
&	        	\m^{\lambda}_{t,s}(\{0\})=0,\\
&	\m^{\lambda}_{t,s}(A)= \m^{\lambda}_s(U(t,s)^{-1}(A)),\quad A\in\B(\H\setminus \{0\}).
 \end{aligned}
\]

Therefore it is clear that \eqref{Char_rep_closed} defines an infinitely divisible probability measure with characteristic triplet given as \eqref{Equ:int_mu_ts_lambda_aRm}.

The relation \eqref{Equ:Lambda_ts} is simple. Indeed  
for every $t> r> s$,  
\[
\begin{aligned}
\lambda_{r,s}\bigl(U(t,r)^*\xi\bigr)&=\lambda_s\bigl(U(r,s)^*U(t,r)^*\xi\bigr)\\&=\lambda_s\bigl((U(t,r)U(r,s))^*\xi\bigr)
=\lambda_s\bigl(U(t,s)^*\xi\bigr)
=\lambda_{t,s}(\xi).
\end{aligned}
\]
 Consequently, 
by the second assertion of Proposition \ref{Prop:Lambda_ts}, we have  \eqref{Equ:SkewMeasure}. 
\end{proof}

\begin{remark}\label{Rem:EntranceLaw}
For all $t>s$, let $\nu_{t,s}$ be the measure on $(\H,\B(\H))$ with Fourier transformation 
$\hat{\nu}_{t,s}(\xi)=\exp(-\lambda_{t,s}(\xi))$ for every $\xi\in\H$.
Fix $s_0$ and set $\nu_t=\nu_{t,s_0}$ for all $t>s_0$. 
Then  
 \eqref{Equ:Lambda_ts} implies that 
 \[
	\nu_t=\nu_r\circ U(t,r)^{-1},\quad t\geq r>s_0.
 \]
Define a transition semigroup $u_{s,t}$ by $u_{s,t}(x,\cdot)=\delta_{U(t,s)x}(\cdot)$ for every $x\in\H$. 
It is easy to see that 
\[
	\int_\H u_{s,t}f(x)\,\nu_s(dx)=\int_\H f(x) \nu_t(dx),
	\quad f\in B_b(\H).
\] 
That is to say,  $(\nu_t)_{t > s_0}$ 
is an entrance law (or an evolution system of measures studied Section \ref{Sec:Evolution_Measures}) for $u_{s,t}$. While the case \eqref{Equ:Lambda_s:closed_el} happens if the setting for $(\lambda_{t,s})_{t>s}$ above extends to $(\lambda_{t,s})_{t\geq s}$. 
Hence $(\nu_t)_{t > s_0}$  could be extended to $(\nu_t)_{t \geq s_0}$.  
In this case, we say the entrance law is closed. 
\end{remark}

\begin{example}\label{Exa:represation:factoring}
Suppose that the additive process $(Z_t)_{t\in\R}$ admits a factoring $((\lambda_s)_{s\in\R}, \pi)$. 
Assume that for every $t\geq s$,  the integral
 $$X_{t,s}=\int_s^t U(t,u)\,dZ_u$$ 
is well defined (see \cite{Det83,Sat06}), then the characteristic functional of  $X_{t,s}$ 
is given by \eqref{Char_rep_closed}.
\end{example}

\subsection{Spectral representation}\label{subsec:ConsRep2}
Now we consider the spectral representation of $\mu_{t,s}$ for each $t\geq s$. We shall always assume that Assumption \ref{Ass:ContinuousMeasure'} holds. Then we have  for all $t\geq s$,  $\mu_{t,s}$ is infinitely divisible and we can set $\mu_{t,s}=[a_{t,s},R_{t,s},\m_{t,s}]$. Moreover, we have the conclusions in Proposition \ref{char_functional_continuity} and Corollary \ref{cora:continuity_R_t,s} concerning the continuity of $a_{t,s},R_{t,s}$ and $\m_{t,s}$ in $s$. 

The main idea is described as follows. 
For every $t\geq s$, $\xi\in\H$, let 
$$\psi_{t,s}(\xi)=-\log \hat{\mu}_{t,s}(\xi).$$
It takes values in the complex plane. 
For every $t\in\R$, 
\[
       F^t_\xi((s,r])=\psi_{r,s}(U(t,r)^*\xi), \quad s\leq r\leq t
\]
defines an additive complex set function on $\A_t^0$. Here $\A_t^0$ is 
the algebra generated by intervals of the form $(s,r]\cap (-\infty, t]$ with $s\leq r\leq t$.
Assume that $(\mu_{t,s})_{t\geq s}$ is natural (ref. Definition \ref{Def:natural_mu_ts}). Then the additive set function $F^t_\xi$ can be extended to be a complex measure 
$F_\xi^t$ on $((-\infty, t], \B((-\infty,t]))$. If in addition there is a $\sigma$-finite measure $\pi$ on 
 $(\R,\B(\R))$ such that
$
     F_\xi^t\ll \pi, 
$
then we get a spectral representation 
\[
\hat{\mu}_{t,s}(\xi)=\exp\left(-\int_s^t \frac{d F^t_\xi}{d\pi}\, d\pi \right),\quad \xi\in\H
\]
for every $s\leq t$.

To make the choice of $\pi$ explicit and the proof strict, we turn to look at the characteristics of $\mu_{t,s}$,  i.e. the Gaussian part $R_{t,s}$, the jump part $\mu_{t,s}$ and the drift part $a_{t,s}$ in detail. In this way, a canonical control measure $\pi$ will be determined naturally. 

Hence the theoretical result of this section can be used to convince the reader that under some conditions, the distributions of stochastic convolution integrals in Example \ref{Exa:represation:factoring} represents  in fact all $(\mu_{t,s})_{t\geq s}$ satisfying \eqref{Equ:SkewMeasure}. We refer to Remark \ref{Remark:convolution-integral-representation} for more discussion. 

\begin{remark}Instead of studying the space inhomogeneous equation \eqref{Equ:SkewMeasure} directly, 
there is another approach. 
Let $t_0\in\R$ be fixed. First we look at the space-homogeneous equation
 \eqref{equ:thm:mu_ts2process:result1}, if $(\tilde{\mu}_{t_0,s})_{t_0\geq s}$ is natural, then it is easy to get a spectral representation for 
 $(\tilde{\mu}_{t,s})_{ t_0\geq t\geq s}$ using similar method in Sato \cite{Sat06}. Therefore we obtain immediately a spectral representation for $(\mu_{t_0,s})_{t_0\geq s}$ since $\mu_{t_0,s}=\tilde{\mu}_{t_0,s}$. 
 That is, there exists a $\sigma$-finite measure $\pi_{t_0}$ on 
 $((-\infty,t_0], \B((-\infty,t_0]))$, and a family of characteristic exponent $(\lambda_{t_0,\sigma})_{ t_0\geq \sigma}$
of infinite divisible probability measures on $(\H,\B(\H))$ such that for all $s\leq t_0$, 
\[
	\hat{\mu}_{t_0,s}(\xi)=\exp\left(- \int_{s}^{t_0}  \lambda_{t_0,\sigma}(\xi)\,\pi_{t_0}(d\sigma)\right),\quad \xi\in\H.
\]
\end{remark}

Let us fix some notations. 
For every set $O\subset \R$, let $\B(O)$, $\B^0(O)$ denote the space of all Borel subsets and bounded Borel subsets of $O$ respectively. 
Let $\{e_i\}_{i=1}^\infty$ be an orthonormal basis of \H. 
Let $Q:=\{q_i\colon i=1,2,\cdots\}$ be a countable dense subset of $\R$. 

\subsubsection{Gaussian part}\label{Spectral_Rep:GaussianPart}

\begin{lemma}\label{lmm:gauss:scattered:measure}
Assume that Assumption \ref{Ass:ContinuousMeasure'} holds. For every $t\in\R$, there is a function 
\[
R(t,\cdot,\cdot)\colon 
    (\xi,B)\mapsto R(t,\xi,B)\in [0,+\infty],\quad \xi\in\H, \ B\in \B((-\infty,t])
\]
such that 
\begin{enumerate}
	\item For every $\xi\in \H$, $R(t,\xi, \cdot)$ is a $\sigma$-finite measure on $((-\infty, t], \B\left((-\infty, t])\right)$.
	\item For every $B\in \B^0((-\infty, t])$, there is a non-negative definite, self-adjoint trace class operator $R(t,B)$ on \H
	         such that \[R(t,\xi, B)=\<R(t,B)\xi,\xi\>,\quad \xi\in\H;\]
	\item For every $\xi\in\H$ and $s\leq r\leq t$, 
		\begin{equation}\label{R_t:(s,r]}
			R(t,\xi, (s,r])=\<R_{r,s}U(t,r)^*\xi,U(t,r)^*\xi\>
		\end{equation}
		and $R(t,\xi, \cdot)$ is atomless
		\begin{equation}\label{R_t:ss}
			R(t,\xi, \{s\})=0.
		\end{equation}
		In particular, taking $r=t$ in \eqref{R_t:(s,r]} we get
		\begin{equation}\label{R_t:(s,t]}
			R(t,\xi, (s,t])=\<R_{t,s}\xi,\xi\>.
		\end{equation}	
	\item For every $\xi\in\H$, $t\geq s$,  $B\in \B((-\infty, s])$, 
		\begin{equation}\label{R_t:R_s}
			R(t,\xi,B)=R(s,U(t,s)^*\xi,B).
		\end{equation}
	That is
	\begin{equation}\label{R_t:R_s'}
		R(t,B)=U(t,s)R(s,B)U(t,s)^*.
	\end{equation}	
\end{enumerate}
\end{lemma}

\begin{proof}
   First we define a finite additive set function 
   $R(t,\xi,\cdot)$ on $\A_t^0$ for every $\xi\in\H$ and $t\in\R$ via
\eqref{R_t:(s,r]} for all  $s\leq r\leq t$. 

   The additive property $R(t,\xi,\cdot)$ on $\A_t^0$ follows easily from the following identity 
   \begin{equation}\label{R_s,r=R_t,s-R_t,r_xi}
          R(t,\xi, (s,r]) = \<R_{t,s}\xi,\xi\> - \<R_{t,r}\xi,\xi\>
   \end{equation}
   which comes from \eqref{Equ:Structure_char_triplet}.   
   By Corollary \ref{cora:continuity_R_t,s}, $s\mapsto \<R_{t,s}\xi,\xi\>$ is continuous over $(-\infty,t]$. 
   Moreover, since $R(t,\xi,(s,t])=\<R_{t,s}\xi,\xi\>$ is positive 
   and non-increasing in $s\leq t$, the function  
   $s\mapsto \<R_{t,s}\xi,\xi\>$ on $(-\infty,t]$ is locally of bounded variation for every $\xi$ and every $t$. 
   Therefore $R(t,\xi, \cdot)$ can be extended to be a $\sigma$-finite measure on the $\sigma$-algebra $\B((-\infty,t])$ generated by $\A_t^0$.  
   
   Identities \eqref{R_t:ss} follows from \eqref{R_s,r=R_t,s-R_t,r_xi} and the continuity of $s\mapsto \<R_{t,s}\xi,\xi\>$.  \eqref{R_t:ss} means that for every $\xi\in \H$, $R(t,\xi, \cdot)$ is atomless.
     To prove \eqref{R_t:R_s}  we only need to show it for sets of the form $B=(u,v]$ with $u\leq v\leq s$.
     By \eqref{R_t:(s,r]}  we obtain
     \[
        \begin{aligned}
             R(t,\xi,(u,v])&=\<R_{v,u}U(t,v)^*\xi, U(t,v)^*\xi\> \\
                                 &=\<R_{v,u}U(s,v)^*U(t,s)^*\xi, U(s,v)^*U(t,s)^*\xi\>\\
                                 &=R(s,U(t,s)^*\xi, (u,v]). 
        \end{aligned}
    \]
    
    It remains to show the second assertion (2). For every $B\in \B^0((-\infty, t])$, we define an operator  $R(t, B)$ in \H by  setting first
        for every $\xi\in\H$, 
    \[
          \<R(t,B)\xi,\xi\>=R(t,\xi,B)
    \]
      and then  for every $\xi,\eta\in\H$
       \begin{equation}\label{Equ:Rtb:polarization}
            \<R(t,B)\xi,\eta\>=\frac14\bigl(
            \< R(t,B)(\xi+\eta), (\xi+\eta)\>
            -            \< R(t,B)(\xi-\eta), (\xi-\eta)  \>\bigr)
       \end{equation}
    using the polarization identity.    
       
       Obviously $R(t,B)$ is positive definite and self-adjoint. Now we show that its trace is bounded. Since $B$ is bounded, there exists some finite interval $(s,t]$ containing $B$. Therefore we have 
       \begin{equation*}
       \begin{aligned}
             \sum_{j=1}^\infty \<R(t,B)e_j,e_j\>
           &=\sum_{j=1}^\infty R(t,e_j,B) \leq   \sum_{j=1}^\infty R(t,e_j,(s,t])\\
           & =\sum_{j=1}^\infty \<R_{t,s}e_j,e_j\> <\infty.
       \end{aligned}
       \end{equation*}
The last inequality is due to the fact that  $R_{t,s}$ is of trace class.   
\end{proof}

For every $t\in \R$, we define a $\sigma$-finite measure $R_{\tr}(t,\cdot)$ on $((-\infty,t],\B((-\infty,t]))$ such that 
\begin{equation}\label{Rtr_def}
     R_{\tr}(t,B):=\tr R(t,B)
                =\sum_{j=1}^\infty R(t,e_j,B),  
\end{equation}
for all $B\in \B^0((-\infty, t])$. Clearly we have for every $t\geq s$, 
\[
     R_{\tr}(t,(s,t])=\tr R_{t,s}.
\]
From \eqref{R_t:R_s} it is easy to get that for every $B\in \B^0((-\infty, s])$, 
\begin{equation}\label{Rtr_t2s}
 R_{\tr}(t,B)=\tr U(t,s)R(s,B)U(t,s)^*.
\end{equation}

\begin{lemma}\label{lmm:Rt<<Rtrs}
 For every $t,s\in\R$, $t\geq s$, and $\xi\in\H$, $R(t,\xi,\cdot)$ is absolutely continuous with respect to $R_{\tr}(s,\cdot)$ on 
 $((-\infty,s],\B((-\infty,s]))$. 
\end{lemma}

\begin{proof}
 Suppose that $R(t,\xi,B)>0$ for some $B\in \B((-\infty,s])$. 
 We may assume that $B$ is bounded.
 By
 \eqref{R_t:R_s} we have
 \[\begin{aligned}
     R(t,\xi,B)&=R(s,U(t,s)^*\xi,B)
                =\<R(s,B)U(t,s)^*\xi,U(t,s)^*\xi\> \\
               &\leq |U(t,s)^*\xi|^2\tr R(s,B)
                = |U(t,s)^*\xi|^2 R_{\tr}(s,B).
   \end{aligned}
\]
So $ R_{\tr}(s,B)>0$ and thus the proof is complete. Here we have used the following simple inequality: for any non-negative definite trace class operator $A$ on \H,
\[
\<A\eta, \eta\> \leq |\eta|^2 \tr A,\quad  \eta\in\H. 
\]
Indeed,  let $\eta_i=\<\eta,e_i\>$ for $i=1,2,\cdots$.   
By Cauchy-Schwartz inequality,
we have 
\[
\begin{aligned}
            \<A\eta,\eta\>
    =  & \sum_{i,j=1}^\infty \eta_i\eta_j \<Ae_i,e_j\>
 \leq  \sum_{i,j=1}^\infty |\eta_i\eta_j| \<Ae_i,e_i\> ^{1/2}\<Ae_j,e_j\>^{1/2}\\
  \leq &
      \left( \sum_{i=1}^\infty |\eta_i|^2 \<Ae_i,e_i\> \right)^{1/2}
      \left( \sum_{j=1}^\infty |\eta_i|^2 \<Ae_j,e_j\> \right)^{1/2} \\
  \leq & |\eta|^2 \tr A. 
\end{aligned}
\]

\end{proof}

Define a Radon measure $\pi_R$ on $(\R,\B(\R))$ by setting 
\[
   \pi_R(ds)=\sum_{i=1}^\infty R_{\tr}(q_i,ds)\Eins_{(q_i-1,q_i]},
\]
where $q_i\in Q$ for all $i=1,2,\cdots$.

\begin{lemma}\label{lmm:Rt<<piR}
  For every $t\in\R$, $\xi\in\H$, $R(t,\xi,\cdot)$ is 
  absolute continuous with respect to $\pi_R$ on $((-\infty,t],\B((-\infty,t]))$.
\end{lemma}

\begin{proof}
   Suppose that for some $B\in \B((-\infty,t])$ we have 
   $R(t,\xi,B)>0$. Then    $R(t,\xi,(q_i-1,q_i]\cap B)>0$ for some $q_i\in Q\cap (-\infty,t]$. 
   By  Lemma \ref{lmm:Rt<<Rtrs}, $R(t,\xi,\cdot)$ is absolute continuous with respect to  $R_{\tr}(q_i,\cdot)$. So we also have $R_{\tr}(q_i,(q_i-1,q_i]\cap B)>0$. This implies $\pi_R(B)>0$. Hence the proof is finished. 
\end{proof}

\begin{lemma}\label{lmm:R2K_R}
  Let $t\in \R$. For all  $ s\leq t$, 
  there exists a positive definite, self-adjoint trace class operator $K_R(t,s)$ on \H such that 
  for every $\xi\in\H$ and $B\in \B^0((-\infty,t])$,  we have 
  \begin{equation}
     R(t,\xi,B)=\int_B \<K_R(t,s)\xi,\xi\>\,\pi_R(ds)
     =\left\<\left(\int_B K_R(t,s)\,\pi_R(ds) \right)\xi,\xi \right\>.
  \end{equation}
 That is, 
 \[
      R(t,B)=\int_B K_R(t,s)\,\pi_R(ds).
 \] 
  Moreover, for every $r\leq t$, we have
  \begin{equation}\label{equ:K_R(t,s)2K_R(r,s)}
    K_R(t,s)=U(t,r)^*K_R(r,s)U(t,r)
  \end{equation}
   for $\pi_R$-almost all $s\leq r$. 
\end{lemma}

\begin{proof}
   By Lemma \ref{lmm:Rt<<piR}, for every $t\in\R$, $\xi\in\H$, 
   $R(t,\xi,\cdot)$ is absolute continuous with respect to $\pi_R$. Hence for every $i,j\geq 1$, there is a locally $\pi_R$-integrable functions $K_{i,j}^+(t,\cdot)$ and $K_{i,j}^-(t,\cdot)$ on $(-\infty,t]$ such that
   for all $B\in\B^0((-\infty,t])$, we have 
   \begin{equation}\label{equ:Kmn:Def1}
      \frac14 \<R(t,B)(e_i+e_j),e_i+e_j\>
			         =\int_B K_{i,j}^+(t,s)\,\pi_R(ds)
   \end{equation}
   and 
      \begin{equation}\label{equ:Kmn:Def2}
      \frac14 \<R(t,B)(e_i-e_j),e_i-e_j\>
			         =\int_B K_{i,j}^-(t,s)\,\pi_R(ds),
   \end{equation}
   where $\{e_j\}_{j=1}^\infty$ is an orthogonal normal basis of \H. 
   
   Let 
   \[
   	K_{i,j}(t,s)=K_{i,j}^+(t,s)-K_{i,j}^-(t,s). 
   \]
   Then by \eqref{Equ:Rtb:polarization}  we have for all $B\in\B^0((-\infty,t])$,
      \begin{equation}\label{equ:Kmn:Def3}
       \<R(t,B)e_i,e_j\>
			         =\int_B K_{i,j}(t,s)\,\pi_R(ds).
   \end{equation}

   Note that $R(t,B)$  is positive definite, symmetric and has finite trace. 
By the symmetricity, we have 
$$\<R(t,B)e_i,e_j\>=\<R(t,B)e_j,e_i\>.$$ 
Therefore, 
   \begin{equation}\label{equ:Kij=Kji}
        \int_B K_{i,j}(t,s)\,\pi_R(ds)
           =\int_B K_{j,i}(t,s)\,\pi_R(ds).
   \end{equation}
By the positive definite property, we get 
   $$\<R(t,B)\xi,\xi\>\geq 0.$$
Hence 
   \begin{equation}\label{equ:Lmn:Positive}
     \int_B \sum_{i,j=1}^\infty K_{i,j}(t,s)\<\xi,e_i\>\<\xi,e_j\>\,\pi_R(ds)\geq 0
   \end{equation}
   for all finite linear combination of 
   $e_1,e_2, \cdots$, i.e. 
   $\xi\in\Span\{e_1,e_2,\cdots\}$. 
Moreover, since $B$ is bounded, we have $$\tr R(t,B)<\infty.$$
   Hence
   \begin{equation}\label{equ:finite:trace}
       \int_B \sum_{i=1}^\infty K_{i,i} (t,s)\,\pi_R(ds)<\infty.
   \end{equation}
Therefore there exists a set $F_{t,R}\in\B((-\infty,t])$ satisfying 
  $$\pi_R((-\infty,t]\setminus F_{t,R})=0$$ such that for all $s\in F_{t,R}$, 
  \[\begin{aligned}
       & K_{i,j}(t,s)=K_{j,i}(t,s),\quad i,j=1,2,\cdots, \\     
       & \sum_{i=1}^\infty K_{i,i}(t,s)<\infty, \\
       & \sum_{i,j=1}^\infty K_{i,j}(t,s)\<\xi,e_i\>\<\xi,e_j\>\geq 0
    \end{aligned}
  \]
  for all $\xi\in\Span\{e_1,e_2,\cdots\}$ with rational coefficients. 
  
  Now for any $\xi\in\H$ define 
  \[
      K_R(t,s)\xi
     = \left\{
     \begin{aligned}
         & \sum_{i,j=1}^\infty K_{i,j}(t,s)\<\xi,e_i\>e_j,& 
                 \quad & s\in F_{t,R};
        				\\
       & 0, &  \quad		& s\notin F_{t,R}.
     \end{aligned}     
     \right.
  \]
Then $K_R(t,s)$ satisfies the desired conditions by aproximation. 
  Note that 
  \eqref{equ:K_R(t,s)2K_R(r,s)} follows from \eqref{R_t:R_s}.
  This completes the proof. 
\end{proof}
  
\subsubsection{Jump part}\label{Spectral_Rep:JumpPart}

\begin{lemma}\label{lmm1:jump_part}
Assume that Assumption \ref{Ass:ContinuousMeasure'} holds. 
  For every $t\in\R$, there is a map
  \[
  \begin{aligned}
     \m(t,\cdot,\cdot)\colon 
         &\B((-\infty, t])\times \B(\H)\to [0,\infty], \\
         &A\times B\mapsto \m(t,A,B)   
  \end{aligned}
  \]
  such that
  \begin{enumerate}
   \item For every $B\in \B(\H) $, 
              $\m(t,\cdot,B)$ is a unique atomless measure on 
              $((-\infty, t]$, $\B\left((-\infty, t])\right)$ such that 
              for all $s\leq t$,
                \begin{equation}\label{equ1:lmm1:jump_part}
                   \m(t,(s,t],B)=\m_{t,s} (B).
                \end{equation}
         In particular if $B$ is contained in $\{x\in\H\colon |x|>\varepsilon\}$\, for some $\varepsilon>0$, then $\m(t,\cdot, B)$ is 
         $\sigma$-finite.                
   \item For every $A\in \B^0((-\infty, t])$, $\m(t,A,\cdot)$ is a L\'evy measure on $(\H,\B(\H))$.
   \item For every $s\leq t$, $A\in \B((-\infty, s])$, 
                 \begin{equation}\label{equ2:lmm1:jump_part}
                   \m(t,A,\cdot)=\m(s,A,\cdot)\circ U(t,s)^{-1}.
                \end{equation}    
  \end{enumerate}
\end{lemma}

\begin{proof} Fix $t\in\R$. 
 For each $B\in \B(\H)$, define for every $s\leq r\leq t$, 
 \begin{equation}\label{equ3:lmm1:jump_part}
     \m(t,(s,r],B)=\m_{r,s}\circ U(t,r)^{-1}(B).
 \end{equation}
 Then condition \eqref{equ1:lmm1:jump_part} is obviously fulfilled. 

Now we show that $\m(t,\cdot,B)$ is an additive set function on  $\A_t^0$.
By \eqref{Equ:Structure_char_triplet} we have 
 for all $s\leq \sigma\leq r\leq t$,
 \[
 \begin{aligned}
    & \m(t,(s,\sigma],B)+\m(t,(\sigma,r],B) \\
  = & \m_{\sigma,s}(U(t,\sigma)^{-1}B) +\m_{r,\sigma}(U(t,r)^{-1}B)  \\
  = & \m_{\sigma,s}(U(r,\sigma)^{-1} U(t,r)^{-1}B)
     +\m_{r,\sigma}(U(t,r)^{-1}B)  \\
  = & ( \m_{\sigma,s}\circ U(r,\sigma)^{-1} 
     +\m_{r,\sigma}) (U(t,r)^{-1}B)  \\  
  = & \m_{r,s}   (U(t,r)^{-1}B)  \\   
  = & \m(t,(s,r],B).
 \end{aligned}
 \]
Here we have used the fact that
  $$U(t,\sigma)^{-1}B=U(r,\sigma)^{-1}U(t,r)^{-1}B.$$
 Indeed, for any $x\in\H$, 
\[\begin{aligned}
                       x\in U(t,\sigma)^{-1}B 
      &\Longleftrightarrow U(t,\sigma)x\in B \\
      &\Longleftrightarrow U(t,r)U(r,\sigma)x\in B \\
      &\Longleftrightarrow U(r,\sigma)x\in U(t,r)^{-1} B \\
      &\Longleftrightarrow x\in U(r,\sigma)^{-1}U(t,r)^{-1}  B.
  \end{aligned}
\]

Therefore $\m(t,\cdot, B)$ can be extended to be a finite-additive positive set function on $\A_t^0$. 

For any $B\subset \{x\in\H\colon |x|>\varepsilon\}$, the function 
$s\mapsto \m(t,(s,t],B)=\m_{t,s}(B)$ is continuous and non-increasing. The continuity implies that $\m(t,\cdot,B)$ is atomless. Hence for any $B\subset \{x\in\H\colon |x|>\varepsilon\}$, 
 $\m(t,\cdot, B)$ on $\A_t^0$ can be extended to be a measure on $((-\infty,t],\B((-\infty,t]))$. 
 
 For all $A\in \B((-\infty,t])$, set 
 \[
 	\m(t,A,\{0\})=0. 
 \]
That is,  $\m(t,A,\cdot)$  is concentrated on $\H\setminus\{0\}$. Moreover, this is consistent with \eqref{equ3:lmm1:jump_part} since for all $s\leq r\leq t$,  
\[
   0\leq \m(t,(s,r],\{0\})\leq \m_{t,s}(\{0\})=0.
\]

Now for all $A\in \B((-\infty,t])$ and $B\in\B(\H)$ let 
\[
	\m(t,A,B):=\m(t,A,B\setminus\{0\})
		    :=\lim_{n\to\infty}\m\left(t,A,B\cap\left\{x\in\H\colon |x|>\frac1n\right\}\right).
\]

If $B$ is contained in $\{x\in\H\colon |x|>\varepsilon\}$ for some $\varepsilon>0$, then for every $s\leq t$, 
$$\m(t,(s,t],B)=\m_{t,s}(B)<\infty.$$ 
So $\m(t,\cdot,B)$ is $\sigma$-finite. 

(2) Since $A\in\B^0((-\infty, t])$ is finite, we have $A\subset (s,t]$ for some $s<t$. Then for every $B\in\B(\H)$, 
\begin{equation}\label{Equ:mtAB:bounded}
0\leq \m(t,A,B)\leq \m(t,(s,t],B)=\m_{t,s}(B).
\end{equation}
Since $\m_{t,s}$ is a L\'evy measure,  by \eqref{Equ:mtAB:bounded} we have 
\begin{equation}\label{equ4:lmm1:jump_part}
        \int_\H (1\wedge |x|^2)\,\m(t,A,dx)
   \leq \int_\H (1\wedge |x|^2)\,\m_{t,s}(dx)<\infty. 
\end{equation}
This proves that $\m(t,A,\cdot)$ is a L\'evy measure. 

(3) We only need to show \eqref{equ2:lmm1:jump_part} for set $A$ of the form $A=(u,v]$ with $u\leq v\leq s\leq t$. 
Let $B\in\B(\H)$. Then
\[\begin{aligned}
   \bigl(\m(s,(u,v],\cdot)\circ U(t,s)^{-1}\bigr) (B) 
 =&\m(s,(u,v],\cdot)( U(t,s)^{-1} B)    \\
 =&\m_{v,u}(U(s,v)^{-1}U(t,s)^{-1}B) \\
 =&\m_{v,u}(U(t,v)^{-1}B) \\
 =&\m(t,(v,u],\cdot)(B). 
  \end{aligned}
\]
\end{proof}

For every $t\in\R$, we define a measure $F(t,\cdot)$ on $((-\infty,t],\B((-\infty,t]))$ by 
\begin{equation}\label{equ1:lmm1-2:jump_part}
   F(t,A)=\int_\H (1\wedge |x|^2)\, \m(t,A,dx),\quad A\in\B((-\infty,t]).
\end{equation}
By \eqref{equ4:lmm1:jump_part} 
it is clear that $F(t,A)<\infty$ for all $A\in\B^0((-\infty,t])$. This shows that $F(t,dx)$ is $\sigma$-finite. 

Define a measure $F(ds)$ on $\R$ by 
\[
F(ds) = \sum_{i=1}^\infty F(q_i,ds)\Eins_{(q_i-1,q_i]},
\]
where $Q=\{q_i\colon i=1,2,\cdots\}$ is a countable dense subset of $\R$. 

\begin{lemma}\label{lmm2:jump_part}
For every $s\leq t$ and $B\in\B(\H)$, 
we have 
\begin{equation}\label{equ1:lmm2:jump_part}
  \m(t,\cdot,B)\ll F(t,\cdot)\ll F(s,\cdot) \ll F(\cdot).
\end{equation}
More precisely, we have 
\begin{enumerate}
\item $ \m(t,\cdot,B)$ is absolute continuous with respect to $F(t,\cdot)$ on $\B((-\infty,t])$.
\item $F(t,\cdot)$ is  is absolute continuous with respect to $F(s,\cdot)$ on $\B((-\infty,s])$.
\item $F(s,\cdot)$ is  is absolute continuous with respect to $F(\cdot)$ on $\B((-\infty,s])$.
\end{enumerate}
\end{lemma}

\begin{proof}
 (1) Suppose that for some $A\in\B^0((-\infty,t])$ we have $\m(t,A,B)>0$. Then we also have $\m(t,A,B\setminus \{0\})>0$ since 
      $\m(t,A,\{0\})=0$. Hence 
      \[
          F(t,A)=\int_\H (1\wedge |x|^2)\,\m(t,A,dx)\geq 
	              \int_{B\setminus \{0\}} (1\wedge |x|^2)\,\m(t,A,dx)>0.
      \]
     This proves that  $ \m(t,\cdot,B)$ is absolute continuous with respect to $F(t,\cdot)$ on $\B((-\infty,t])$.

   (2) Suppose that for some  $A\in\B^0((-\infty,s])$ we have $F(t,A)>0$. Then 
        \[
           \begin{aligned}
	       0<F(t,A)&=\int_\H (1\wedge |x|^2)\,\m(t,A,dx) \\
	                     &=\int_\H (1\wedge |x|^2)\,\m(s,A,\cdot)\circ U(t,s)^{-1}(dx) \\ 
	                     &=\int_\H (1\wedge |U(t,s)x|^2)\,\m(s,A,dx) \\
	                     &\leq (1\vee |U(t,s)|^2) \int_\H (1\wedge |x|^2)\,\m(s,A,dx) \\
	                     &=(1\vee |U(t,s)|^2) F(s,A).
	  \end{aligned}
        \]
     This implies $F(s,A)>0$. Therefore $F(t,\cdot)\ll F(s,\cdot)$ on $\left((-\infty, s], \B((-\infty,s])\right)$. 
     
     (3)   Suppose that  $F(s,A)>0$ for some $A\in\B((-\infty,s])$. Then 
            $$F(s, A\cap (q_i-1,q_i])>0$$for some $q_i\in Q$ with $q_i\leq s$. Since we have proved    
                    $F(s,\cdot)\leq F(q_i,\cdot)$ previously, we get $$F(q_i, A\cap (q_i-1,q_i])>0.$$ This proves that $F(A)>0$. 
            Hence $F(s,\cdot)\ll F(\cdot)$ on $((-\infty,s],\B((-\infty,s]))$. 
\end{proof}

The following lemma is inspired by \cite[Lemma 2.2]{RR89}. 

\begin{lemma}\label{lmm3:jump_part}
For every $t\in\R$, there exists a unique $\sigma$-finite measure $M(t,\cdot)$ on the product space $\left((-\infty,t]\times\H,\B((-\infty,t])\times \B(\H)\right)$ such that
\begin{equation}\label{equ1:lmm3:jump_part}
  M(t,A\times B)=\m(t,A,B),\quad A\in\B((-\infty,t]), B\in \B(\H). 
\end{equation}
Moreover, there exists a function $J(t,\cdot,\cdot)\colon (-\infty, t]\times \B(\H)\to [0,\infty]$ such that 
\begin{enumerate}
 \item  For every $s\in (-\infty,t]$, $J(t,s,\cdot)$ is a L\'evy measure on $(\H,\B(\H))$.
 \item For every $B\in\B(\H)$, $J(t,\cdot,B)$ is a Borel measurable function on $(-\infty,t]$.
 \item For every   $\B((-\infty,t])\times \B(\H)$-measurable function $h\colon (-\infty,t]\times \H\to [0,\infty]$, 
           \begin{equation}\label{lmm3:jump_part:int_A_H_J(t,s,dx):general}
                 \int_{(-\infty,t]\times\H} h(s,x) M(t,ds\times dx) = \int_{(-\infty,t]}\left[\int_\H h(s,x) J(t,s,dx) \right] F(ds).
	  \end{equation}
	 In particular, for all $\B(\H)$-measurable function $f\colon \H\to [0,\infty]$, we have
	  \begin{equation}\label{lmm3:jump_part:int_A_H_J(t,s,dx)}
                 \int_{\H} f(x)\, \m(t,A,dx) = \int_A\int_\H f(x) J(t,s,dx)\, F(ds).
	  \end{equation} 
\item For every $t\geq r$, 
\begin{equation}\label{lmm3:jump_part:J_tJ_r}
   J(t,s,dx)=J(r,s,\cdot)\circ U(t,r)^{-1}(dx)
\end{equation}
for $F$-almost all $s\in(-\infty,r]$.	  
\end{enumerate}
\end{lemma}

\begin{proof}
Define for every $A\in \B((-\infty,t])$ and $B\in\B(\H)$, 
\begin{equation}
F_0(t,A,B)=\int_B (1\wedge |x|^2)\,\m(t,A,dx).
\end{equation}
For every 
$A\in \B^0((-\infty,t])$, $F_0(t,A,\cdot)$ is a finite measure on $(\H,\B(\H))$. 
And for every $B\in\B(\H)$, $F_0(t,\cdot,B)$ is a $\sigma$-finite measure on $((-\infty,t],\B((-\infty,t]))$.
Note that we also have $F(t,A)=F_0(t,A,\H)$.

Therefore, 
there is a measure 
$G(t,\cdot)$ on the product space $((-\infty,t]\times\H$,$\B((-\infty,t])\times \B(\H))$ such that 
\[
   G(t,A\times B) = F_0(t,A,B) = \int_A q(t,s,B) F(t,ds),
\]
where $q(t,s,dx)$ satisfies the following conditions (see \cite[Proposition 2.4]{RR89} for a brief proof):
\begin{enumerate}[(a)]
\item For every $s\leq t$, $q(t,s,\cdot)$ is a probability measure on $(\H,\B(\H))$. 
\item For every $B\in\B(\H)$, $q(t,\cdot,B)$ is $\B((-\infty,t])$-measurable. Moreover, $q(t,s,\cdot)$ is unique in the sense of $F(t,\cdot)$-almost everywhere. 
\end{enumerate}

By Lemma \ref{lmm2:jump_part} we have $F(t,\cdot)\ll F(\cdot)$.  
We choose one version of the Radon-Nikodym derivative 
$$f_t(s):=\frac{dF(t,\cdot)}{dF(\cdot)} (s),\quad s\leq t $$
such that $f_t$ is finite everywhere. 

For every $s\leq t$, we define a measure $J(t,s,dx)$ concentrating on $\H\setminus\{0\}$ by setting
\[
	J(t,s,\{0\})=0
\]
and for all $B\in\B(\H\setminus\{0\})$
\[
    J(t,s,B)=f_t(s) \cdot \int_B \frac{1}{1\wedge |x|^2}\cdot q(t,s,dx).
\]
Thus we have 
\[
     	\int_\H (1\wedge |x|^2) J(t,s,dx)
      = f_t(s)  \cdot \int_\H q(t,s,dx) 
      = f_t(s) <\infty. 
\]
It shows that for every $s\leq t$, $J(t,s,\cdot)$ is a L\'evy measure on $(\H,\B(\H))$. It is clear that for every $B\in\B(\H)$, $J(t,\cdot,B)$ is measurable. 

Let $M(t,\cdot)$ be the measure on $\B((-\infty,t])\times \B(\H)$ such that for all 
$A\in\B((-\infty,t])$ and $B\in\B(\H)$, 
\begin{equation}\label{M_tAB}
	M(t,A\times B)=\int_A F(ds) \int_B J(t,s,dx).
\end{equation}
By definition, it is standard to get \eqref{lmm3:jump_part:int_A_H_J(t,s,dx):general}. 

Now we show \eqref{equ1:lmm3:jump_part}. For every 
$ A\in\B((-\infty,t])$, $B\in \B(\H)$ we get 
\[
\begin{aligned}
M(t,A\times B)
&=\int_A f_t(s)\,F(ds)\int_B  \frac{1}{1\wedge |x|^2}\cdot q(t,s,dx) \\
  &=\int_B \frac{1}{1\wedge |x|^2} \int_A  q(t,s,dx) F(t,ds) \\
 &=\int_{A\times B} \frac{1}{1\wedge |x|^2} G(t,ds,dx)\\
                     &=\int_{B} \frac{1}{1\wedge |x|^2}F_0(t,A,dx)\\
                     &=\m(t,A,B). 
\end{aligned}
\]
So \eqref{equ1:lmm3:jump_part} holds.

By taking 
$h(s,x)=f(x) \Eins_{A}(s)$ in  \eqref{lmm3:jump_part:int_A_H_J(t,s,dx):general} and using \eqref{equ1:lmm3:jump_part}, we get 
identity \eqref{lmm3:jump_part:int_A_H_J(t,s,dx)} imediately. 

It remains to show  \eqref{lmm3:jump_part:J_tJ_r} holds almost surely.
By \eqref{equ2:lmm1:jump_part}  and \eqref{lmm3:jump_part:int_A_H_J(t,s,dx)} we get 
 for all $t\geq r$ and $A\in\B((-\infty,r])$, 
 \[
 \begin{aligned}
 	      \int_A\int_\H f(x)\, J(t,s,dx)\, F(ds) 
	 =&\int_{\H} f(x)\, \m(t,A,dx) \\
	 =&\int_{\H} f(x)\, \m(r,A,\cdot)\circ U(t,r)^{-1}(dx) \\
	 =&\int_{\H} f(U(t,r)x)\, \m(r,A,\cdot) \\
	 =&\int_A\int_{\H} f(U(t,r)x)\, J(r,s,dx)\, F(ds) \\
	 =&\int_A\int_{\H} f(x)\, J(r,s,\cdot)\circ U(t,r)^{-1}(dx)\, F(ds)	 
 \end{aligned}
 \]
This shows that \eqref{lmm3:jump_part:J_tJ_r} holds for $F$-almost all $s\in(-\infty,r]$.
\end{proof}

\subsubsection{Drift part}\label{Spectral_Rep:DriftPart}
Assume that Assumption \ref{Ass:ContinuousMeasure'} holds. 
For every $t\geq s$ let 
$\mu_{t,s}=[a_{t,s}, R_{t,s},\m_{t,s}]$. 
\begin{definition}\label{Def:natural_mu_ts}
We say that $(\mu_{t,s})_{t\geq s}$ is \emph{natural} if  for each $\xi\in\H$, $t\in\R$, 
the function $s\mapsto \<a_{t,s},\xi\>$ on $(-\infty,t]$ is locally of bounded variation, i.e. $s\mapsto \<a_{t,s},\xi\>$ is of bounded variation on any finite interval of $(-\infty,t]$.
\end{definition}

\begin{lemma}\label{lmm1:drift_part}
Suppose that $(\mu_{t,s})_{t\geq s}$ is natural. 
Then for every $t\in\R$, there is a map
\[\begin{aligned}
a(t,\cdot,\cdot)\colon & \H\times \B((-\infty,t])\to [0,\infty]\\
& (\xi,A)\mapsto a(t,\xi,A)
\end{aligned}
\]
such that 
\begin{enumerate}
\item For every $\xi\in\H$, $a(t,\xi,\cdot)$ is a singed measure on $((-\infty,t],\B((-\infty,t]))$ such that 
		\begin{equation}\label{equ1:lmm1:drift_part}
			a(t,\xi,(s,t])=\<a_{t,s},\xi\>.
		\end{equation}
\item For every $A\in\B((-\infty,t])$, $a(t,\cdot,A)$ is a linear functional on $\H$. That is, for all $c_1,c_2\in\R$ and $\xi_1,\xi_2\in\H$, 
		\begin{equation}\label{equ2:lmm1:drift_part}
			a(t,c_1\xi_1+c_2\xi_2,A)=c_1a(t,\xi_1,A)+c_2a(t,\xi_2,A).
		\end{equation}
Therefore, there is an element $a(t,A)\in\H$ such that 
\[
        a(t,\xi,A)=\<a(t,A),\xi\>
\]		
and hence for all $s\leq t$ 
\begin{equation}\label{a_st:a_ts}
	a(t,(s,t])=a_{t,s}.	
\end{equation}
\end{enumerate}
\end{lemma}

\begin{proof}
For every interval $(s,r]\subset (-\infty,t]$, $\xi\in\H$, we define
		\begin{equation}\label{equ3:lmm1:drift_part}
				a(t,\xi,(s,r])=\<a_{t,s},\xi\>-\<a_{t,r},\xi\>.
		\end{equation}
Identity \eqref{equ1:lmm1:drift_part}		
follows from \eqref{equ3:lmm1:drift_part} by noting that 
$a_{t,t}=\delta_0$	(cf. \eqref{Equ:a_ttR_ttM_tt=0}) and the continuity of $r\mapsto a_{t,r}$ with $r\leq t$.

Obviously \eqref{equ3:lmm1:drift_part} defines an additive set function on the algebra $\A_t^0$. 
Since the function  $s\mapsto \<a_{t,s},\xi\>$ is continuous and of locally bounded variation, the set function 
$a(t,\xi,\cdot)$ can be extended to be a signed measure on $\B((-\infty,t])$. 

From \eqref{equ3:lmm1:drift_part} we get that $a(t,\cdot,A)$ is a linear functional for sets of the form $A=(s,r]$ with $s\leq r\leq t$. By Dynkin's  $\pi$-$\lambda$ theorem, we obtain the linearity  claimed by \eqref{equ2:lmm1:drift_part}.

The existence of $a(t,A)$ is due to Riesz representation theorem. The identity \eqref{a_st:a_ts}
follows from \eqref{equ1:lmm1:drift_part}.
\end{proof}

The following lemma shows that the connection between $a(t,\cdot,\cdot)$ and $a(s,\cdot,\cdot)$ is not 
as simple as the case for the Gaussian and jump parts (cf. \eqref{R_t:R_s} and \eqref{equ2:lmm1:jump_part}). 

\begin{lemma}\label{lmm2:drift_part}
 For every $\xi\in\H$, $u\leq v\leq s\leq t$, we have
 \begin{equation}\label{equ1:lmm2:drift_part}
 \begin{aligned}
     a(t,\xi,(u,v])=& a(s,U(t,s)^*\xi, (u,v]) \\
     & +
     \int_{\H\setminus\{0\}} U(t,v)x\left[
     					\frac{1}{1+|U(t,v)x|^2} -  \frac{1}{1+|U(s,v)x|^2}
     \right]\,\m_{v,u}(dx).
\end{aligned}     
 \end{equation}
\end{lemma}

\begin{proof}
By \eqref{equ3:lmm1:drift_part} and  \eqref{Equ:Structure_char_triplet}, we have 		
\begin{equation}\label{equ2:lmm2:drift_part}
\begin{aligned}
&a(t,\xi,(u,v])
=\<a_{t,u},\xi\>-\<a_{t,v},\xi\> \\
=& \left\<\xi,U(t,v)a_{v,u} +
           \int_{\H\setminus \{0\}} U(t,v)x\left[ 
		\frac{1}{1+|U(t,v)x|^2}-\frac{1}{1+|x|^2}
	 \right]\m_{v,u}(dx)\right\>.
\end{aligned}	 
\end{equation}	 
Similarly we have
\begin{equation}\label{equ3:lmm2:drift_part}
\begin{aligned}
&a(s,U(t,s)^*\xi,(u,v])\\
=& \left\<U(t,s)^*\xi,U(s,v)a_{v,u} +
           \int_{\H\setminus \{0\}} U(s,v)x\left[ 
		\frac{1}{1+|U(s,v)x|^2}-\frac{1}{1+|x|^2}
	 \right]\m_{v,u}(dx)\right\>. \\
=& \left\<\xi,U(t,v)a_{v,u} +
           \int_{\H\setminus \{0\}} U(t,v)x\left[ 
		\frac{1}{1+|U(s,v)x|^2}-\frac{1}{1+|x|^2}
	 \right]\m_{v,u}(dx)\right\>. 	 
\end{aligned}	 
\end{equation}	
Comparing \eqref{equ2:lmm2:drift_part} and \eqref{equ3:lmm2:drift_part} we arrive at
\eqref{equ1:lmm2:drift_part}.
\end{proof}

Due to the complicated relationship \eqref{equ1:lmm2:drift_part} 
between $a(t,\cdot,\cdot)$ and $a(s,\cdot,\cdot)$, 
we need the following hypothesis.
\begin{assumption}\label{assumption_drift_gamma}
Assume that there is a $\sigma$-finite atomless measure $\gamma$ on $(\R,\B(\R))$ such that 
for every $t\in\R$, there is a measurable function $b(t,\cdot)$ on $(-\infty,t]$ 
taking values in \H such that 
\[
     a(t,A)=\int_A b(t,s)\,\gamma(ds),\quad A\in\B((-\infty,t]).
\]
\end{assumption}

\begin{remark}
We note here why method using similar ideas for the Gaussian and jump part does not work. 
For every $t\in\R$ and $\xi\in\H$, let 
$|a|_{\var}(t,\xi,\cdot)$ denote the total variation measure of $a(t,\xi,\cdot)$ on $((-\infty,t],\B(-\infty,t]))$.
Define for every $t\in\R$,  for every $A\in\B((-\infty,t])$, 
\begin{equation}
  \gamma(t,A):=\sum_{i=1}^\infty |a|_{\var}(t,e_i,A).
\end{equation}

Assume further that for every $A\in \B^0((-\infty,t])$, $\gamma(t,A)<\infty$. 
We may try to define for every $A\in\B((-\infty,t])$, 
\[
{\gamma}(A):=\sum_{i=1}^\infty \gamma(q_i,A)\Eins_{(q_i-1,q_i]}.
\]
But because the absolute continuity of  $\gamma(t,\cdot)$ with respect to $\gamma(s,\cdot)$ for $s<t$ is not known, we are not able to show that  $a(t,\xi,\cdot)\ll \gamma(\cdot)$ for every $t\in\R$ and $\xi\in\H$. 
\end{remark}

\subsubsection{Combining Gaussian, jump and drift parts}\label{subsubsec:comb}
For every $t\in\R$, $A\in\B^0((-\infty,t])$, define 
\begin{equation*}
\begin{aligned}
   \psi(t,A)(\xi)= & -i\<a(t,A),\xi\>
     +\frac12\<R(t,A)\xi,\xi\> \\
     &\quad\quad 
      -\int_{\H\setminus\{0\}} \left( \e^{i\<\xi,x\>}
     -1-\frac{1}{1+|x|^2}\right)\,\m(t,A,dx),\quad \xi\in\H.
\end{aligned}
\end{equation*}
By Lemma \ref{lmm:gauss:scattered:measure}, Lemma \ref{lmm1:jump_part} and Assumption \ref{assumption_drift_gamma}, there is an
infinitely divisible measure 
$\mu(t,A)$ on \H associated with $\psi(t,A)$: 
\[
    \widehat{\mu(t,A)}(\xi):=\exp(-\psi(t,A)(\xi)),\quad \xi\in\H.
\]

In terms of the characteristic exponent of $\mu(t,A)$ we have
\[
    \mu(t,A)=[a(t,A),R(t,A),\m(t,A,\cdot)]. 
\]
In particular, according to \eqref{R_t:(s,t]}, \eqref{equ1:lmm1:jump_part} and 
\eqref{a_st:a_ts}
we have respectively 
\[			R(t,\xi, (s,t])=\<R_{t,s}\xi,\xi\>,\quad
	                   \m(t,(s,t],\cdot)=\m_{t,s} (\cdot),\quad 
			a(t,(s,t])=a_{t,s}		                   
\]
for all $s\leq  t$. Hence
\begin{equation}\label{Equ:mu_srt}
    \mu(t,(s,t])=\mu_{t,s},\quad s\leq t.
\end{equation}

Note that for all $A,B\in\B^0((-\infty,t])$ such that $A\cap B=\emptyset$, we have 
\begin{equation}\label{additive_psi(t,A,xi)}
   \psi(t,A\cup B)(\xi)=\psi(t,A)(\xi)+\psi(t,B)(\xi). 
\end{equation}
This implies 
\begin{equation}\label{Equ:mu_A+B}
    \mu(t,A)*\mu(t,B)=\mu(t,A\cup B).
\end{equation}
Obviously the finite additivity in \eqref{additive_psi(t,A,xi)} can be extended to be $\sigma$-additivity. 

Recall that for all $t\geq s$, $\xi\in\H$,  we have set $\psi_{t,s}(\xi)=-\log \hat{\mu}_{t,s}(\xi)$. Therefore 
\begin{equation}\label{psi_st_psi_ts_xi}
    \psi(t,(s,t])(\xi)=\psi_{t,s}(\xi).
\end{equation}

\begin{lemma}\label{lmm:psi_t_A_s_A:Uts_xi} 
	For every $t\geq s$ and all $A\in\B^0((-\infty,s])$, 
	\begin{equation}\label{equ:psi_t_A:psi_s_A}
	 \psi(t,A)(\xi) = \psi(s,A)(U(t,s)^*\xi),\quad \xi\in\H.
	\end{equation}
Therefore
	\begin{equation}\label{equ:mu_t_A:mu_s_A}
	 \mu(t,A) = \mu(s,A)\circ U(t,s)^{-1}.
	\end{equation}
In particular, for all $s\leq r\leq t$, we have
	\begin{equation}\label{equ:mu_t_sr=}
	 \mu(t,(s,r]) = \mu_{r,s}\circ U(t,r)^{-1}.
	\end{equation}
\end{lemma}

\begin{remark}
This lemma shows that Equation 
 \eqref{Equ:mu_A+B} generalizes 
 \eqref{Equ:SkewMeasure}.
 Indeed, let $s\leq r\leq t$. Substituting $A=(r,t]$, $B=(s,r]$ into  \eqref{Equ:mu_A+B} and using  equations \eqref{Equ:mu_srt} and 
 \eqref{equ:mu_t_sr=}, we get
 \eqref{Equ:SkewMeasure}.
\end{remark}

To show Lemma \ref{lmm:psi_t_A_s_A:Uts_xi}, we first note the following simple fact. 

\begin{lemma}\label{lmm:mu_Uts_xi} 
	For every $t\geq s$ and any probability measure $\mu$ on $(\H,\B(\H))$, we have 
	\begin{equation}\label{Equ:mu_Uts_xi}
		\hat{\mu}(U(t,s)^*\xi) = \hat{\mu}(U(r,s)^*U(t,r)^*\xi),\quad \xi\in\H.
	\end{equation}
\end{lemma}
\begin{proof}
	For all $\xi\in\H$, we have 
	\[\begin{aligned}
		\hat{\mu}(U(t,s)^*\xi) &= \int_\H \e^{i\<U(t,s)^*\xi,x\>}\,d\mu(x) \\
						    &= \int_\H \e^{\<\xi,U(t,r)U(r,s)x\>}\,d\mu(x) \\
&= \int_\H \e^{i\<U(r,s)^*U(t,r)^*\xi,x\>}\,d\mu(x)
						    = \hat{\mu}(U(r,s)^*U(t,r)^*\xi).
	\end{aligned}					    
	\]
\end{proof}

\begin{proof}[Proof of Lemma \ref{lmm:psi_t_A_s_A:Uts_xi}]
We only need to show \eqref{equ:psi_t_A:psi_s_A} for set $A$ of the form $A=(v,u]$ with 
$v\leq u\leq s\leq t$. By \eqref{additive_psi(t,A,xi)}, \eqref{psi_st_psi_ts_xi} and \eqref{Equ:SkewMeasure:Symbol} we have
\begin{equation}\label{Equ:psi_t_vu_xi_psi_uv_Utu_xi}
\begin{aligned}
       \psi(t,(v,u])(\xi)
=& \psi(t,(v,t])(\xi)  - \psi(t,(u,t])(\xi) \\
=& \psi_{t,v}(\xi) - \psi_{t,u}(\xi) 
= \psi_{u,v}(U(t,u)^*\xi). 
\end{aligned}
\end{equation}
Analogously we have 
\begin{equation}\label{Equ:psi_s_vu_xi_psi_uv_Usu_xi}
       \psi(s,(v,u])(\xi)
       = \psi_{u,v}(U(s,u)^*\xi). 
\end{equation}
By Lemma \ref{lmm:mu_Uts_xi} we obtain 
\[
	\psi_{u,v}(U(t,u)^*\xi)=\psi_{u,v}(U(s,u)^*U(t,s)^*\xi). 
\]
Therefore, by comparing \eqref{Equ:psi_t_vu_xi_psi_uv_Utu_xi} and \eqref{Equ:psi_s_vu_xi_psi_uv_Usu_xi} we get 
\[
     \psi(t,(v,u])(\xi)=     \psi(s,(v,u])(U(t,s)^*\xi).
\]
This completes the proof. 
\end{proof}

We define a \emph{control measure} $\pi$ on $(\R,\B(\R))$ by 
\[
     \pi(A)=\gamma(A)+\pi_R(A)+F(A),\quad A\in\B(\R).
\]
It is clear that $\pi$ is atomless. Moreover,  $\gamma$, $\pi_R$ and $F$ all are absolute continuous with respect to $\pi$. Let
\[
h^{(1)}(s)=\frac{d\gamma}{d\pi}(s),\quad
h^{(2)}(s)=\frac{d\pi_R}{d\pi}(s),\quad
h^{(3)}(s)=\frac{d F}{d\pi}(s),\quad s\in\R.
\]

For all $s\leq t$, let 
\[
\begin{aligned}
     \rho(t,s)(\xi):=& - i\< h^{(1)}(s)b(t,s),\xi \> +  \frac12 \left< h^{(2)}(s) K_R(t,s)\xi,\xi \right> \\
                            &  {} - \int_{\H\setminus\{0\}} \left(
                                    \e^{i\<\xi,x\>} - 1 -\frac{i\<\xi,x\>}{1+|x|^2} 
                                \right) h^{(3)}(s) J(t,s,dx),\quad \xi\in\H.
\end{aligned}
\]
It is clear that 
$\rho(t,s)$ is the characteristic functional of some infinitely divisible measure on \H.

By
Lemma \ref{lmm:R2K_R},
Lemma \ref{lmm3:jump_part}
and
Assumption \ref{assumption_drift_gamma} 
we have the following spectral representation theorem. 


\begin{theorem}
Suppose that Assumption \ref{Ass:ContinuousMeasure'},
Assumption \ref{assumption_drift_gamma} hold and  
that $(\mu_{t,s})_{t\geq s}$ is natural. 
Then  for every $t\in\R$ and  $A\in\B^0((-\infty,t])$ 
\begin{equation}\label{Equ:widehat-mu(t,A)-summary}
\widehat{\mu(t,A)}(\xi) = \exp\left(-\psi(t,A)(\xi)\right)
= \exp\left(  -\int_A \rho(t,s)(\xi)\,d\pi(s) \right),\quad \xi\in\H.
\end{equation}
In other words, 
\[
\begin{aligned}
\mu(t,A)=
&\left[ \int_A b(t,s)h^{(1)}(s)\,d\pi(s),
\int_A K_R(t,s)h^{(2)}(s)\,d\pi(s), \right. \\
&\left. \quad\quad\quad\quad\quad\quad\quad\quad\quad\quad\quad
\int_A 
		     J(t,s,dx) h^{(3)}(s)\,d\pi(s)
\right].
\end{aligned}
\]
In particular, for all $s\leq t$, we have 
\[
\hat{\mu}_{t,s}(\xi) = \exp\left(  -\int_s^t \rho(t,r)(\xi)\,d\pi(r) \right),\quad \xi\in\H.
\]
Moreover, for any fixed $t\geq r$, $\xi\in\H$, we have
\begin{equation}\label{Equ:rho_ts_rho_rs_Uts}
\rho(t,s)(\xi)=\rho(r,s)(U(t,r)^*\xi)
\end{equation}
for $\pi$-almost all $s\leq r$.
\end{theorem}

\begin{proof}We only need show \eqref{Equ:rho_ts_rho_rs_Uts}. 
By \eqref{Equ:widehat-mu(t,A)-summary} and \eqref{equ:psi_t_A:psi_s_A} we have 
for all $\xi\in\H$,  $t\geq r$ and $A\in(-\infty,r]$, 
\[
\begin{aligned}
	   \int_A \rho(t,s)(\xi)\,d\pi(s)
        =&\psi(t,A)(\xi) 
        =\psi(r,A)(U(t,r)^*\xi)\\
        =&\int_A \rho(r,s)(U(t,r)^*\xi) \,d\pi(s).
\end{aligned}
\]
Hence for all $A\in\B^0((-\infty,r])$, it holds
\[
\int_A \rho(t,s)(\xi)\,d\pi(s)= \int_A \rho(r,s)(U(t,r)^*\xi) \,d\pi(s).
\]
This proves that 
\[\rho(t,s)(\xi)=\rho(r,s)(U(t,r)^*\xi)\]
holds for $\pi$-almost all $s\leq r$.
\end{proof}

\begin{remark}\label{Remark:convolution-integral-representation}
If we have (comparing with \eqref{Equ:rho_ts_rho_rs_Uts}) 
\[
	\rho(t,s)=\rho(s,s)(U(t,s)^*\xi)=:\rho_s(U(t,s)^*\xi)
\]
for all $\xi\in\H$ and $t\geq s$, then 
\[
\hat{\mu}_{t,s}(\xi) = \exp\left(  -\int_s^t \rho_r(U(t,r)^*\xi)\,d\pi(r) \right),\quad \xi\in\H.
\]
It is clear that there is an additive process  $(Z_t)_{t\in\R}$   with factoring $((\rho_r)_{r\in\R},\pi)$. Then by Example \ref{Exa:represation:factoring}, for all $t\geq s$, $\mu_{t,s}$ is the distribution of $\int_s^t U(t,r)\,dZ_r$.

That is, under some conditions, $(\mu_{t,s})_{t\geq s}$ satisfying \eqref{Equ:SkewMeasure} corresponds one to one the distribution of stochastic convolution integrals of $U(t,s)$ with respect to some additive process. 
\end{remark}

\subsubsection{Space homogeneous case}\label{subsubsec:time_homo_factoring}
Let us consider the special case of the skew convolution equation \eqref{Equ:SkewMeasure} when $U\equiv I$, that is, 
\begin{equation}\label{Equ:mu_ts=mu_tr_mu_rs}
      \mu_{t,s}=\mu_{t,r}*\mu_{r,s},\quad t\geq r\geq s.
\end{equation}
We also use Assumption \ref{Ass:ContinuousMeasure'}. 
We are going to show a factoring of $\mu_{t,s}$. It is an infinite dimensional generalization of the factoring of finite dimensional additive process studied by Sato \cite{Sat06}.

Let $\mu_{t,s}=[a_{t,s},R_{t,s},\m_{t,s}]$.  From \eqref{Equ:mu_ts=mu_tr_mu_rs} we have the following equations which
are \eqref{Equ:Structure_char_triplet} in the case $U\equiv I$:
\begin{equation}\label{Equ:Structure_char_triplet_U=I}
\begin{aligned}
a_{t,s}&=a_{t,r}+a_{r,s},\\
R_{t,s}&=R_{t,r}+R_{r,s},\\
\m_{t,s}&=\m_{t,r}+\m_{r,s}
\end{aligned}
\end{equation}
for every $t\geq r\geq s$. 

As proved in Lemma \ref{lmm:gauss:scattered:measure} there is a function 
\[
R(\cdot,\cdot)\colon \H\times \B(\R) \ni (\xi,B) \mapsto R(\xi,B)\in [0,\infty]
\]
such that 
\begin{enumerate}
\item For all $\xi\in \H$, $R(\xi,\cdot)$ is a $\sigma$-finite measure on $(\R,\B(\R))$ such that 
for all $s<t$, 
$$R(\xi,(s,t])=\<R_{t,s}\xi,\xi\>.$$
\item For all $B\in \B^0(\R)$, 
\[
  \< R(B)\xi,\xi\> = R(\xi,B) ,\quad \xi\in\H
\]
defines a non-negative definite, self-adjoint trace class operator $R(B)$ on \H. 
\end{enumerate}

As in Section \ref{Spectral_Rep:GaussianPart} we set for every $B\in\B^0(\R)$, 
\[
\tilde{\pi}_R(B):=\tr R(B) =  \sum_{j=1}^\infty \<R(B)e_j,e_j\> = \sum_{j=1}^\infty R(e_j,B).
\]
According to Lemma \ref{lmm:R2K_R}, for every $s\in\R$ there is a non-negative definite, self-adjoint trace class operator $K(s)$ on \H such that 
\[
     R(\xi,B)=\<R(B)\xi,\xi\> = \int_B \<K(s)\xi,\xi\>\,\tilde{\pi}_R(ds),\quad B\in \B^0(\R).
\]
In other words, 
\[
R(B)=\int_B K(s) \tilde{\pi}_R (ds),\quad B\in \B^0(\R). 
\]

Similar to the proof of Lemma \ref{lmm1:jump_part}, there is a map 	
\[
  \begin{aligned}
     {\m}(\cdot,\cdot)\colon 
         &\B(\R)\times \B(\H)\to [0,\infty], \\
         &A\times B\mapsto {\m}(A,B)   
  \end{aligned}
  \]
  such that
  \begin{enumerate}
   \item For every $ B\in \B(\H) $, 
              ${\m}(\cdot,B)$ is a unique measure on 
              $(\R,\B(\R))$ satisfying 
               \begin{equation}\label{equ1:lmm1:jump_part'}
                   {\m}((s,t],\cdot)=\m_{t,s}(\cdot) 
                \end{equation}
         on $\B(\H)$ for all $t\geq s$. 
         In particular if $B\subset \{x\colon |x|>\varepsilon\}$ for some $\varepsilon>0$, then ${\m}(\cdot, B)$ is 
         $\sigma$-finite.                
   \item For every $A\in \B^0(\R)$, ${\m}(A,\cdot)$ is a L\'evy measure on $(\H,\B(\H))$.   
  \end{enumerate}

Let 
\[
\tilde{F}(ds) = \int_\H (1\wedge |x|^2)\, \m(ds,dx).
\]
Then for all $B\in \B(\H)$, $\m (ds,B)\ll \tilde{F}(ds)$. 
Moreover, there is a function 
\[
    {J}(\cdot,\cdot)\colon \R\times \B(\H)\to [0,\infty]
\]
such that 
\begin{enumerate}
\item For every $s\in\R$, $J(s,\cdot)$ is a L\'evy measure on $\B(\H)$.
\item For every $B\in\B(\H)$, $J(\cdot,B)$ is a Borel measurable function 
on $\R$ satisfying 
\[
    \int_{\R\times \H} h(s,x)\,{\m}(ds,dx)
  =\int_{\R}\left(\int_{\H} h(s,x)\,{J}(s,dx)\right) \tilde{F}(ds)
\]
for all $\B(\R)\times \B(\H)$-measurable function $h\colon \R\times \H\to [0,\infty]$. 
\end{enumerate}

Suppose that $(\mu_{t,s})_{t\geq s}$ is natural. That is, 
for every $\xi\in\H$ and every $t\in\R$, $s\mapsto \<a_{t,s},\xi\>$ is of locally finite variation. Then there is a function 
\[
a(\cdot,\cdot)\colon \H\times \B(\R)\ni (\xi,A)\mapsto
a(\xi,A)\in [0,\infty]
\]
such that
\begin{enumerate}
\item For every $\xi\in\H$, $a(\xi,\cdot)$ is a measure on $(\R,\B(\R))$.
\item For every $A\in\B(\R)$, $a(\cdot,A)$ is a linear functional on $(\H,\B(\H))$. 
\end{enumerate}
Hence there exists an element $a(A)\in\H$ such that 
\[
\<a(A),\xi\> = a(\xi,A),\quad \xi\in\H. 
\]

For every $\xi\in\H$, let $|a|_{\var}(\xi,\cdot)$ be the total variational measure of $a(\xi,\cdot)$.  We shall use the following assumption. 

\begin{assumption}\label{assumption:summable:total_var}
Assume that for every $A\in\B^0(\R)$, 
\[
    \tilde{\gamma}(A):=\sum_{i=1}^\infty |a|_{\var}(e_i,A)<\infty.
\]
\end{assumption}

We have the following proposition. 

\begin{prop} 
Assume that $(\mu_{t,s})_{t\geq s}$ is natural and Assumption
\ref{Ass:ContinuousMeasure'},  
Assumption \ref{assumption:summable:total_var} hold. 
Then there is a measurable map $\R\ni s\mapsto b(s)\in\H$ such that for all $A\in \B(\R)$, 
\[
a(A)=\int_A {b}(s)\,\tilde{\gamma}(ds).
\]
\end{prop}

\begin{proof}
It is clear that for all $i=1,2,\cdots$, 
\(
     {a}(e_i,\cdot)\ll \tilde{\gamma}(\cdot).
\) Hence  there exists a measurable function ${b}_i$ on $\R$ such that 
\[
  {a}(e_i,A)=\int_A {b}_i(s)\,\tilde{\gamma}(ds),\quad A\in\B(\R)
\]
for all $i=1,2,\cdots$.

By assumption \ref{assumption:summable:total_var}, we have 
\[
  \int_A \sum_{i=1}^\infty |{b}_i(s)|\,\tilde{\gamma}(ds)
=\sum_{i=1}^\infty\int_A  |{b}_i(s)|\,\tilde{\gamma}(ds)
\leq \sum_{i=1}^\infty |a|_{\var}(e_i,A)<\infty. 
\]
Hence there exists a measurable set $\tilde{N}_0$ in $\R$ such that 
$\tilde{\gamma}(\R\setminus \tilde{N}_0)=0$ and for all $s\in  \R\setminus \tilde{N}_0$, 
\[
\sum_{i=1}^\infty {b}_i^2(s)<\infty. 
\]
Therefore we are allowed to define 
\[
   b(s)=\left\{
   \begin{aligned}
   	  	 &\sum_{i=1}^\infty {b}_i(s)e_i,
	          & \quad  
	           	   s\in&\R\setminus \tilde{N}_0; \\
	 &0, &  \quad	s\in&\tilde{N}_0 .
   \end{aligned}
\right.	 
\]

For all $\xi\in\H$, $A\in\B(\R)$, we have 
\[
\begin{aligned}
       \<\xi,a(A)\>
 &=\sum_{i=1}^\infty \<\xi,e_i\>\<e_i,{a}(A)\>
    =\sum_{i=1}^\infty \<\xi,e_i\>{a}(e_i,A)
 =\sum_{i=1}^\infty \<\xi,e_i\> 
	    \int_A {b}_i(s)\,\tilde{\gamma}(ds)\\
 &  =\left\<\xi,  \int_A \sum_{i=1}^\infty 
	    {b}_i(s)e_i\,\tilde{\gamma}(ds)\right\>
   =\left\<\xi,  \int_A {b}(s)\,\tilde{\gamma}(ds)\right\>.
\end{aligned}
\]
This shows that $b(s)$ is what we are looking for. 
\end{proof}

For every $A\in \B^0(\R)$,  we define 
\[
\begin{aligned}
\psi_A(\xi)&=-i\<a(A),\xi\>+\frac12 \<R(A)\xi,\xi\> \\
                  &\quad\quad - 
                  \int_{\H\setminus\{0\}} 
                  \left(
		          \e^{i\<\xi,x\>} - 1 -\frac{i\<\xi,x\>}{1+|x|^2} 
		 \right)\,\m(A,dx),\quad \xi\in\H.
\end{aligned}
\]
Thus for all $A\in\B^0(\R)$, 
\[
\hat{\mu}_A(\xi)=\exp(-\psi_A(\xi)),\quad \xi\in\H
\]
defines an infinitely divisible measure $\mu_A$ on $(\H,\B(\H))$ 
with 
\[
\mu_A=[a(A),R(A),\m(A,\cdot)].
\]
Moreover, it is obvious that for all $s<t$, 
\(
\mu_{t,s}=\mu_{(s,t]}
\)
and for all $A,B\in\B^0(\R)$ such that $A\cap B=\emptyset$, we have
\(
\mu_{A\cup B}=\mu_A*\mu_B. 
\)

Now we define a canonical control measure $\tilde{\pi}$ on $(\R,\B(\R))$ by 
\[
\tilde{\pi}(A):=\tilde{\gamma}(A)+\tilde{\pi}_R(A)+\tilde{F}(A),\quad A\in\B^0(\R).
\]
It is clear that all the measures $\tilde{\gamma}$, $\tilde{\pi}_R$ and $\tilde{F}$ are absolutely continuous with respect to $\tilde{\pi}$. Let 
\[
\tilde{h}^{(1)}(s):=\frac{d\tilde{\gamma}}{d\tilde{\pi}}(s),\quad
\tilde{h}^{(2)}(s):=\frac{d\tilde{\pi}_R}{d\tilde{\pi}}(s),\quad
\tilde{h}^{(3)}(s):=\frac{d\tilde{F}}{d\tilde{\pi}}(s),\quad s\in\R.
\]

For every $s\in\H$, set 
\[
\begin{aligned}
\rho_s(\xi)=
             &-i    
                \left<{b}(s)\tilde{h}^{(1)}(s),    \xi \right> 
		    +\frac12	\left< \left ( K(s)\tilde{h}^{(2)}(s) \right) \xi, \xi \right> \\
		    & - \int_{\H\setminus \{0\}}  \left(
		          \e^{i\<\xi,x\>} - 1 -\frac{i\<\xi,x\>}{1+|x|^2} 
		     \right) \tilde{h}^{(3)}(s){J}(s,dx),\quad \xi\in\H . 
\end{aligned}		     
\]
It is clear that 
$\rho_s(\cdot)$ is the characteristic exponent of some infinitely divisible measure  of the following form 
\[
  [{b}(s)\tilde{h}^{(1)}(s),
{K}(s)\tilde{h}^{(2)}(s),
\tilde{h}^{(3)}(s)J(s,\cdot)].
\]

Then we get the following theorem. 


\begin{theorem}\label{Thm:space-homo-process-rept}
Suppose that $(\mu_{t,s})_{t\geq s}$ is natural and Assumptions 
\ref{Ass:ContinuousMeasure'}, \ref{assumption:summable:total_var} hold. 
Then for all $A\in\B^0(\R)$ we have 
\[
     \hat{\mu}_A(\xi) 
     = \exp(-\psi_A(\xi))
     =\exp\left(-\int_A \rho_s(\xi)\,\tilde{\pi}(ds)\right),\quad \xi\in\H.
\]

In other words, 
for every $A\in\B^0(\R)$ we have
\[
\mu_A=
\left[
\int_A {b}(s)\tilde{h}^{(1)}(s)\,\tilde{\pi}(ds),
\int_A K (s)\tilde{h}^{(2)}(s)\,\tilde{\pi}(ds),
\int_A  \tilde{h}^{(3)}(s) {J}(s,dx) \,\tilde{\pi}(ds)
\right].
\]
In particular, we have 
\[
\hat{\mu}_{t,s}(\xi) 
     = \exp(-\psi_{t,s}(\xi))
     =\exp\left(-\int_s^t \rho_s(\xi)\,\tilde{\pi}(ds)\right),\quad \xi\in\H.
\]
and 
\[
\mu_{t,s}=
\left[
\int_s^t {b}(s)\tilde{h}^{(1)}(s)\,\tilde{\pi}(ds),
\int_s^t  K (s)\tilde{h}^{(2)}(s)\,\tilde{\pi}(ds),
\int_s^t  \tilde{h}^{(3)}(s) {J}(s,dx) \,\tilde{\pi}(ds)
\right].
\]
\end{theorem}

\begin{remark}
It is clear that measures $(\mu_{t,s})_{t\geq s}$ satisfying the convolution equation \eqref{Equ:mu_ts=mu_tr_mu_rs} and Assumption \ref{Ass:ContinuousMeasure'} is associated with an stochastic additive process (cf. the arguments in Subsection \ref{Subsec:additive_process}). 
  For stochastic continuous additive process, the representation of the characteristic function of the increments can be found for example in 
  \cite[Theorem 3.2.17]{Sko91}, \cite[Theorem II.5.2]{JS87} and \cite{Sat99}. 
Theorem \ref{Thm:space-homo-process-rept} above says that if we have in addition that $(\mu_{t,s})_{t\geq s}$ is natural  \ref{assumption:summable:total_var} hold, then the additive process is locally homogeneous (see  \cite[Section 3.4]{Sko91}). 
\end{remark}

\subsubsection{Differentiable condition}\label{Spectral_Rep:DiffCond}	
We shall consider the special case where the control measure $\pi$ is the Lebesgue measure. We need the following fact. 

\begin{lemma}\label{lemma:zero_one_side_derivative} 
Let $f$ be a continuous function on $[a,b]$. If for each $x\in(a,b)$ either the left derivative or the right derivative vanishes, then $f$ is constant. 
\end{lemma}

For the proof we refer to \cite[Theorem 1]{MV86} (or the references therein, e.g. \cite[Page 365 (3rd Ed.) or Page 341 (2nd Ed.)]{Hob57}). 

The main result of this subsection is the following proposition. 

\begin{prop}\label{Prop:Construction_ld}
Let $(\mu_{t,s})_{t\geq s}$ satisfy
\eqref{Equ:SkewMeasure}, or equivalently $\psi_{t,s}(\xi)=-\log \hat{\mu}_{t,s}(\xi)$ satisfy
 \eqref{Equ:SkewMeasure:Symbol} for all $\xi\in\H$. 
Suppose that for every  \(\xi\in \H\) and $t\geq s$, 
		\begin{enumerate}
\item The function \(s\mapsto \psi_{t,s}(\xi) \) is continuous and left differentiable at \(s=t\). Denote the left derivative by $-\lambda_t(\xi)$, i.e. 
			\begin{equation}\label{M_ss:Derivative}
			   -\lambda_t(\xi):=\left. \frac{d^-}{ds}\psi_{t,s}(\xi) \right\vert_{s=t}
			   	=\lim_{r\uparrow t}\frac{\psi_{t,r}(\xi)}{r-t}.
			\end{equation}
\item The function \(s\mapsto \lambda_s\bigl(U(t,s)^*\xi\bigr)\) is continuous.
		\end{enumerate}
	Then  for every $t\in\R$, $\lambda_t(\cdot)$ is 
	negative definite and Sazonov continuous on \H such that
	for every $t\geq s$, 
		\begin{equation}\label{M_ts:Integral}
			   \hat{\mu}_{t,s}(\xi)=\exp\left(-\int_s^t \lambda_r\bigl(U(t,r)^*\xi\bigr)\,dr\right),\quad \xi\in \H.
		\end{equation}
		\end{prop}

\begin{proof}
		For every \(\xi\in \H\) and \(r\leq t\), by \eqref{Equ:SkewMeasure:Symbol} we get 
          \[
		   \begin{aligned}
			  \frac{d^-}{dr} \psi_{t,r}(\xi) 
		       &= \lim_{r'\uparrow r} \frac{\psi_{t,r'}(\xi)-\psi_{t,r}(\xi)}{r'-r}   
		       = \lim_{r'\uparrow r} \frac{\psi_{t,r}(\xi)+\psi_{r,r'}(U(t,r)^*\xi)-\psi_{t,r}(\xi)}{r'-r}\\
		       &= \lim_{r'\uparrow r} \frac{\psi_{r,r'}(U(t,r)^*\xi)}{r'-r}
		       =-\lambda_r\bigl(U(t,r)^*\xi\bigr).  
		   \end{aligned}
		\] 
By the assumption, for every $\xi\in\H$, $r\mapsto \lambda_r(U(t,r)^*\xi)$,  $r\leq t$, is continuous. Hence we have
\begin{equation}\label{ZeroLeftDiff}
			\frac{d^-}{dr}\Phi_{t,r}(\xi)=0, \quad r\leq t, \ \xi\in\H, 
\end{equation}
where \[\Phi_{t,r}(\xi):=\psi_{t,r}(\xi)-\int_r^t  \lambda_{u}(U(t,u)^*\xi)\,du, \quad r\leq t, \ \xi\in\H.  \]
By Lemma \ref{lemma:zero_one_side_derivative} we get that 
$\Phi_{t,r}(\xi)$ is constant for every $r\in[s,t]$. 
But $\Phi_{t,t}(\xi)=0$.  So we also have $\Phi_{t,s}(\xi)=0$. It follows that 
$$\psi_{t,s}(\xi)=\int_s^t\lambda_r(U(t,r)^*\xi)\,dr.$$
Since $\psi_{t,s}(\xi)=-\log \hat{\mu}_{t,s}(\xi)$,  we obtain \eqref{M_ts:Integral}. 

From the negative definiteness and Sazonov continuity of $\psi_{t,s}(\cdot)$ we get the corresponding properties of $\lambda_t(\cdot)$.  
\end{proof}

\begin{remark}
For every $\xi\in\H$, the assumption that $s\mapsto \lambda_s(U(t,s)^*\xi)$ is continuous for $s\leq t$, is used to ensure that the map
\(s\mapsto \int_s^t \lambda_u(U(t,u)^*\xi)\,du\) is continuous and has (left)-derivative $-\lambda_s(U(t,s)^*\xi)$. 
This continuity assumption on $\lambda_\cdot(U(t,\cdot)^*\xi)$ holds if we assume that for every $\varepsilon>0$, $s\leq t$, and 
 for every bounded set $B\subset \H$, there exists a $\delta>0$ such that 
 $\sup_{x\in B}|\lambda_{s+h}(x)-\lambda_s(x)|<\varepsilon$ provided $|h|<\delta$ and $h<t-s$. Indeed, note that 
\[
\begin{aligned}
	& |\lambda_{s+h}(U(t,s+h)^*\xi)-\lambda_s(U(t,s)^*\xi)| \\
  \leq & |\lambda_{s+h}(U(t,s+h)^*\xi)-\lambda_{s}(U(t,s+h)^*\xi)| 
  		+ |\lambda_{s}(U(t,s+h)^*\xi)-\lambda_s(U(t,s)^*\xi)|. 
\end{aligned}
\] 
Hence $|\lambda_{s+h}(U(t,s+h)^*\xi)-\lambda_s(U(t,s)^*\xi)|$ can be made arbitrarily small, since the first term on the right hand side of the inequality above can be made small by the assumption that $s\mapsto \lambda_s(x)$ is continuous uniformly in $x$ on bounded sets; the second term can be made small by the strong continuity of $U(t,\cdot)$. 
\end{remark}

\begin{remark}\label{successful_representation}
Proposition \ref{Prop:Construction_ld} generalizes  \cite[Lemma 2.6]{BRS96} which 
study homogeneous generalized Mehler semigroups (see \eqref{GMehlerSemigroup}) using differentiability condition.
For the homogeneous case, there are some generalizations of  \cite[Lemma 2.6]{BRS96}.
Neerven \cite{Van00} relaxed the differentiability condition for general Gaussian Mehler semigroups on Banach space. 
Dawson et al. \cite[Theorem 2.1]{DL04} (see also \cite[Theorem 2.3]{DLSS04}) 
used entrance laws to characterize $\mu_t$ and dropped the differentiability condition for homogeneous generalized Mehler semigroups on Hilbert spaces.
For measure-valued skew convolution semigroups, the sufficiency and necessity of the representation appeared in \cite[Theorem 2]{Li96} and \cite[Theorem 3.1]{Li02} respectively, for the homogeneous case and the inhomogeneous case by using entrance laws. 
Proposition \ref{Prop:Construction_ld} can be seen as an attempt to 
 use entrance laws to characterize $\mu_{s,t}$. 
But we do not know how to find a natural measure $\pi$. 
\end{remark}

\section{Evolution systems of measures}\label{Sec:Evolution_Measures}
Let \((p_{s,t})_{t\geq s}\) be defined as in \eqref{2GMS} on a separable Hilbert space \(\H\) with \((U(t,s))_{t\geq s}\) 
and \((\mu_{t,s})_{t\geq s}\) satisfying \eqref{Equ:SkewMeasure}.

In general, for a family of  non-autonomous operators 
\((p_{s,t})_{t\geq s}\) on \(\H\), 
we cannot expect to have a stationary invariant measure for them. But we can try to look for a family of probability measures
 \((\nu_{t})_{t\in\R}\) on \((\H,\B(\H))\)
such that
\begin{equation}\label{Def:SysEvolutionMeasure}
  \int_\H p_{s,t}f(x)\,\nu_s(dx)=\int_\H f(x)\,\nu_t(dx),\quad s\leq t
\end{equation}
for all \(f\in \BB_b(\H)\). 
Such a family of probability measures is called an evolution 
system of measures for \((p_{s,t})_{t\geq s}\) (see \cite{DaPR08}). It can be regarded as a space-time invariant measure for \((p_{s,t})_{t\geq s}\). 
We note that evolution system of measures is called entrance law in \cite{Dyn89}.

We shall first show some basic properties and then study the existence and uniqueness of evolution system of measures.

\subsection{Basic properties}

\begin{lemma}\label{Lemma:SysEvolutionMeasure}
 A family of probability measures \((\nu_t)_{t\in\R}\) on \H is an evolution system of measures for \((p_{s,t})_{t\geq s}\) if and only if for every \(t\geq s\), 
\begin{equation}\label{Criteria:SysEvolutionMeasure:measure:equation}
       \mu_{t,s}*\bigl(\nu_s\circ U(t,s)^{-1}\bigr)
 =\nu_t,
\end{equation}
or equivalently, 
\begin{equation}\label{Criteria:SysEvolutionMeasure}
 \hat{\mu}_{t,s}(\xi) \hat{\nu}_s\bigl(U(t,s)^*\xi\bigr)
 =\hat{\nu}_t(\xi),\quad \xi\in \H.
\end{equation}
\end{lemma}
\begin{proof} 
We only need to note that for all $f\in B_b(\H)$ and $s\leq t$, we have 
\[
\begin{aligned}
	  &\int_\H p_{s,t}f(x)\,\nu_s(dx) \\
	=&\int_\H\int_\H f(U(t,s)x+y) \,\mu_{t,s}(dy)\nu_s(dx) \\
	=&\int_\H\int_\H f(x+y) \,\mu_{t,s}(dy)(\nu_s \circ U(t,s)^{-1}) (dx) \\
	=&\int_\H f(z)\, (\mu_{t,s }* (\nu_s \circ U(t,s)^{-1})(dz).
\end{aligned}
\]
\end{proof}

\begin{remark}
Equation \eqref{invariant_measure:equ:Xu} is the homogeneous version of Equation 
\eqref{Criteria:SysEvolutionMeasure:measure:equation}.
\end{remark}

\begin{prop}
Let $(\nu_t)_{t\in\R}$ be a family of probability measures 
on $(\H, \B(\H))$ such that \eqref{Criteria:SysEvolutionMeasure} holds for 
a family of probability measures $(\mu_{t,s})_{t\geq s}$ on  $(\H, \B(\H))$ and an evolution family of operators in $\H$. 
Suppose in addition that for all $\xi\in\H$ and all $t,s\in\R, t\geq s$,  we have $\hat{\nu}_t(U(t,s)^*\xi)\neq 0$.
Then $(\mu_{t,s})_{t\geq s}$ satisfies  
 \eqref{Equ:SkewMeasure}. 
\end{prop}

\begin{proof}
For any $t\geq r\geq s$, by \eqref{Criteria:SysEvolutionMeasure:measure:equation}
 and \eqref{Equ:mu*nu_T_inverse} we have
    \[
    \begin{aligned}
            \mu_{t,s}*\bigl(\nu_s\circ U(t,s)^{-1}\bigr) 
    & =\nu_t \\
 & = \mu_{t,r}*\bigl(\nu_r\circ U(t,r)^{-1}\bigr) \\
 & = \mu_{t,r}*\bigl(
             [
                   \mu_{r,s}*(\nu_s\circ U(r,s)^{-1} )
               ]
                \circ U(t,r)^{-1}\bigr) \\
 & = \mu_{t,r}*\bigl(
             [
                   \mu_{r,s}  
		                   	\circ U(t,r)^{-1} ] *
           [          (\nu_s\circ U(r,s)^{-1} ) 
				            \circ U(t,r)^{-1}       
               ]
                \bigr) \\
  &  =  \mu_{t,r}*
                \bigl(   \mu_{r,s}\circ U(t,r)^{-1} \bigr) * \bigl(\nu_s\circ U(t,s)^{-1} \bigr).
    \end{aligned}
    \]
 So for all $\xi\in\H$, we have
 \[
 	    \hat{\mu}_{t,s}(\xi) \cdot \hat{\nu}_s( U(t,s)^{*}\xi)
	 =\hat{\mu}_{t,r}(\xi) \cdot
 		 \hat{\mu}_{r,s}( U(r,s)^{*}\xi ) \cdot \hat{\nu}_s( U(t,s)^{*}\xi).
 \]   
Since $\hat{\nu}_t(U(t,s)^*\xi)\neq 0$ by assumption,  we have  
\[
 	    \hat{\mu}_{t,s}(\xi) 
	 =\hat{\mu}_{t,r}(\xi)
		\cdot   \hat{\mu}_{r,s}( U(r,s)^{*}\xi ).
 \]   
This proves \eqref{Equ:SkewMeasure}.
\end{proof}

We have parital result for the weak continuity of $(\nu_t)_{t\in\R}$ in $t$. 

\begin{prop}\label{Prop:stochastic_continuity:nu_t}
For any $t\in\R$, 
suppose that $\mu_{t+\varepsilon,t}$ converges weakly to $\delta_0$
as $\varepsilon\downarrow 0$. 
Then $\nu_{t+\varepsilon}$ converges weakly to $\nu_{t}$ as $\varepsilon \downarrow 0$. 
\end{prop}
\begin{proof}
By \eqref{Criteria:SysEvolutionMeasure:measure:equation} we have for all $\varepsilon>0$,
\begin{equation}\label{Criteria:SysEvolutionMeasure:measure:equation_t0s_new1}
  \mu_{t+\varepsilon,t}*\bigl(\nu_t\circ U(t+\varepsilon,t)^{-1}\bigr)
	 =\nu_{t+\varepsilon}. 
\end{equation}

By assumption we have
$\mu_{t+\varepsilon,t}\Rightarrow \delta_0$.  
By the bounded convergence theorem, it is easy to show  (cf. the proof of \eqref{Equ:mu_ts_U_t+epsilon})
\begin{equation}\label{equ:convergence_nu_t_U}
\nu_t\circ U(t+\varepsilon,t)^{-1}\Rightarrow \nu_t\quad 
\textrm{as}\ \varepsilon\downarrow 0. 
\end{equation}
Then by applying the first result of Lemma 
\ref{lmm:convolution:continuity} to \eqref{Criteria:SysEvolutionMeasure:measure:equation_t0s_new1}, 
the proof is finished immediately. 
\end{proof}

Similar to Theorem \ref{thm:mu_ts2process}, there exists certain stochastic process associated with $\nu_t$.

\begin{theorem}
	Suppose that Assumption \ref{Ass:ContinuousMeasure'} and   \eqref{Equ:SkewMeasure}, \eqref{Criteria:SysEvolutionMeasure:measure:equation} hold. Then there is a stochastic  process $(X_{t,-\infty})_{t\in\R}$, such that for every $t\in\R$, 
	 $X_{t,-\infty}$ is distributed as $\nu_t$.
\end{theorem}

\begin{proof} 
 Let $\Omega=\H^{(-\infty,\infty)}$ be the collection of all functions 
 $\omega=(\omega(t))_{t\in(-\infty,\infty)}$ from $(-\infty,\infty)$ into $\H$. 
 Let $X_t(\omega)=\omega(t)$, $t\in (-\infty,\infty)$, be the canonical process on $\Omega$. 
 Let $\F$ be the $\sigma$-algebra generated by cylinder sets on $\H^{(-\infty,\infty)}$. For any $n\in\N$, 
 $-\infty<t_1\leq t_2\leq \cdots\leq t_n<\infty$,  $B_j\in\B(\H)$, $j=1,2,\cdots,n$, define 
\begin{equation}\label{Equ:nu_t0tn:def}
\begin{aligned}
    &\nu_{t_1,t_2\cdots,t_n}(B_1\times B_2\times \cdots \times B_n)\\
   =&\int_{\H} \Eins_{B_1}(y_1) \nu_{t_1}(dy_1) \int_\H \Eins_{B_2}(U(t_2,t_1)y_1+y_2)\mu_{t_2,t_1}(dy_2) \\
     &\times \int_\H \Eins_{B_3}(U(t_3,t_1)y_1+U(t_3,t_2)y_2+y_3)\mu_{t_3,t_2}(dy_3) \times \cdots  \\
     &\times \int_\H \Eins_{B_n}(U(t_n,t_1)y_1+U(t_n,t_2)y_2+\cdots+U(t_n,t_{n-1})y_{n-1}+y_n)\mu_{t_n,t_{n-1}}(dy_n). 
%
   %
\end{aligned}
\end{equation}
Then, as in the proof of Theorem \ref{thm:mu_ts2process}, $\nu_{t_1,t_2,\cdots,t_n}$ is extended to a probability measure on $(\H^{\otimes n}, \B(\H^{\otimes n}))$ and it is easy to check that the family of probability measures $\{\nu_{t_1,t_2,\cdots,t_n}\}_{t_1\leq t_2\leq \cdots\leq t_n}$ satisfies the consistency condition.  Again, we shall only show the point by the following example. That is, we are going to prove 
$$\nu_{t_1,t_2,t_3,t_4}(B_1\times \H\times B_3\times B_4 )=\nu_{t_1,t_3,t_4}(B_1 \times B_3 \times B_4)$$
holds for  all 
$-\infty<t_1\leq t_2\leq t_3\leq t_4<\infty$, and  all 
 $B_j\in\B(\H)$, $j=1,2,3,4$ with $B_2=\H$. 
By  \eqref{Equ:SkewMeasure}  we have
\[
\begin{aligned}
	       &  \nu_{t_1,t_2,t_3,t_4}(B_1\times \H\times B_3 \times B_4) \\
	     =&\int_\H 
	     			 \Eins_{B_1}(y_1) \nu_{t_1}(dy_1)
			\int_\H \mu_{t_4,t_3}(dy_4)
					\int_\H\int_\H 					    
					      \Eins_{B_3}(U(t_3,t_1)y_1+U(t_3,t_2)y_2+y_3) 
			  \\
	     	&		  \times   \Eins_{B_4 }(U(t_4,t_1)y_1+U(t_4,t_2)y_2+U(t_4,t_3)y_3+y_4)
						     \mu_{t_2,t_1}(dy_2) 
							\mu_{t_3,t_2}(dy_3)
  					      \\
	     =&\int_\H 
	     			 \Eins_{B_1}(y_1) \nu_{t_1}(dy_1)
			\int_\H \mu_{t_4,t_3}(dy_4)
					\int_\H\int_\H 					    
					      \Eins_{B_3}(U(t_3,t_1)y_1+y_2+y_3) 
			  \\
	     	&			  \times   \Eins_{B_4 }(U(t_4,t_1)y_1+U(t_4,t_3)(y_2+y_3)+y_4)
						     (\mu_{t_2,t_1}\circ U(t_3,t_2)^{-1})(dy_2) 
							\mu_{t_3,t_2}(dy_3)
  					      \\
	     =&\int_\H 
	     			 \Eins_{B_1}(y_1) \nu_{t_1}(dy_1)
			\int_\H \mu_{t_4,t_3}(dy_4)
					\int_\H\int_\H 					    
					      \Eins_{B_3}(U(t_3,t_1)y_1+z) 
			  \\
	     	&			  \times   \Eins_{B_4 }(U(t_4,t_1)y_1+U(t_4,t_3)z+y_4)
						   \left[  (\mu_{t_2,t_1}\circ U(t_3,t_2)^{-1})* 
							\mu_{t_3,t_2} \right](dz)
  					      \\					      
	     =&\int_\H 
	     			 \Eins_{B_1}(y_1) \nu_{t_1}(dy_1)
			\int_\H \mu_{t_4,t_3}(dy_4)
					\int_\H\int_\H 					    
					      \Eins_{B_3}(U(t_3,t_1)y_1+z) 
			  \\
	     	&			  \times   \Eins_{B_4 }(U(t_4,t_1)y_1+U(t_4,t_3)z+y_4)
						   \mu_{t_3,t_1}(dz)
  					      \\	
%
=&\nu_{t_1,t_3,t_4}(B_1 \times B_3 \times B_4).	      
\end{aligned}
\]

Therefore, by Kolmogorov's extension theorem there is a unique probability measure $\P$ on $\F$ such that 
for all 
$-\infty<t_1\leq t_2\leq \cdots\leq t_n<\infty$,  and $B_j\in\B(\H)$, $j=1,2,\cdots,n$, we have 
\begin{equation}\label{P2mu_t0tn'}
\begin{aligned}
 & \P(X_{t_1,-\infty}\in B_1,X_{t_2,-\infty}\in B_2, \cdots, X_{t_n,-\infty}\in B_n) \\
 = & \nu_{t_1,t_2,\cdots,t_n}(B_1\times B_2\times \cdots\times B_n).
\end{aligned} 
\end{equation}
In particular, $X_{t,-\infty}$ has the distribution $\nu_t$.
\end{proof}

Similar to Example \ref{Example:X_t:process:mu_ts}, we have the following example. 
\begin{example}\label{Example:stoch_integral:nu_t}
For any $t\in\R$, 
    consider a stochastic process $(X_{t,-\infty})_{t\in\R}$ , defined by
    \[
     	X_{t,-\infty}:=\int_{-\infty}^{t} U(t,\sigma)\,dZ_\sigma,
    \]
where $U$ and $Z$ are the same as in \eqref{NonAutoOUprocess}. Let $\mu_{t,s}$ be the distribution of 
$$X_{t,s}:=\int_{s}^{t} U(t,\sigma)\,dZ_\sigma$$ 
and $\nu_t$ the distribution of $X_{t,-\infty}$. 

Note that  for all $t\geq s$, 
\[
	X_{t,-\infty}  = U(t,s) X_{s,-\infty} + X_{t,s}.
\]
Hence we get identity 
\eqref{Criteria:SysEvolutionMeasure:measure:equation} immediately 
since $X_{s,-\infty}$ is independent of  $X_{t,s}$. This proves that  $(\nu_t)_{t\in\R}$ is an evolution system of measures for the transition function $(p_{s,t})_{t\geq s}$ of $X(t,s,x):=U(t,s)x+X_{t,s}$, $x\in\H,\,t\geq s$.
\end{example}

Concerning the infinite divisibility of $\nu_t$, we have the following simple result. 
\begin{prop}
Suppose that Assumption \ref{Ass:ContinuousMeasure'} holds. 
If for some $s_0\in\R$, $\nu_{s_0}$ is infinitely divisible, then $\nu_t$ is infinitely divisible for all $t\geq s_0$. 
\end{prop} 

\begin{proof}
Since $\nu_{s_0}$ is infinitely divisible, 
for any $n\in\N$, there is some probability measure $\nu_{s_0}^{(n)}$ on $(\H,\B(\H))$ such that 
$\nu_{s_0}=(\nu_{s_0}^{(n)})^{*n}$ .
So by \eqref{Equ:mu*nu_T_inverse}, we have $$\nu_{s_0}\circ U(t,s_0)^{-1}=(\nu_{s_0}^{(n)}\circ U(t,s_0)^{-1})^{*n}.$$
for all $t\geq s_0$. 
This proves that 
$\nu_{s_0}\circ U(t,s_0)^{-1}$ is also infinitely divisible. 
By Theorem \ref{Thm:mu_ts_IID}, $\mu_{t,s_0}$ is infinitely divisible. Hence 
by \eqref{Criteria:SysEvolutionMeasure:measure:equation}, we get that $\nu_t$ is the convolution of two infinitely divisible measures. So $\nu_t$ itself is also an infinitely divisible measure. 
\end{proof}

For any two probability measures $\mu$ and $\nu$, we say $\mu$ is a \emph{factor} of $\nu$ if $\nu=\mu*\sigma$ for some probability measures $\sigma$. 

\begin{theorem}\label{Thm:SysEvolutionMeasure:Structure}
Suppose that \((\nu_t^{(1)})_{t\in\R}\) is an evolution system of measures for \((p_{s,t})_{t\geq s}\).
Let \((\nu_t^{(2)})_{t\in\R}\) be another system of  probability measures and assume that for any $t\in\R$,
$\nu_t^{(1)}$ is a factor of $\nu_t^{(2)}$ satisfying 
\begin{equation}\label{Equ:nu_t_sigma_t}
    \nu_t^{(2)}=\nu_t^{(1)}*\sigma_t
\end{equation}
with probability measure \(\sigma_t\) on \((\H,\B(\H))\) satisfying
 \[ \sigma_t = \sigma_s\circ U(t,s)^{-1}.   
 \]
Then \((\nu_t^{(2)})_{t\in\R}\) is also an evolution system of  measures for \((p_{s,t})_{t\geq s}\). 
\end{theorem}

\begin{proof}
     For every \(\xi\in\H\), by \eqref{Criteria:SysEvolutionMeasure}
and \eqref{Equ:nu_t_sigma_t}, we have 
	\[
		\begin{aligned}
			     \hat{\nu}_t^{(2)}(\xi)
                        &=\hat{\nu}_t^{(1)}(\xi)\hat{\sigma}_t(\xi)
			 = \hat{\mu}_{t,s}(\xi) \hat{\nu}_s^{(1)}\bigl(U(t,s)^*\xi\bigr)\hat{\sigma}_t(\xi)\\
			&= \hat{\mu}_{t,s}(\xi) \hat{\nu}_s^{(1)}\bigl(U(t,s)^*\xi\bigr)\hat{\sigma}_s( U(t,s)^* \xi)
			= \hat{\mu}_{t,s}(\xi) \hat{\nu}_s^{(2)}\bigl(U(t,s)^*\xi\bigr).
		\end{aligned}
	\]
	Hence the assertion follows by Lemma \ref{Lemma:SysEvolutionMeasure}. 
\end{proof}
 
\subsection{Existence and uniqueness} 
Assume that for every $t\geq s$, $\mu_{t,s}$ is infinitely divisible and  has the form 
\[\mu_{t,s}=[a_{t,s},R_{t,s}, \m_{t,s}],\]
 where $a_{t,s}\in\H$, $R_{t,s}$ is a non-negative definite, self-adjoint  trace class operator on \H, and $\m_{t,s}$ is a L\'evy measure on \H.  

By \eqref{Equ:Structure_char_triplet}, we get that for every fixed \(t\in\R\), 
\((\m_{t,s})_{ t\geq s}\) is a family of L\'evy measures decreasing in $s$ in the sense that for all $A\in\B(\H\setminus\{0\})$, we have 
\[ \m_{t,s}(A)\geq \m_{t,r}(A) \]
 for all $s\leq r\leq t$. 
It allows us to define $\m_{t,-\infty}$ for every $t\in\R$ by setting
$\m_{t,-\infty}(\{0\})=0$ and
\[
   \m_{t,-\infty}(A)=\lim_{s\rightarrow -\infty}\m_{t,s}(A),
\quad A\in \B(\H\setminus \{0\}). 
\]

From \eqref{Equ:Structure_char_triplet} we also get that  for every 
$x\in \H$ and $t\in\R$, $\<R_{t,s}x,x\>$ is decreasing in $s$. 
More precisely, for every $x\in \H$ and $s\leq r\leq t$, we have
\[
   \<R_{t,s}x,x\>=\<R_{t,r}x,x\> + \<R_{r,s} U(t,r)^*x, U(t,r)^*x\>
                \geq \<R_{t,r}x,x\>.
\]
From the inequality above we obtain that 
for every $s\leq r\leq t$
\[
        \tr R_{t,r}\leq \tr R_{t,s}.
\]
Therefore, the limit $\lim_{s\to -\infty}\<R_{t,s}x,x\>$ which might be infinity exists. 

We shall use condition 
\[\sup_{s\leq t}\tr R_{t,s}<\infty.\]
Then for every $t\in\R$, there exists a constant $C_t>0$ such that
 $$\sup_{s<t}\<R_{t,s}x,x\>\leq C_t|x|^2$$ for every $x\in \H$.
Therefore, by the polarization identity, we see that for every $t\in\R$, $x,y\in\H$, 
the limit $\lim_{s\to -\infty}\<R_{t,s}x,y\>$ exists. Fixing $x\in\H$ and letting $y\in\H$ vary, we get a functional 
$\lim_{s\to -\infty}\<R_{t,s}x,\cdot\>$. 
So by Riesz's  representation theorem we obtain that for every $x\in \H$, there exists an element $x_t^*\in\H$ 
for every $t\in \R$ such that for all $y\in\H$, 
\[\lim_{s\to -\infty}\<R_{t,s}x,y\>=\<x_t^*,y\>.\]
By the property of $R_{t,s}$, we see that the mapping from $x$ to $x_t^*$ is a trace class operator and we denote it by $R_{t,-\infty}$. 
That is, 
for every $t\in\R$ there is a trace class operator $R_{t,-\infty}$ on \H such that 
\[
    \<R_{t,-\infty}x,y\>
=\lim_{s\rightarrow -\infty} \<R_{t,s}x,y\>,\quad x,y\in\H.
\]

\begin{theorem}\label{EMeasure:Levy2GMS}
Suppose that for every $t\in\R$, the following three hypothesises hold:
\begin{description}
\item[(H1)]  \(\sup_{s\leq t}\tr R_{t,s}<\infty\);
\item[(H2)] $\sup_{s\leq t} \int_\H (1\wedge |x|^2)\, \m_{t,s}(dx)<\infty $;
\item[(H3)]
	\(a_{t,-\infty}:=\lim_{s\rightarrow -\infty} a_{t,s}
        \).
\end{description}        
Then for every $t\in\R$, $\m_{t,-\infty}$ is a L\'evy measure, 
\(R_{t,-\infty}\) is a non-negative definite, self-adjoint  trace class operator 
such that
\[\tr R_{t,-\infty}=\sup_{s\leq t}R_{t,s}<\infty.\]
Moreover, the system of measures $(\nu_t)_{t\in\R}$ given by
\[\nu_t=[a_{t,-\infty}, R_{t,-\infty}, \m_{t,-\infty}], \quad t\in\R\] 
is an evolution system of measures for \((p_{s,t})_{t\geq s}\).
\end{theorem}

\begin{proof}
 Suppose that (H1), (H2) and (H3) hold. 
 Apparently,  for every $t\in\R$, 
 \(R_{t,-\infty}\) is a non-negative definite, self-adjoint  trace class operator satisfying
 \[\tr R_{t,-\infty}=\sup_{s\leq t}R_{t,s}<\infty.\]
For each $t\in\R$, 
	\[
	 \begin{aligned}
	    \int_\H (1\wedge |y|^2)\,\m_{t,-\infty}(dy) 
	    = \sup_{s\leq t} \int_\H (1\wedge |x|^2)\, \m_{t,s}(dx)<\infty.
	 \end{aligned}
	\]
 This shows that  \(\m_{t,-\infty}\) is a L\'evy measure. 

Now we show that $(\nu_t)_{t\in\R}$ is an evolution system of measures. 
By \eqref{Equ:SkewMeasure}, for every $t\geq s\geq r$, we have 
\begin{equation}\label{Equ:PfEM1.2}
 \begin{aligned}     
\mu_{t,s}*\bigl(\mu_{s,r}\circ U(t,s)^{-1}\bigr)
   =\mu_{t,r}.
 \end{aligned}
\end{equation}
Note that 
\(
       \mu_{t,s} 
    =  [a_{t,s},R_{t,s},\m_{t,s}]
\) converges weakly to 
 \(  [a_{t,-\infty},R_{t,-\infty},\m_{t,-\infty}]=\nu_t\) as $s\rightarrow -\infty$ (cf. \cite[Lemma 3.4]{FR00}). 
Hence letting $r\rightarrow -\infty$ in both sides of \eqref{Equ:PfEM1.2} we obtain
\[
\mu_{t,s}* \big(\nu_{s}\circ U(t,s)^{-1}\bigr) = \nu_t.
\]
This proves that $(\nu_t)_{t\in\R}$ is an evolution system of measures for 
\((p_{s,t})_{t\geq s}\) by Lemma \ref{Lemma:SysEvolutionMeasure}.
\end{proof}

The following theorem is the converse to  Theorem \ref{EMeasure:Levy2GMS}.

\begin{theorem}\label{EMeasure:Levy2GMS:Converse}   
Let  \((\tilde{\nu}_t)_{t\in\R}\) be an evolution system of measures for \((p_{s,t})_{t\geq s}\). Then 
\begin{enumerate}
\item Hypotheses (H1) and (H2) hold and 
\begin{equation}\label{Equ:Convergence:R_ts_M_ts_to_R_t_infty_M_t_infty}
 [0,R_{t,s},\m_{t,s}] \Rightarrow  [0,R_{t,-\infty},\m_{t,-\infty}],\quad s\to -\infty.
\end{equation}

\item 
There exists some probability measure $\tilde{\sigma}_t$ such that
  \begin{equation}\label{Equ:converge:a:nu:U_ts}
  	\delta_{a_{t,s}}*\bigl(\tilde{\nu}_{s}\circ U(t,s)^{-1}\bigr) \Rightarrow \tilde{\sigma}_t 
   \end{equation} 
 as $s\rightarrow -\infty$. 
Moreover,  for every $t\in\R$,
\begin{equation}\label{Equ:tilde_nu_R_t_infty_convolution_tilde_sigma}
		\tilde{\nu}_t=[0,R_{t,-\infty},\m_{t,-\infty}]*\tilde{\sigma}_t, \quad t\in\R.
\end{equation}
\item  Assume in addition that the following hypotheses holds 
          \begin{description}
        \item[(H4)] For every $t\in\R$, 
                 \(\tilde{\nu}_s\circ U(t,s)^{-1}\Rightarrow \sigma_t\)
	       as $s\rightarrow -\infty$.
           \end{description}
If 
for every $t\in\R$, $\sigma_t$ is an infinitely divisible distribution, 
then the limit in (H3) exists and $\nu_t$ is a factor of $\tilde{\nu}_t$:
 \begin{equation}\label{EvolutionMeasures:Connection}
  	\tilde{\nu}_t=\nu_t*\sigma_t ,\quad t\in\R.
 \end{equation} 
Moreover, 
\begin{equation}\label{Sigma_ts}
  \sigma_t=\sigma_s\circ U(t,s)^{-1},\quad t\geq s.
\end{equation} 
Especially, if $\sigma_t\equiv \delta_0$ for all $t\in\R$, then we have $\tilde{\nu}_t=\nu_t$, 
\(t\in\R\). That is, the evolution system of measures is unique. 
\item On the other hand, if the limit in (H3) exists, then the limit in (H4) exists,  and hence \eqref{EvolutionMeasures:Connection}, \eqref{Sigma_ts} hold. 
\end{enumerate}
\end{theorem}

\begin{proof} (1)
Since $(\tilde{\nu_t})_{t\in\R}$ is an evolution system of measures for $(p_{s,t})_{t\geq s}$, by Lemma \ref{Lemma:SysEvolutionMeasure} we have 
for every $t\geq s$, 
\begin{equation}\label{Equ:PfEM1.5}
 \begin{aligned}
   \tilde{\nu_t}
=\mu_{t,s}* \bigl(\tilde{\nu}_{s}\circ U(t,s)^{-1}\bigr) 
&= [a_{t,s},R_{t,s},\m_{t,s}] *\bigl(\tilde{\nu}_{s}\circ U(t,s)^{-1}\bigr) \\
&=\delta_{a_{t,s}}*N_{R_{t,s}}*M_{t,s}*\bigl(\tilde{\nu}_{s}\circ U(t,s)^{-1}\bigr).
 \end{aligned}
\end{equation}
Here we have set 
$$N_{R_{t,s}}:=[0,N_{R_{t,s}},0],\quad 
	M_{t,s}:=[0,0,\m_{t,s}].
$$
Consider $s=-n, \ n\in\N$, for \eqref{Equ:PfEM1.5}. The sequence 
\((\delta_{a_{t,-n}}*N_{R_{t,-n}}*M_{t,-n})_{n\geq 1}\) is shift compact by \cite[Theorem III.2.2]{Par67} (see Theorem \ref{Theorem_Par67_III_2_1}), i.e. there exists a sequence  $(y_{t,-n})_{n\geq 1}$ in $\H$ for every  $t\in\R$ such that the sequence 
\[\delta_{y_{t,-n}}*(\delta_{a_{t,-n}}*N_{R_{t,-n}}*M_{t,-n})
=[y_{t,-n}+a_{t,-n}, R_{t,-n}, \m_{t,-n}]
\]
is weakly relatively compact. 
It implies that (see \cite[Theorem VI.5.3]{Par67}) 
\begin{equation}\label{Equ:sup_m_tn_x_geq_1}
	\sup_{n\in\N} \m_{t,-n}(\{x\in\H\colon |x|\geq 1\})<\infty.
\end{equation}
and
\begin{equation}\label{Equ:sup_R_tn_m_tn_x_leq_1}
\begin{aligned}
\sup_n\left( \tr R_{t,-n}+\int_{|x|< 1}|x|^2\,\m_{t,-n}(dx)\right)<\infty.
\end{aligned}
\end{equation}

From \eqref{Equ:sup_R_tn_m_tn_x_leq_1} we have 
\[
	\sup_{s\leq t}\tr R_{t,s}<\infty.
\]
Combining \eqref{Equ:sup_m_tn_x_geq_1} and \eqref{Equ:sup_R_tn_m_tn_x_leq_1} we have 
\[
	 \sup_{s\leq t} \int_\H (1\wedge |x|^2)\, \m_{t,s}(dx)<\infty . 
\]

So we can define naturally a L\'evy measure $\m_{t,-\infty}$ and a trace class operator $R_{t,-\infty}$ for each $t\in\R$. Hence by Theorem \ref{EMeasure:Levy2GMS} we have \eqref{Equ:Convergence:R_ts_M_ts_to_R_t_infty_M_t_infty}.
This proves (1).

(2) 
By \eqref{Equ:PfEM1.5} we have for all $t\geq s$
 \begin{equation*}
 \begin{aligned}
   \tilde{\nu_t}=
\mu_{t,s}* \bigl(\tilde{\nu}_{s}\circ U(t,s)^{-1}\bigr)  =
\left(\delta_{a_{t,s}}*\bigl(\tilde{\nu}_{s}\circ U(t,s)^{-1}\bigr)\right)*N_{R_{t,s}}*M_{t,s}.
 \end{aligned}
\end{equation*}
On the other hand, we have shown in \eqref{Equ:Convergence:R_ts_M_ts_to_R_t_infty_M_t_infty} that 
$N_{R_{t,s}}*M_{t,s}$ converges weakly to an infinitely divisible measure 
$[0,R_{t,-\infty},\m_{t,-\infty}]$ as $s\to -\infty$. 
Therefore 
by Corollary \ref{cora:convolution:continuity''},
measure $\delta_{a_{t,s}}*\bigl(\tilde{\nu}_{s}\circ U(t,s)^{-1})$  converges weakly as $s\to -\infty$. This proves 
\eqref{Equ:converge:a:nu:U_ts} and \eqref{Equ:tilde_nu_R_t_infty_convolution_tilde_sigma}. 

%

%

%





(3) Applying Corollary \ref{cora:convolution:continuity''}, 
and using the hypothesis (H4) and that 
$\sigma_t$ is infinitely divisible for every $t\in\R$, 
it follows from 
\eqref{Equ:converge:a:nu:U_ts}, 
that the limit of $a_{t,s}$ as $s\to -\infty$ exist. So (H3) holds and 
\[
\tilde{\sigma}_t=\delta_{a_{t,-\infty}}*\sigma_t.
\] 
By\eqref{Equ:tilde_nu_R_t_infty_convolution_tilde_sigma} we get for every $t\in\R$, 
\[
	\tilde{\nu}_t = [0,R_{t,-\infty},\m_{t,-\infty}]*\tilde{\sigma}_t 
	=	 [a_{t,-\infty},R_{t,-\infty},\m_{t,-\infty}]*  \sigma_t
	=
	\nu_t*\sigma_t .
\]
This proves \eqref{EvolutionMeasures:Connection}. 

Now we show  \eqref{Sigma_ts}. 
For every \(r\leq s\leq t\) and \(\xi\in\H\), we have 
\begin{equation}\label{Sigma_EM}
 \hat{\tilde{\nu}}_r(U(t,r)^*\xi)
=\hat{\tilde{\nu}}_r\bigl(U(s,r)^*U(t,s)^*\xi\bigr).
\end{equation}
Letting \(r\rightarrow -\infty\) in \eqref{Sigma_EM}, we get 
\(\hat{\sigma}_t(\xi)=\hat{\sigma}_s\bigl(U(t,s)^*\xi\bigr)\)
by ({H4}). This is equivalent to \eqref{Sigma_ts}.

(4) 
If (H3) holds, then $\delta_{a_{t,s}}\Rightarrow \delta_{a_{t,-\infty}}$ as $s$ tends to $-\infty$. 
Note that any dirac measure is infinitely divisible. 
So by applying
Corollary \ref{cora:convolution:continuity''} to
\eqref{Equ:converge:a:nu:U_ts}, 
we get that the limit in (H4) exist. Therefore, the assertions in (3) hold.   
\end{proof}

\begin{remark}
 Similarly, the invariant measure for $(p_t)_{t\geq 0}$ defined in \eqref{GMehlerSemigroup} is of the form $\nu*\mu_\infty$, where $\nu$
 is a measure on \H that is invariant under the action of the semigroup $T_t$, and $\mu_\infty$ is the centered Gaussian measure with covariance operator
 $Q_\infty$ which is the proper limit of the covariance operator of $\mu_t$.  We refer to \cite[Theorem 5.22]{Hai09} for details. 
\end{remark}

By applying Theorem \ref{EMeasure:Levy2GMS:Converse}   we have the following result on the uniqueness. It is a generalization of \cite[Theorem 3.12]{Woo09}.  

\begin{cora}\label{Cora:Evolu-Measure-uniquness:tight}
Let (H1),(H2) hold and $(\tilde{\nu})_{t\in\R}$ be an evolution system of measure for $(p_{s,t})_{t\geq s}$. 
Suppose that for every $t\in\R$, there is a sequence $(s_n)_{n\geq 1}$ 
satisfying $s_n\leq t$ for all $n\geq 1$ and $s_n\to -\infty$ as $n\to\infty$ such that the following conditions are fulfilled
\begin{enumerate}
\item There exist some constant $M, \omega>0$ such that
	 \begin{equation}\label{Equ:Utsn<=Me_omega}
		\|U(t,s_n)\|\leq M\e^{-\omega(t-s_n)}.
	 \end{equation}
\item	The sequence of probability measures $(\tilde{\nu}_{s_n})_{n\geq 1}$ 
	is uniformly tight.  
\end{enumerate}
Then 
 (H3) holds and $(\tilde{\nu}_t)_{t\in\R}=(\nu_t)_{t\in\R}$.
That is, $(\nu_t)_{t\in\R}=([a_{t,-\infty}, R_{t,-\infty}, \m_{t,-\infty} ])_{t\in\R}$
is the unique evolution system of measures for $(p_{s,t})_{t\geq s}$.
\end{cora}
\begin{proof}
    By the proof of Theorem \ref{EMeasure:Levy2GMS:Converse}, it is sufficient to 
show that for a sequence 
$(s_n)_{n\geq 1}$ 
satisfying $s_n\leq t$ for all $n\geq 1$ and $s_n\to -\infty$ as $n\to\infty$
such that 
\begin{equation}\label{Equ:convergence_nuU_sn}
	\hat{\tilde{\nu}}_{s_n}\circ U(t,s_n)^{-1}\Rightarrow 
	\delta_0,\quad 
	\quad \textrm{as}\ n\to \infty. 
\end{equation}

Let $\varepsilon,\eta>0$ be arbitrary. Because  
$(\tilde{\nu}_{s_n})_{n\geq 1}$ 
is uniformly tight, there is a compact set $K_{\eta}\subset \H$ such that 
for all $n\geq 1$, 
\begin{equation}\label{Equ:nu_sn_uniform_tight}
	\tilde{\nu}_{s_n}(\H\setminus K_{\eta})<\eta.
\end{equation}
Set for all $n\geq 1$
\[
	C_n=M^{-1}\e^{\omega(t-s_n)}.
\]
It is clear that $C_n\to\infty$ as $n\to \infty$. 
Hence there exists some $N_0>0$ such that for all $n\geq N_0$, the compact set $K_\eta$ is contained in 
\[
	K_n:=\{x\in\H\colon 
					|x|\leq \varepsilon C_n
	\}.
\]
So 
\begin{equation}\label{Equ:nu_sn_Kn<eta}
     \tilde{\nu}_{s_n}(\{\H\setminus K_n\})
     \leq \tilde{\nu}_{s_n}(\{\H\setminus K_\eta\}) 	<\eta,
     \quad n\geq N_0.
\end{equation}

Because 
\[
	\| U(t,s_n)\|\leq M\e^{-\omega (t-s_n)}=1/C_n, 
\]
by applying \eqref{Equ:nu_sn_Kn<eta} we obtain  for all 
$n\geq N_0$ (cf. Lemma \ref{Lemma:measure:contraction})
\[
\begin{aligned}
         & \tilde{\nu}_{s_n}
         \circ U(t,s_n)^{-1}(\{x\in\H\colon |x|>\varepsilon\}) \\
       =&\tilde{\nu}_{s_n}(\{x\in\H\colon | U(t,s_n) x|>\varepsilon\}) \\
   \leq& \tilde{\nu}_{s_n} (\{x\in\H\colon |x|>\varepsilon C_n)\\
   =& \tilde{\nu}_{s_n}(\{\H\setminus K_n\}) 
    <\eta.
\end{aligned}	
\]
By Lemma \ref{lmm:Ass:ContinuousMeasure'} we obtain \eqref{Equ:convergence_nuU_sn}. Therefore by Theorem \ref{EMeasure:Levy2GMS:Converse} we get 
 (H3) and hence $(\tilde{\nu}_t)_{t\in\R}=(\nu_t)_{t\in\R}$ is the unique evolution system of measures for $(p_{s,t})_{t\geq s}$.
\end{proof}

The advantage of using sequence in Corolally \ref{Cora:Evolu-Measure-uniquness:tight} can be seen from the following corollary.

\begin{cora}
Assume that
the evolution function $(U(t,s))_{t\geq s}$ and the family of probability measures  $(\mu_{t,s})_{t\geq s}$ are $T$-periodic for some $T>0$. That is, for every $t\geq s$, 
\[U(t+T,s+T)=U(t,s), \quad \mu_{t+T,s+T}=\mu_{t,s}.\]
Assume that (H1) and (H2) hold and 
 there exist some constants $M, \omega>0$ such that
$$\|U(t,s)\|\leq M\e^{-\omega(t-s)}.$$
 Then $([a_{t,-\infty}, R_{t,-\infty}, \m_{t,-\infty} ])_{t\in\R}$ is the unique evolution system of measures with period $T$ for $(p_{s,t})_{t\geq s}$.
\end{cora}

\begin{proof}
Let $(\tilde{\nu}_t)_{t\in\R}$ be any 
evolution system of measures with period $T$.
Take any $s_0\leq t$ and set $s_n=s_0-nT$ for all $n\geq 1$. Then 
for all $n\geq 1$, $\tilde{\nu}_{s_n}=\tilde{\nu}_{s_0}$. It is obvious that 
$(\tilde{\nu}_{s_n})_{n\geq 1}$ 
	is uniformly tight.  
So Corollary \ref{Cora:Evolu-Measure-uniquness:tight}
applies. Therefore $\nu_t:=[a_{t,-\infty}, R_{t,-\infty}, \m_{t,-\infty} ]$ exists for all $t\in\R$. Moreover it is easy to show that $(\nu_t)_{t\in\R}$ is 
$T$-periodic.
\end{proof}

\begin{example}\label{Exa:semi-Levy-process}
Let us use the framework in Example \ref{Example:stoch_integral:nu_t}.
Assume that for some constant $T>0$ we have 
\begin{enumerate}
\item $(U_{t,s})_{t\geq s}$ is $T$-periodic evolution family of operators on \H. 
\item $(Z_t)_{t\in\R}$ be a semi-L\'evy process with period $T>0$. That is, $(Z_t)_{t\in\R}$ is a stochastic continuous additive process such that  for all $t\geq s$, $Z_{t+T}-Z_{s+T}$ has the same distribution with $Z_t-Z_s$ (we refer to \cite{MS03} for more details). 
\end{enumerate}
 It is clear that   $(\mu_{t,s})_{t\geq s}$ is $T$-periodic. Moreover, the evolution system of measures $(\nu_t)_{t\in\R}$ for $(p_{s,t})_{t\geq s}$ is also $T$-periodic. So if \eqref{Equ:Utsn<=Me_omega} holds, then 
  by Corollary \ref{Cora:Evolu-Measure-uniquness:tight} 
 $(\mu_{t,s})_{t\geq s}$ is the unique $T$-periodic evolution system of measures.
 \end{example}

 Clearly Example \ref{Exa:semi-Levy-process} above generalizes the 
 stochastic equations with time dependent periodic coefficients driven by Gaussian and L\'evy processes  
considered in \cite{DaPL07} and \cite{Kna09} respectively.


\section{Harnack inequalities and applications}\label{Sec:Harnack inequalities}
Let \((p_{s,t})_{t\geq s}\) be defined as in \eqref{2GMS}.  That is, 
$p_{s,t}f(x)=(\mu_{t,s}*\delta_{U(t,s)x})f$ for every $x\in\H$ and $f\in\BB_b(\H)$.
Suppose that for all $t\geq s$, $\mu_{t,s}=[a_{t,s},R_{t,s},\m_{t,s}]$  is an infinitely divisible measure on $(\H,\B(\H))$ satisfying \eqref{Equ:SkewMeasure}.

For each $t\geq s$, set
$$\mu_{t,s}^g=[0,R_{t,s},0],\quad \mu_{t,s}^j=[a_{t,s},0,\m_{t,s}],$$
and for every $f\in\BB_b(\H)$, $x\in\H$, set 
$$p_{s,t}^g f(x):=(\mu_{t,s}^g * \delta_{U(t,s)x}) (f)=\int_\H f(U(t,s)x+y)\,\mu_{t,s}^g(dy), $$
$$p_{s,t}^j f(x):=(\mu_{t,s}^j *\delta_{x}) (f)=\int_\H f(x+y)\,\mu_{t,s}^j(dy).$$

With these notations, we have the following decomposition for $p_{s,t}$.

\begin{prop}\label{Lem:Decomposition}
  For every $t\geq s$, $x\in\H$ and $f\in\BB_b(\H)$,  
  $$p_{s,t}f(x)=p_{s,t}^g(p_{s,t}^j )f(x).$$
\end{prop}

\begin{proof} Note that $\mu_{t,s}=\mu_{t,s}^g*\mu_{t,s}^j$. Then we get 
\[
\begin{aligned}
 p_{s,t}f(x)=&(\mu_{t,s}*\delta_{U(t,s)x})(f)
                     =(\mu_{t,s}^g*\mu_{t,s}^j*\delta_{U(t,s)x})(f) \\
                 =&\bigl( (\mu_{t,s}^g*\delta_{U(t,s)x})*\mu_{t,s}^j\bigr)(f)  \\
                   =&\int_\H \mu_{t,s}^g*\delta_{U(t,s)x} (dy) \int_\H f(y+z)\,\mu_{t,s}^j  (dz)\\
                 =&\bigl(\mu_{t,s}^g*\delta_{U(t,s)x}  \bigr) (p_{s,t}^j f)
                     =p_{s,t}^g (p_{s,t}^j f)(x).
\end{aligned}
\]
\end{proof}

Define for every $t\geq s$,  
\begin{equation}\label{Equ:def:Gamma_ts}
 \Gamma_{t,s}=R_{t,s}^{-1/2}U(t,s)
\end{equation}
with domain $$\mathscr{D}(\Gamma_{t,s})=\{x\in\H\colon U(t,s)x\in R_{t,s}^{1/2}(\H)\}.$$
If  $x\notin \mathscr{D}(\Gamma_{t,s})$ then we set 
$|\Gamma_{t,s}x|:=\infty$.
Let $\BB_b^+(\H)$ denote the space of all bounded positive measurable functions on $\H$.

\begin{theorem}\label{Thm1}
  For every $\alpha>1$, \(t\geq s\) and $f\in \BB_b^+(\H)$, we have 
\begin{equation}\label{HIgj}
(p_{s,t} f(x))^\alpha \leq
\exp\left(\frac{\alpha |\Gamma_{t,s}(x-y) |^2}{2(\alpha-1)}\right) p_{s,t} f^\alpha(y),
\quad
x,y\in \H.
\end{equation}
\end{theorem}

\begin{proof}
It is sufficient to consider the case \(U(t,s)(\H)\in R_{t,s}^{1/2}(\H)\), 
 since otherwise the inequality \eqref{HIg} becomes trivial because the right hand side of \eqref{HIg} is infinite by the definition of $|\Gamma_{t,s}(\cdot)|$. 
 
We claim that
we only need to show the following Harnack inequality for $p_{s,t}^g$
\begin{equation}\label{HIg}
(p_{s,t}^g f(x))^\alpha \leq
\exp\left(\frac{\alpha |\Gamma_{t,s}(x-y) |^2}{2(\alpha-1)}\right) p_{s,t}^gf^\alpha(y),
\quad
x,y\in \H.
\end{equation}
Indeed,  by Proposition \ref{Lem:Decomposition}, we have $p_{s,t}=p_{s,t}^gp_{s,t}^j$. 
If  \eqref{HIg} holds, then
by applying inequality \eqref{HIg} to $p_{s,t}^g$ and Jensen's inequality to $p_{s,t}^j$
we obtain 
   \[
    \begin{aligned}
     \bigl(p_{s,t} f(x)\bigr)^\alpha 
  = & \bigl(p_{s,t}^{g} (p_{s,t}^{j} f)(x)\bigr)^\alpha \\
  \leq & \exp\left(\frac{\alpha |\Gamma_{t,s}(x-y) |^2}{2(\alpha-1)}\right) \bigl(p_{s,t}^{g} (p_{s,t}^{j} f)^\alpha  \bigr) (y)\\
  \leq& \exp\left(\frac{\alpha |\Gamma_{t,s}(x-y) |^2}{2(\alpha-1)}\right)  \bigl(p_{s,t}^{g} (p_{s,t}^{j} f^\alpha)  \bigr) (y)\\
     = & \exp\left(\frac{\alpha |\Gamma_{t,s}(x-y) |^2}{2(\alpha-1)}\right)  \bigl(p_{s,t} f^\alpha  \bigr) (y).  
    \end{aligned}
   \]

Applying the Cameron-Martin formula 
for Gaussian measures (see \cite[Theorem 2.21]{DZ92}) 
we get
\begin{equation}\label{rho_ts}
\begin{aligned}
               & \rho_{t,s}(x-y,z)  : 	=  \frac{d N(U(t,s) (x-y), R_{t,s} )}{d N(0,R_{t,s} )}(z)\\	
                                    = & \exp\left(
				             \left\< R_{t,s}^{-1/2}U(t,s)(x-y), R_{t,s}^{-1/2}z\right\>
						 	- \frac{1}{2}|R_{t,s} ^{-1/2}U(t,s) (x-y)|^2
			             \right).
\end{aligned}			             
\end{equation}
By changing variables and using H\"older's inequality we obtain
		\begin{equation*}
			\begin{aligned}
				&p_{s,t}^g f(x)\\
			=&\int_\H f(U(t,s)x+z)\,\mu_{t,s}^g (dz)\\
			=&\int_\H f( U(t,s)y+U(t,s)(x-y)+z )\,  N(0,R_{t,s} )
												 (dz)\\
			=&\int_\H f( U(t,s)y+z' )
						\frac{d N(U(t,s) (x-y), R_{t,s} )}{d N(0,R_{t,s} )}(z')
							\,  N(0,R_{t,s} )
												 (dz')\\									 
			=&\int_\H f(U(t,s) y +z') 	\rho_{t,s}(x-y,z')
					\,\mu_{t,s}^g (dz')\\
			\leq&  \exp\left(-\frac12|\Gamma_{t,s} (x-y)|^2\right) 
			\left( \int_\H f^\alpha(U(t,s) y+z')\,\mu_{t,s}^g (dz') \right)^{1/\alpha}\cdot \\
			     &	\quad\quad
			     	 \left( \int_\H 
				 	\exp\left(\frac{\alpha}{\alpha-1}
						\< R_{t,s}^{-1/2}U(t,s)(x-y), R_{t,s}^{-1/2}z'\>
						\right) \,\mu_{t,s}^g(dz')
			       \right)^{(\alpha-1)/\alpha} \\
			=& \exp\left(\frac{1}{2(\alpha-1)}|\Gamma_{t,s} (x-y)|^2\right) 
					 \left(p_{s,t}^g f^{\alpha}(y)\right)^{1/\alpha}.
			\end{aligned}
		\end{equation*}
This proves \eqref{HIg}. Here we have used the fact that for any $h\in\H$, it holds (cf. \cite[Proposition 1.2.5, Page 11]{DZ02})
	\begin{equation}\label{Equ:Gauss_integral_exp_function}
			\int_\H \exp(\<h,x\>)\,d\mu_{t,s}^g(x) = \exp\left(\frac12 \<R_{t,s}h,h\> \right).
		\end{equation}
	\end{proof}

Applying the previous theorem, we have the following result. 

\begin{theorem}\label{Thm2} 
Fix $t\geq s$. The implications 
\((1)\Rightarrow(2)\Rightarrow(3)\Rightarrow(4)\Rightarrow(5)\)
of the following statements hold. 
\begin{enumerate}
 \item  \begin{equation}\label{Condition:NullControll2}
 		 U(t,s)(\H)\subset R_{t, s}^{1/2}(\H),
	\end{equation}
 \item $\|\Gamma_{t,s}\|<\infty$ and for every $\alpha>1$ and $f\in
          \BB_b^+(\H)$,  
    \begin{equation}\label{HI2}
       (p_{s,t}f(x))^\alpha\leq\exp\left(\frac{\alpha(\|\Gamma_{t,s}\|\cdot |x-y|)^2}{2(\alpha-1)}\right)
           p_{s,t}f^\alpha(y),\quad x,y\in \H;
    \end{equation}

\item   $\|\Gamma_{t,s}\|<\infty$ and  there  exists $\alpha>1$ such that \eqref{HI2}
holds for all $f\in \BB_b^+(\H);$

\item $\|\Gamma_{t,s}\|<\infty$
and for every $f\in \BB_b^+(\H)$ with $f> 1$, 
	\begin{equation}\label{LHI}		
		p_{s,t}\log f(x)\le \log p_{s,t} f(y) +\frac{\|\Gamma_{t,s}\|^2}2 |x-y|^2,\quad
	x,y\in \H;
\end{equation}
\item $p_{s,t}$ is strong Feller. 
\end{enumerate}
In particular, if  $\m_{t,s}\equiv 0$, then these statements are equivalent to each other. 
\end{theorem}

\begin{proof}
If \eqref{Condition:NullControll2} hold, then $\|\Gamma_{t,s}\|$ is bounded. Hence
by Theorem \ref{Thm1}, we get (2) from (1). 
That (2) implies (3) is trivial. 
The implications (3)$\Rightarrow$(4)$\Rightarrow$(5) are consequences of Harnack inequalities,  as proved in \cite{Wang09C}. 

It remains to show that (5) implies (1) in the case $\m_{t,s}\equiv 0$.
Note that $$p_{s,t}f(x) =\int_\H f(y)\, N(U(t,s)x, R_{t,s})(dy).$$
If \eqref{Condition:NullControll2} doesn't hold, then there exists
$x_0\in \H$ such that $U(t,s) x_0 \notin R_{t,s}^{1/2}(\H).$
Take $x_n=\frac{1}{n}x_0\in\H$, $n=1,2,\cdots$. 
By the Cameron-Martin theorem (see e.g. \cite{DZ92}),  we know that 
for each $n=1,2,\cdots$, the Gaussian measure $\mu_n:=N(U(t,s)x_n, R_{t,s})$ is orthogonal to $\mu_0:=N(0, R_{t,s})$  since
$U(t,s)x_n\notin R_{t,s}^{1/2}(\H).$  
That is, for all $n=1,2,\cdots,$ there exists $A_n\in \B(\H)$ such that
$\mu_n(A_n)=1$, $\mu_0(A_n)=0$. 
Set $A:=\cup_{n\geq 1} A_n$. Then $\mu_0(A)=0$,
$\mu_n(A)=1$ since $\mu_0(A)\leq \sum_{n=1}^\infty\mu_0(A_n)=0$ and $\mu_n(A)\geq \mu_n(A_n)=1$. 

Take $f=\Eins_A$. 
Since $x_n$ tends to $0$ as $n\to\infty$ and  $p_{s,t}$ is strong Feller, 
\(p_{s,t}f(x_n)\) should converge to \(p_{s,t}f(0) \) as $n\to\infty$. But this is impossible
because we have 
$p_{s,t}f(x_n)=1$, $p_{s,t}f(0)=0$. 
Therefore we  have \eqref{Condition:NullControll2}.
\end{proof}

\begin{remark}
If $R_{t,s}$ has the form \eqref{Pi_ts}, then \eqref{Condition:NullControll2} is equivalent to the null controllability of 
 a non-autonomous control system \eqref{NonAutoControlSys} (see Section \ref{Sec:NullControll} for details). For this reason, condition  \eqref{Condition:NullControll2} is also called null-controllability condition.  
This gives an equivalent description of the strong Feller property. 
\end{remark}

\begin{remark} In \cite{DaP95} the fact that the null controllability implies the strong Feller property was proved 
 for autonomous Ornstein-Uhlenbeck processes driven by a Wiener process and with deterministic perturbation.  
 Our result generalizes his result. 
\end{remark}

In fact \eqref{Condition:NullControll2} implies more. Let  $UC^\infty(\H)$ denote the space of all 
infinitely Fr\'echet differentiable functions with uniform continuous derivatives  on $\H$.
\begin{prop}
	Suppose  that \eqref{Condition:NullControll2} holds. Then for every $f\in\BB_b(\H)$ and every $t>s$,  
	$p_{s,t}f\in UC^\infty(\H)$. 
\end{prop}
\begin{proof}
	In view of  the decomposition $p_{s,t}=p_{s,t}^gp_{s,t}^j$ shown in Proposition \ref{Lem:Decomposition},  we only need to show that $p_{s,t}^g\in UC^\infty(\H)$ for every $g\in\BB_b(\H)$.  The rest of the proof is the same as in \cite[Theorem 6.2.2]{DZ02}.
\end{proof}

We have the following quantitative estimate for the strong Feller property. This result is shown in 
 \cite{ORW09} for L\'evy driven Ornstein-Uhlenbeck process by a coupling method. 

\begin{prop}\label{Prop:StrongFellerEstimate}
Let $t>s$ and $x,y\in\H$. 
Then 
\begin{equation}\label{StrongFellerEstimate}  
\begin{aligned}
       & |p_{s,t} f(x)-p_{s,t} f(y)|^2\\
\leq & \left(\e^{|\Gamma_{t,s}(x-y)|^2}-1\right) 
\min \left\{p_{s,t}f^2(z)-(p_{s,t}f(z))^2\colon z=x,y \right\}.
\end{aligned}
\end{equation}
\end{prop}

\begin{proof}
 Let $h=p_{s,t}^j f$. Then by Proposition \ref{Lem:Decomposition} we have $p_{s,t}f=p_{s,t}^gh$.
 Moreover, 
 $$h^2=(p_{s,t}^j f)^2\leq p_{s,t}^j f^2$$
 by Jensen's inequality. 
  So, for every $z\in\H$, we have
 \begin{equation}\label{StrongFellerEstimate:f2h}
\begin{aligned}
		          & p_{s,t}^g h^2 (z) -\bigl(p_{s,t}^gh(z)\bigr)^2 \\
		    \leq  &   p_{s,t}^g p_{s,t}^jf ^2 (z) -\bigl(p_{s,t}^g p_{s,t}^j f(z)\bigr)^2      
		    =p_{s,t}f^2(z)-(p_{s,t}f(z))^2. 
\end{aligned}
\end{equation}
Note also that $x,y$ play the same role in \eqref{StrongFellerEstimate}. So, according to \eqref{StrongFellerEstimate:f2h} we only need to show the following inequality
\begin{equation}\label{StrongFellerEstimate:Gauss}  
\begin{aligned}
        |p_{s,t}^g h(x)-p_{s,t}^g h(y)|^2
\leq  \left(\e^{|\Gamma_{t,s}(x-y)|^2}-1\right) 
\left( p_{s,t}^g h^2(y)-(p_{s,t}^g h(y))^2 \right).
\end{aligned}
\end{equation}

According to formula \eqref{rho_ts} for $\rho_{t,s}(x-y,z)$, we have
\[
      p_{s,t}^g h(x)=\int_\H h(U(t,s)x+z)\, \mu_{t,s}^g(dz)\\
  = \int_\H \rho_{t,s}(x-y,z) h(U(t,s)y+z)\, \mu_{t,s}^g(dz).
\] 
Therefore we have 
\[
\begin{aligned}
    &  |p_{s,t}^g h(x)-p_{s,t}^g h(y)|^2 \\
 =& \left(
     \int_\H [\rho_{t,s}(x-y,z)-1]\cdot [h(U(t,s)y+z)-p_{s,t}^g h(y)]\, \mu_{t,s}^g(dz)
  \right)^2\\
 \leq & \int_\H (\rho_{t,s}(x-y,z)-1)^2\,\mu_{t,s}^g(dz)\int_\H 
\bigl[ h(U(t,s)y+z)-p_{s,t}^g h(y)\bigr]^2 \, \mu_{t,s}^g(dz)\\
=& \left( \int_\H  \rho^2_{t,s}(x-y,z)\, \mu_{t,s}^g(dz) -1 \right)  
 	      \cdot  \left(
  \int_\H h^2(U(t,s)y+z)  \, \mu_{t,s}^g(dz) - (p_{s,t}^g h(y))^2 \right) 	\\
=&	\left(\e^{|\Gamma_{t,s}(x-y)|^2}-1\right) 
\left( p_{s,t}^g h^2(y)-(p_{s,t}^g h(y))^2 \right).        
\end{aligned}
\]
Note that here we have used \eqref{Equ:Gauss_integral_exp_function} again to obtain 
\[
	\int_\H  \rho^2_{t,s}(x-y,z)\, \mu_{t,s}^g(dz) =  \e^{|\Gamma_{t,s}(x-y)|^2}.
\]
\end{proof}

Now we apply the Harnack inequality \eqref{HIgj} to study the hyperboundedness of the transition function $p_{s,t}$. In \cite{GL08} hypercontractivity is studied for the Gaussian case via log-Soboblev inequality. 

\begin{theorem}\label{Thm:Hypercontractive}
 Let $(\nu_t)_{t\in\R}$ be an evolution system of measures for $p_{s,t}$.  For every $s\leq t$, $\alpha>1$, and $\varepsilon>0$, let
\[C_{s,t}(\alpha,\varepsilon):=\int_\H \left[ \int_\H
		\exp\left(        -\frac{\alpha |\Gamma_{t,s}(x-y) |^2}{2(\alpha-1)}
			 \right) \nu_s(dy)
	\right]^{-(1+\varepsilon)}\,\nu_s(dx).\]
Then 
\begin{equation}\label{Equ:Hyper}
\| p_{s,t} f \|_{L^{\alpha(1+\varepsilon)}(\H,\nu_s)}
\leq C_{s,t}(\alpha,\varepsilon)^{-\alpha(1+\varepsilon)} \|f\|_{L^\alpha(\H,\nu_t)}.
\end{equation}
\end{theorem}
\begin{proof}
  From the Harnack inequality \eqref{HIgj} we have
\[
(p_{s,t} f(x))^\alpha 
\exp\left[-\frac{\alpha |\Gamma_{t,s}(x-y) |^2}{2(\alpha-1)}\right]
\leq  p_{s,t}f^\alpha(y),
\quad
x,y\in \H.
\]
Integrating both sides of the inequality above with respect to $\nu_s(dy)$ and using the fact that $(\nu_t)_{t\in\R}$ is an evolution system of measures, we get
\[
(p_{s,t}|f|)^{\alpha}(x) 
   \int_\H	\exp\left( -\frac{\alpha |\Gamma_{t,s}(x-y) |^2}{2(\alpha-1)}
			 \right) \,\nu_s(dy) \leq 
    \int_\H |f|^\alpha\,\nu_t(dy).
\]
Hence 
\[
(p_{s,t}|f|)^{\alpha(1+\varepsilon)}(x) \leq 
\left[ \int_\H
		\exp\left(        -\frac{\alpha |\Gamma_{t,s}(x-y) |^2}{2(\alpha-1)}
			 \right) \nu_s(dy)
	\right]^{-(1+\varepsilon)}
\|f\|^{\alpha(1+\varepsilon)}_{L^{\alpha}(\H,\nu_t)}.
\]
Integrating both sides of the inequality above with respect to $\nu_s(dx)$, we get \eqref{Equ:Hyper}.
\end{proof}


\section{Semi-linear equations}\label{Sec:Semi-linear}
Fix $s\in\R$ and consider the following equation for $t\geq s$, 
\begin{equation}\label{Equ:Semi-linear} 
\left\{
\begin{aligned}
d X(t,s,x)  &= A(t)X(t,s,x)\,d t + F(t, X(t,s,x)) d t + R^{1/2} d W_t, \\
   X(s,s,x) &= x\in\H,
\end{aligned}
\right.
\end{equation} 
where 
\begin{enumerate}
\item $(A(t))_{t\in\R}$ is a family of operators on $\H$ associated with an evolution family $(U(t,s)_{t\geq s})$ (See Section \ref{Sec:Intro});
\item $R$ is a trace class operator on $\H$;
\item $(W_t)_{t\in\R}$ is an \H-valued  cylindrical Wiener process on some filtered probability space $(\Omega, (\F_t)_{t\in\R},\F, \P)$;
\item $F$ is a measurable map from $[s,+\infty)\times \H$ to $R^{1/2}(\H)$ satisfying 
		\begin{equation}\label{linear:growth}
		   |R^{-1/2}F(t,x)|^2\leq k_1+k_2|x|^2,\quad t\in\R,\ x\in\H
		\end{equation}
	for  some constants $k_1,k_2>0$.
\end{enumerate}

\begin{prop}\label{prop:semi-linear:existence}
Equation \eqref{Equ:Semi-linear} martingale solution. 
\end{prop}

\begin{proof}
  For every $r\in[s,t]$, 
  let 
  	\[ W_U(r,s) : = \int_s^{r} U(r,\sigma) R^{1/2}\,d W_\sigma \]
and  set $$\tilde X(r,s,x):=U(r,s)x+W_U(r,s).$$
Moreover,  for every $r\in[s,t]$, $[s',t']\subset [s,t]$, define
\[
\begin{aligned} 
    &\psi_x(r,s):=R^{-1/2}F(r, \tilde{X}(r,s,x))=R^{-1/2}F(r, U(r,s)x+W_U(r,s)), \\
    &\tilde W_r^x=W_r-\int_s^r \psi_x(\sigma,s)\, d \sigma, \\
    &M_{t',s'}^x=\exp\left(\int_{s'}^{t'} \<\psi_x(\sigma,s),  d W_\sigma\>
    	-\frac{1}{2} \int_{s'}^{t'} |\psi_x(\sigma,s)|^2\, d \sigma\right).
\end{aligned}
\]
We first show that $\E M_{t,s}^x=1$. 
By \eqref{linear:growth}, for every $r\in[s,t]$, 
\[|\psi_x(r,s)|^2\leq k_1+2k_2(|U(r,s)x|^2+|W_U(r,s)|^2).\]  
Hence we get 
\[
\begin{aligned}
	&\E\exp\left(\frac{1}{2} \int_{s}^{t} |\psi_x(\sigma,s)|^2\, d \sigma\right)\\
\leq   &\E\exp\left(\frac{k_1\vee k_2}{2} \int_{s}^{t} (1+2|U(\sigma,s)x|^2)\, d \sigma\right)
	\cdot
		\E\exp\left(\frac{1}{2} \int_{s}^{t} |W_U(\sigma,s)|^2\, d \sigma\right).
\end{aligned}
\] 
Since  $\int_{s}^{t} |W_U(\sigma,s)|^2\, d \sigma$ is Gaussian distributed, applying Fernique's Theorem, 
for a fine partition $s=t_0<t_1<\cdots<t_{n-1}<t_n=t$, we have
\[\E\exp\left( \frac{1}{2} \int_{t_{i-1}}^{t_i} |\psi_x(\sigma,s)|^2\, d \sigma \right)<+\infty.\]
This implies that for each $i=1,2,\cdots, n$, $M_{t_{i-1},t}^x$ for $t\in[t_{i-1},t_i]$,  is a martingale. Noting that
\(M_{t,s}^x=M_{t_{n},t_{n-1}}^x\cdots M_{t_1,t_0}^x\),
we get $\E M_{t,s}^x=1$.   
  
Consequently, we can define a new probability measure $\Q_x:= M_{t,s}^x\P$ on $(\Omega, \F_t)$.
By \cite[Theorem 10.14]{DZ92}, $\tilde W_r^x$ is also a Wiener process with respect to $\Q_x$. 
Hence 
\[
\begin{aligned}
\tilde{X}(t,s,x)&=U(t,s)x+\int_s^t U(t,r)R^{1/2}\,dW_r \\
&=U(t,s)x +\int_s^t U(t,r)F(r, \tilde{X}(r,s,x))\,dr+\int_s^t U(t,r)R^{1/2}\,d\tilde{W}_r. 
\end{aligned}
\]
This shows that $\tilde{X}(t,s,x)$ is a martingale solution of \eqref{Equ:Semi-linear} on $(\Omega, (\F_t)_{t\geq s}, \F, \Q_x)$.
\end{proof}

We shall need the following fact.

\begin{lemma}\label{lemma:exp_estimate}
Let $s\in\R$. Set
 $$\lambda:=\operatorname{tr} \int_s^{s+1} U(s+1,\sigma)RU(s+1,\sigma)^*\,d \sigma.$$ Then  
$$C_0:=\sup_{r\in [s,s+1]}\E \exp\left(|W_U(r,s)|^2/(4\lambda)\right)<\infty$$ and   
for every $\kappa>0$ and $t\in [s,s+(1\wedge (4\lambda\kappa)^{-1})]$, 
\begin{equation}\label{GGG} 
\E \exp\left(\kappa \int_s^t |W_U(r, s) |^2\,d r\right) <C_0^{4\lambda\kappa (t-s)}.
\end{equation}
\end{lemma}
\begin{proof}
Note that the covariance operator of $$W_U(r,s)=\int_s^r U(r,\sigma)R^{1/2}\,dW_\sigma$$ is given by
$$\int_s^{r} U(r,\sigma)RU(r,\sigma)^*\,d \sigma.$$
By Fernique's Theorem (see \cite[Propostion 2.16]{DZ92}), it follows that $C_0<\infty$.
Moreover, 
\[
\begin{aligned}  
     &\E \exp\left(\kappa \int_s^t |W_U(r, s) |^2\,d r\right) 
           =\E\exp\left(\frac 1 {t-s} \int_s^t \kappa (t-s) |W_U(r,s) |\,d r\right) \\
    \le& \frac 1 {t-s} \int_s^t \E \exp\left(\kappa (t-s) |W_U(r,s)|^2\right)\,d r\\
   \le & \frac 1 {t-s} \int_s^t \left[\E\exp\left(|W_U(r,s)|^2/(4\lambda)\right)\right]^{4\lambda
\kappa (t-s)}\,d r
\le  C_0^{4\lambda\kappa (t-s)}.
\end{aligned}
\]
\end{proof}

From Lemma \ref{lemma:exp_estimate} we get that for every $p>0$, there exists $t_p>0$ such that for every $t\in [s,s+t_p]$, 
\[C_{p,k_2}(t,s):=\E\exp\left(
				2p(2p+1)k_2\int_s^t|W_U(r,s)|^2\,ds
												\right)<\infty.
\]
In particular, if $k_2=0$ then
$C_{p,0}(t,s)=1$ for all $ t\ge s.$ 

\begin{lemma} \label{L4.2} For any $t>s$, $p>1$, $\delta>0$
and $x\in\H$, 
\begin{equation*} \begin{split} &\E (M_{t,s}^x)^p \le
(C_{p,k_2}(t,s))^{1/2}
\exp\Big( \frac { p(2p-1)} 2\int_s^t(k_1+2k_2 |U(r,s)x |^2)\,d r\Big),\\
&\E(M_{t,s}^x)^{-\delta}\le  (C_{\delta,k_2}(t,s))^{1/2} \exp\Big( \frac {
\delta(2\delta+1)} 2 \int_s^t(k_1+2k_2 |U(r,s)x |^2) \, d r\Big).\end{split}\end{equation*}
\end{lemma}

\begin{proof} From the proof of Proposition \ref{prop:semi-linear:existence} we obtain that for every $\kappa\in\R$
\[
t\mapsto \exp\left(\kappa \int_s^t \<\psi_x(r,s), dW_r\>- \frac {\kappa^2} 2 \int_s^t  |\psi_x(r,s) |^2\,d r\right)\]
is  a martingale. Therefore, 
\[
    \begin{aligned}
           \E(M_{t,s}^x)^p
        =&\E \exp\left(p  \int_s^t\< \psi_x(r,s),dW_r\>- p^2 \int_s^t  | \psi_x(r,s) |^2\,d r\right) \\
            &    \quad\quad\quad\quad\quad\quad\quad\quad\
            \cdot
                  \exp\left(  \frac{p(2p-1)}{2} \int_s^t  | \psi_x(r,s) |^2\,d r \right)\\
        \leq & \left[ \E \exp\left(2 p  \int_s^t\< \psi_x(r,s),dW_r\>
                -\frac12 (2 p)^2 \int_s^t  | \psi_x(r,s) |^2\,d s\right) \right]^{1/2}\\
            & \quad\quad\quad\quad\quad\quad\quad\quad\
            \cdot \left[
                \E \exp\left(   p(2 p-1)  \int_s^t  | \psi_x(r,s) |^2\,d s \right)
            \right]^{1/2}\\
            = & \left[\E \exp\left(   p(2 p-1)  \int_s^t  | \psi_x(r,s)|^2\,d s \right)
            \right]^{1/2}.
    \end{aligned}
\]
This implies the first inequality claimed in the lemma, since by \eqref{linear:growth}
$$ | \psi_x(r,s) |^2\le k_1 + 2k_2 |W_U(r,s) |^2 +2k_2 |U(r,s)x |^2.$$
Similarly, the second inequality follows by
\[
    \begin{aligned}
           \E(M_{t,s}^x)^{-\delta}
        =&\E \exp\left(-\delta  \int_s^t\< \psi_x(r,s),dW_r\>- \delta^2 \int_s^t  | \psi_x(r,s) |^2\,d r\right)\\
             & \quad\quad\quad\quad\quad\cdot  \exp\left(  \frac{\delta(2\delta+1)}{2} \int_s^t  | \psi_x(r,s) |^2\,d r \right)\\
        \leq & \left[ \E \exp\left(-2 \delta  \int_s^t\< \psi_x(r,s),dW_r\>
                -\frac12 (2 \delta)^2 \int_s^t  | \psi_x(r,s) |^2\,d s\right) \right]^{1/2}\\
            & \quad\quad\quad\quad\quad
            \cdot \left[
                \E \exp\left(   \delta(2 \delta+1)  \int_s^t  | \psi_x(r,s) |^2\,dr \right)
            \right]^{1/2}\\
            = &  \left[
                \E \exp\left(   \delta(2 \delta+1)  \int_s^t  | \psi_x(r,s)|^2 \,dr \right)
            \right]^{1/2}.
    \end{aligned}
\]
\end{proof}

By the proof of Proposition \ref{prop:semi-linear:existence}, we see that 
$\tilde{X}(t,s,x)$ is a solution of \eqref{Equ:Semi-linear}. Hence we define the ``transition semigroup" of $X(t,s,x)$ by
\begin{equation}\label{PPP}
P_{s,t}^Ff(x) = \E_{\Q_x} f(\tilde X(t,s,x)),\quad f\in\BB_b(\H).
\end{equation} 

We have the following result. 
\begin{theorem} For any $t>0$,  $\alpha>1$, $x,y\in \H$, $p, q>1$ with $\alpha/(pq)>1$,
and  $f\in \BB_b^+(\H)$
\begin{equation}\label{HI:Semi-linear}
	(P_{s,t}^F f)^\alpha(x)\leq  N P_{s,t}^F f^\alpha (y).
\end{equation} 
Here we set $\Gamma_{t,s}^F:=R^{-1/2}U(t,s)$ and 
\[
\begin{aligned}
  N&:= \left(C_{\frac{p}{p-1},k_2}(t,s)\right)^{\alpha p/(2(p-1))} \cdot
                         \left(C_{\frac{1}{q-1}, k_2}(t,s)\right)^{\alpha q/(2(q-1))}     \cdot
          \exp\left( \frac{\alpha q  | \Gamma_{t,s}^F(x-y) |^2}{2(\alpha-q)} \right. \\
         &\quad\quad\quad \left.+ \alpha \left[\frac{p+1}{p-1}+\frac{q+1}{q(q-1)}\right]\int_s^t
                         \left[k_1+k_2( |U(r,s) x |^2+ |U(r,s)
                         y |^2)\right]\,d r\right).
\end{aligned}
\] 
Assume that for every $s\leq r\leq t$, $P^F_{s,t}=P^F_{s,r}P^F_{r,t}$. If $ \|\Gamma_{t,s}^F\|<\infty$ for every $t\geq s$, then
$P_{s,t}^F$  is strong Feller.
\end{theorem}

 \begin{proof} 
  Recall that $\tilde X(t,s,x)$ is a mild solution to 
$$\,d \tilde X(t,s,x) = A(t) \tilde X(t,s,x) d t + R^{1/2} d W_t,\quad \tilde{X}(s,s,x)=x.$$ 
 Let $P_{s,t}^0$ be the semigroup of $\tilde X(t,s,x)$ under $\P$. Then by Theorem \ref{Thm1} we have
\begin{equation}\label{HI:Gauss}
    (P_{s,t}^0f)^\alpha(x)\leq P_{s,t}^0f^\alpha (y)
                        \exp\left( \frac{\alpha | \Gamma_{t,s}^F(x-y) |^2}{2(\alpha-1)} \right)
    , \quad f\in \BB_b^+(\H),
\end{equation}

     For simplicity, we set
        $p':=\frac{p}{p-1}$,         
        $q':=\frac{q}{q-1}$,     
        $\theta=\alpha/(pq)$. 
     By \eqref{HI:Gauss} we have
\[
\begin{aligned}
    P_{s,t}^F f (x) &= \E_{\Q_x} f(\tilde X(t,s,x))=\E M_{t,s}^x f(\tilde X(t,s,x)) \\
             &   \leq  (\E f^p (\tilde X(t,s,x)))^{1/p} (\E (M_{t,s}^x)^{p'} )^{1/p'}
                =  (P_{s,t}^0  f^p (x))^{1/p} (\E(M_{t,s}^x)^{p'} )^{1/p'}\\
             &   \leq \left[    P_{s,t}^0 f^{\theta p}(y)
                    \exp\left(  \frac{\theta |\Gamma_{t,s}^F(x-y)|^2}{2(\theta-1)}   \right)
                        \right]^{1/(\theta p)} (\E(M_{t,s}^x)^{p'} )^{1/p'}.
\end{aligned}
\]
On the other hand, for every $g\in \BB_b^+(\H)$,
\[
\begin{aligned}
    P_{s,t}^0g(y)&\leq \E_\P g(\tilde X(t,s,y)) =\E_{\Q_y}g(\tilde X(t,s,y)) (M_{t,s}^y)^{-1}\\
    &\leq (P_{s,t}^F g^q (y))^{1/q} (\E(M_{t,s}^y )^{1-q'} )^{1/q'}.
\end{aligned}    
\]
So, taking $g=f^{\theta p}$ we obtain
\[
    ( P_{s,t}^F f )^\alpha (x)
        \leq P_{s,t}^F f^{\alpha}(y) \exp\left( \frac{\alpha |\Gamma_{t,s}^F(x-y)|^2}{2p(\theta-1)}   \right)
        (\E(M_{t,s}^x)^{p'} )^{\alpha/p'}
         (\E (M_{t,s}^y)^{1-q'} )^{\alpha/q'}.
\]
This implies the desired Harnack inequality according to Lemma \ref{L4.2}.

		Now we show that $P_{s,t}^F$ is strongly Feller. 
		Let $f\in \BB_b^+(\H).$ By 
		\eqref{GGG} and \eqref{HI:Semi-linear}, for any $\alpha>1$ there exist constants
		$t_\alpha, c_\alpha>0$ and a positive function $H_\alpha(r,s)$,  $r\in (s,s+t_\alpha)$ such
		that
		\begin{equation}\label{LLL}
		P_{s,r}^Ff(x)\le (P_{s,r}^Ff^\alpha(y))^{1/\alpha} \e^{c_\alpha (r-s)+
		| x-y| ^2H_\alpha(r,s)},\quad  r\in (s,s+t_\alpha).
		\end{equation} 
		We take $t_\alpha<t-s$. 
		Then, using the assumption that $P_{s,t}^F$ is a semigroup,  for every $r\in(s,s+t_\alpha)$, we get
		\begin{equation}\label{equ:limsup}
		\begin{aligned}
		     & \varlimsup_{x\to y} P_{s,t}^Ff(x)=\varlimsup_{x\to y} P_{s,r}^FP_{r,t}^Ff(x)\\
		 \le&
			\varlimsup_{\alpha\to 1}\varlimsup_{r\to s}\varlimsup_{x\to y}\big[
			P_{s,r}^F(P_{r,t}^Ff)^\alpha(y)\big]^{1/\alpha} \e^{c_\alpha (r-s)
			+| x-y| ^2H_\alpha(r,s)}\\
		\le& \varlimsup_{\alpha\to 1}\varlimsup_{r\to s}\varlimsup_{x\to y}\big[
			P_{s,t}^F f^\alpha(y)\big]^{1/\alpha} \e^{c_\alpha (r-s) +| x-y| ^2H_\alpha(r,s)}= P_{s,t}^F f
			(y).
		\end{aligned}\end{equation}
		On the other hand, (\ref{LLL}) also implies for every $r\in (s,s+t_\alpha)$
		\begin{equation*}
		\begin{aligned}
			P_{s,t}^Ff(x)&\ge 
		   		\big[P_{s,r}^F(P_{r,t}^Ff)^{1/\alpha}(y)\big]^\alpha\e^{-\alpha c_\alpha (r-s)-\alpha
					H_\alpha(r,s)| x-y| ^2}\\
			&\ge \big[P_{s,t}^Ff^{1/\alpha}(y)\big]^\alpha\e^{-\alpha c_\alpha (r-s)-\alpha
			H_\alpha(r,s)| x-y| ^2}.
		\end{aligned}\end{equation*}
			So, first letting $x\to y$ then $r\to s$ and finally $\alpha\to 1$, 
			we arrive at
			\begin{equation}\label{equ:liminf}
			\varliminf_{x\to y} P_{s,t}^Ff(x)\ge P_{s,t}^Ff(y).
			\end{equation}
	 From \eqref{equ:limsup} and \eqref{equ:liminf} we see $P_{s,t}^Ff$ is
			continuous. So, $P_{s,t}^F$ is strongly Feller. 
\end{proof}


\section{Appendix: null controllability}\label{Sec:NullControll}
Consider the following non-autonomous linear control system 
\begin{equation}\label{NonAutoControlSys} 
\left\{
    \begin{aligned}
	dz(t)&=A(t)z(t)dt+C(t) u(t)\,dt,\\
    	 z(s)&=x,
    \end{aligned} 
  \right.
\end{equation}
where \((A(t))_{t\in\R}\) is a family 
of linear operators on \(\H\) with dense domains 
 and \((C(t))_{t\in\R}\) is a family of bounded linear operators on \(\H\).
Let \((U(t,s))_{t\geq s}\) be an evolution family on \(\H\)
associated with \((A(t))_{t\in\R}\). Consider the mild solution of \eqref{NonAutoControlSys}
\[
  z(t,s,x)=U(t,s)x+\int_s^t U(t,r)C(r)u(r)\,dr.\quad x\in \H, \ t\geq s.
\]
 $z(t,s,x)$ is interpreted as the state of the system and $u$ as a strategy to control the system. 
If there exists $u\in L^2([s,t],\H)$ such that \(z(t,s,x)=0\), then 
we say the system \eqref{NonAutoControlSys} can be transferred to 0 at time $t$ from initial state $x\in\H$ at time $s$.
If for every initial state $x\in\H$ the system \eqref{NonAutoControlSys} can be transferred to 0 then we say the system \eqref{NonAutoControlSys} is null controllable at time $t$. We refer to 
\cite{Zab08} (see also \cite[Appendix B]{DZ92}) for details on the null controllability of autonomous control systems. 

Set for every $t\geq s$
\begin{equation}\label{Pi_ts}
\Pi_{t,s}x:=\int_s^t U(t,r)C(r)C(r)^*U(t,r)^*\,dr,\quad  x\in\H.
\end{equation}

\begin{prop}\label{Thm:NullControl}
   Let \(x\in\H\) and $t\geq s$. The system \eqref{NonAutoControlSys} can be transferred to 0 at time $t$ from $x$ if and only if \(U(t,s)x\in \Pi_{t,s}^{1/2}(\H) \).
Moreover, the minimal energy among all strategies transferring $x$ to 0 at time $t$ is given by \(|\Pi_{t,s}^{-1/2}U(t,s)x|^2\), i.e. 
\begin{equation}\label{MinimalEnergy}
\begin{aligned}
 & |\Pi_{t,s}^{-1/2}U(t,s)x|^2\\
=& \inf\left\{\int_s^t |u(r)|^2\,dr\colon z(t,s,x)=0, z(s,s,x)=x, u\in L^2([s, t],\H) \right\}.
\end{aligned}
\end{equation}
\end{prop}

\begin{proof} For every $t\geq s$ define a linear operator 
\[
L_{t,s}\colon L^2([s,t],\H)\to \H,\quad 
u\mapsto 
L_{t,s}u:=\int_s^t U(t,r)C(r)u(r)\,dr.
\]
The adjoint $L^*_{t,s}$ of $L_{t,s}$ is given by
\[
(L^*_{t,s}x)(r)=C^*(r)U(t,r)^*x,\quad x\in\H, \ r\in [s,t].
\]
It is easy to check  that  \[\Pi_{t,s}=L_{t,s}L^*_{t,s}.\]
Then by \cite[Corollary B.4]{DZ92}, we know that 
\(L_{t,s}(L^2([s,t],\H)=\Pi_{t,s}(\H). \)
Hence the first assertion of the theorem is proved since the initial state $x$ can be transferred to $0$ if and only if 
\(U(t,s)x\) is contained in the image space of \(L_{t,s}\) due to the fact that 
 $z(t,s,x)=U(t,s)x+L_{t,s}u$.

By \cite[Corollary B.4]{DZ92}  we also get 
\begin{equation}\label{Equ:Proof:Inverse=}
 |\Pi_{t,s}^{-1/2}y|=|L_{t,s}^{-1}y|,\quad y\in L_{t,s}(L^2([s,t],\H)).
\end{equation}
Here the inverse is understood as a pseudo--inverse. 
Taking  
$y=U(t,s)x$ in \eqref{Equ:Proof:Inverse=}, we obtain \eqref{MinimalEnergy}.
\end{proof}

From Proposition \ref{Thm:NullControl}, we get the following corollary.

\begin{cora}\label{Cora:NullControl}
   The system \eqref{NonAutoControlSys} is null controllable at time $t$ if and only if 
\begin{equation}\label{NullControllabilityCond}
   U(t,s)(\H)\subset \Pi_{t,s}^{1/2}(\H). 
\end{equation}
\end{cora}

From \eqref{MinimalEnergy}, it is easy to get upper bounds of \(|\Pi_{t,s}^{-1/2}U(t,s)x|^2\) by choosing proper null control functions $u$.
The following proposition is analogous to \cite[Proposition 2.1]{ORW09}.

\begin{prop}\label{Prop:EnergyEstimate} Let $t>s$. Assume that 
for every $r\in[s,t]$, the operator $C(r)$ is invertible. 
Then for every strictly positive function
$\xi\in C([s,t])$,
\begin{equation}\label{2.1}
       |\Pi_{t}^{-1/2} U(t,s) x|^2 
  \leq
      \frac{ \int_{s}^{t}  |C(r)^{-1}U(r,s) x|^{2}\,\xi_{r}^{2}    \,d r }
        { \left( \int_{s}^{t}   \xi_{r} \, d r \right) ^2},
   \quad x\in\H. 
\end{equation}
Especially if $C(r)\equiv C$ and 
\(|C^{-1}U(r,s)x|^2 \leq h(r) |C^{-1}x|^2\) for every \(x\in\H\), 
then 
\begin{equation}\label{SpecialEnergyEstimate}
   |\Pi_{t}^{-1/2} U(t,s) x|^2  \leq
            \frac{|C^{-1}x|^2}{ \int_{s}^{t}  h(r)^{-1}\, d r  },\quad x\in \H.
\end{equation}
\end{prop}

\begin{proof}
We only need to consider the case where $U(t,s) x\in \Pi_{t,s}^{1/2}(\H)$ and the function 
$[s,t]\ni r\mapsto \xi_r C(r)^{-1} U(r,s) x$ belongs to
$L^2([0,t],\H)$. Then the following function 
\begin{equation*}
u(r):=-\frac{\xi_{r}}{\int_{s}^{t} \xi_{r}\, d r } C(r)^{-1} U(r,s) x
,\quad r\in[s,t],
\end{equation*}
is a null control of the system \eqref{NonAutoControlSys}. And hence the estimate \eqref{2.1} follows from 
\eqref{MinimalEnergy}.
The second estimate \eqref{SpecialEnergyEstimate} 
follows by taking $\xi(r)=h(r)^{-1}$ for all $r\in[s,t]$.
\end{proof}

\textbf{Acknowledgment} The authors thank Alexander Grigor$'$yan for discussing  which press the first author to find a simpler  proof of the Harnack inequality for the Ornstein-Uhlenbeck processes with L\'evy noise.
This leads to the new proof of the Harnack inequality presented in this paper. 
The first author also thank Zenghu Li for useful discussion. 
 
\bibliographystyle
{amsalpha}
 
\addcontentsline{toc}{section}{References} 

\def\cprime{$'$}
\providecommand{\bysame}{\leavevmode\hbox to3em{\hrulefill}\thinspace}
\providecommand{\MR}{\relax\ifhmode\unskip\space\fi MR }
\providecommand{\MRhref}[2]{%
  \href{http://www.ams.org/mathscinet-getitem?mr=#1}{#2}
}
\providecommand{\href}[2]{#2}

\end{document}